\documentclass[12pt,a4paper]{amsart}
\newif\ifpersonal
\personaltrue 
\usepackage{graphicx}
\usepackage{epstopdf}
\usepackage{amssymb}
\usepackage{amsthm}
\usepackage{amsmath}
\usepackage{pstcol, pst-node}
\usepackage{mathrsfs}
\usepackage{amscd}
\usepackage{enumerate,comment}
\usepackage[utf8]{inputenc} 

\usepackage[all]{xy}
\newtheorem{thm}{Theorem}

\newtheorem{cor}[thm]{Corollary}
\newtheorem{prop}[thm]{Proposition}
\newtheorem{lem}[thm]{Lemma}
\newtheorem{rem}[thm]{Remark}
\newtheorem{claim}[thm]{Claim}

\theoremstyle{definition}
\newtheorem{defn}[thm]{Definition}

\newtheorem{prop-def}[thm]{Proposition-Definition}
\newtheorem{example}[thm]{Example}
\newtheorem{rem-defn}[thm]{Remark-Definition}

\theoremstyle{remark}
\begin{document}
\title[Correspondence theorems]{Correspondence
 theorems for tropical curves I}
\author{Takeo Nishinou}
\thanks{email : nishinou@rikkyo.ac.jp}
\address{Department of Mathematics, Rikkyo University,
  Nishi-Ikebukuro, Toshima, Tokyo, Japan} 
\begin{abstract}
In this paper, we study the deformation theory of degenerate algebraic curves
 on singular varieties which appear as the degenerate limit of families of varieties.
For this purpose, we systematically develop a new method to calculate 
 the obstruction cohomology class of degenerate algebraic curves.
This enables us to judge whether a given degenerate curve can be 
 deformed to a smooth curve or not in variety of situations.
In this paper, we apply it to curves of genus one on degeneration of toric varieties.
In particular, we obtain the necessary and sufficient condition for the
 realizability of tropical curves of genus one, extending various results
 obtained so far \cite{BPR, CFPU, Katz, N1, S, Tyo}.   
\end{abstract}
\maketitle
\section{Introduction}
In this paper, we study the deformation theory of degenerate algebraic curves
 on degenerate varieties.
Studying algebraic curves in varieties (usually assumed to be smooth, 
 or with mild singularities) plays crucial roles in many situations.
However, directly constructing curves with particular properties on smooth varieties
 is often difficult since such varieties are so smooth that one cannot
 readily find a foothold.
One way to overcome such difficulty is to deform the varieties to
 singular ones, where one can find more clues, and deduce information
 of curves on the original varieties from that on the singular varieties.
Such a study has been developed by many people, but it has become particularly
 prosperous since it was linked to its combinatorial counter part, 
 the geometry of tropical curves, see \cite{M, M2}.

In this paper, we consider
 the correspondence between 
 tropical curves in real affine 
 spaces and algebraic curves in toric varieties.
Previously this was studied in the case of so-called regular 
 tropical curves (see \cite[Definition 58]{N1}) by several authors \cite{CFPU, M, N1, NS, Tyo}.
In particular, in the paper \cite{N1}, 
 we solved this problem for arbitrary regular tropical curves
 (see also \cite{CFPU, Tyo}).
Superabundant tropical curves are those which 
 have larger dimensional parameter
 space than expected.
From the view point of algebraic curves, 
 this corresponds to the situation where
 the obstruction cohomology group 
 is non-trivial in the tangent-obstruction complex.

The general strategy to study the correspondence
 between tropical and algebraic curves 
 consists of two main steps:
\begin{enumerate}
\item First, construct a degenerate curve in a suitable ambient space, which is the central
 fiber of a degenerating family of varieties.
\item Then, deform the curve to a general fiber of the family.
\end{enumerate}

When the tropical curve is regular, the second problem is absent by definition of regularity.
Furthermore, the first problem can also be solved unconditionally. 
However, when the tropical curve is superabundant, 
 both problems become quite nontrivial.
We solved the first problem for curves of genus one in \cite{N1}
 (we will recall the result in Subsection \ref{subsec:4-vloop}).
In this paper, our main purpose is to examine whether 
 the obstruction cohomology group 
 really obstructs the smoothing of a curve in the degeneration limit (maximally
 degenerate curve and its variant, see Definition \ref{def:type}
 and Remark \ref{rem:type0}) or not.
It turns out that sometimes the obstruction class
 is effective and sometimes not, 
 and it can be read off from the combinatorics of tropical curves.

The first result in this direction was obtained by Speyer \cite{S}
 in the case of genus one, 
 where he discovered a sufficient condition for the existence of 
 a smoothing (the well-spacedness condition. See also \cite{BPR, Katz}).
In this paper, we use a different method which allows us to carry out
 a finer 
 and more general study of 
 superabundant curves.
This method clarifies the relation between obstructions and
 the combinatorics of tropical curves,
 so that degenerate cases (e.g., tropical curves with
 higher ($\geq 4$) valent vertices) and cases with curves of higher genus can be treated
 on the same footing.

In fact, it turns out that higher valent vertices play a critical role 
 in the study of correspondence between tropical curves and
 algebraic curves, though this point has not been appreciated much
 in literatures.
In the cases of curves of positive genus, tropical curves with higher valent
 vertices appear as general members in the set of those curves
 corresponding to classical algebraic curves
 (in other words, tropical curves which are 
 smoothable in the sense of Definition \ref{def:smoothable}).
The existence of higher valent vertices modifies the well-spacedness condition of 
 Speyer (see Definition \ref{well-spaced}),
 and in some cases these vertices cancel the obstruction for the smoothing
 by totally different mechanism compared to the known well-spacedness 
 condition (see Theorem \ref{thm:highvalloop}).\\

The main difficulty in our study
 is the existence of the obstructions to the deformation,
 which was absent
 in regular cases.
In modern deformation theory, deformation problems are formulated
 using differential graded Lie algebras.
The local model of the parameter space for deformed curves
 is obtained as the zero locus of 
  the \emph{Kuranishi map} whose existence 
  can be proved in such a general
  setting with moderate assumptions, see for example \cite{Ma}.
However, usually such an abstract formalism
 is not useful for the explicit calculation of the
 Kuranishi map.
Actually, to the best of the author's knowledge,
 there are few cases (for example, \cite{F, H})
 where the Kuranishi maps are 
 actually computed (so that the 
 corresponding deformation problem is explicitly solved).
See Subsection \ref{subsec:Kuranishi} for general remarks
 about the Kuranishi map in our situation.

The obstruction cohomology class 
 is usually represented as a \v{C}ech cocycle.
For example, suppose we have a pair of complex manifolds
 $Y\subset X$ and consider a deformation of $Y$ in $X$.
Taking suitable coordinate neighborhoods of $Y$ and $X$, 
 $Y$ can locally be expressed by vector valued analytic functions.
A local deformation of $Y$ is given by a perturbation of such functions.
The differences of these local deformations constitute a
 \v{C}ech 1-cocycle, and this represents the obstruction class
 (see Subsection \ref{subsec:Kuranishi} for more about this construction).
However, 
 in this calculation, we need to construct suitable coordinate neighborhoods, 
 and it is already not easy in general.
Even if we could do it and could construct the 1-cocycle, 
 it will be difficult to see whether it is zero as a cohomology class
 or not.
Thus, this point of view is not useful for actual calculation
 (and this is the reason the calculation of the Kuranishi map
 is difficult).

Our proposal is that degeneration technique is useful for the 
 calculation of the Kuranishi maps
 (both of the authors of \cite{F, H} deal with degenerate situations, but
 the use of degeneration in this paper is rather different from these papers).
The fundamental idea is very primitive.
Since we deal with analytic objects, essentially everything is determined by 
 the data in a germ of one point.
When we have
 a local deformation of $Y$ in some coordinate neighborhood as in the 
 above paragraph (so it is given by a set of vector valued analytic functions),
 if one can tell that whether or not this locally given set of functions
 can be analytically continued without developing divergence or multivaluedness, 
 then the deformation problem is solved.
Of course, calculating analytic continuation is usually very difficult, 
 and this way of thinking will be of little use as it is.

However, if one can at least partially calculate the analytic continuation, 
 one may be able to deduce some information about the obstruction class.
In particular, since the poles and multivaluedness prevents us from 
 constructing global deformation, if we can capture these behaviors of 
 the analytic continuation, we will be able to relate it to the 
 obstruction 1-cocycle.

This actually happens in  
 degenerate situations.
In degenerate situations, many of the difficulties appeared in the argument
 above resolve.
For example, both the ambient space and the curve in it degenerate to 
 the union of simple objects such as toric varieties and rational curves, 
 constructing coordinate neighborhoods becomes quite easy.
Also, since rational curves have global parameterizations, 
 the calculation of analytic continuation becomes viable on them, too.
Moreover, the obstruction 1-cocycle represented by the data of poles
 is much more tractable than the original \v{C}ech representative, 
 so that one can easily see whether the corresponding cohomology
 class is zero or not.

Besides these merits in practical calculation, our degeneration technique
 has the advantage that we need to know the ambient space
 only in a neighborhood of the curve, since once we know the data
 of the embedding of the curve 
 in a small neighborhood (in fact, in a germ of a single point), 
 then in principle the analytic continuation is completely determined.
Because of this notice, the idea in this paper can be applied to 
 various situations.
See \cite{N3, N5, NY} 
 for the calculation of the Kuranishi maps in other situations
 where the ambient spaces are
 degenerations of K3 surfaces, projective hypersurfaces, 
 complex tori, etc..\\

Now we give an outline of the contents of each section.
In Section \ref{sec:pre},
 we recall basic definitions concerning tropical curves
 with an emphasis on their relation to algebraic curves on toric varieties. 
We also recall the combinatorial 
 description of the dual space of the obstructions
 for degenerate algebraic curves on the central fiber of 
 toric degenerations of toric varieties given in
 \cite[Theorem 43]{N1}.
 
In Section \ref{sec:I},
 we calculate the obstructions for degenerate algebraic curves
 corresponding to tropical curves of genus one
 which are immersions.
The description of the dual obstruction cohomology classes of 
 superabundant curves (Theorem \ref{thm:obstruction})
 turns out to be very well fitted to this calculation.
In the calculation, we make use of degeneration technique in a
 crucial way.  
As we noted above, it is used to formulate the problem in a
 practical way in the first place.
Also, it is used during the actual calculation to reduce the arguments
 on general situations to those on lines in projective spaces.
Since the latter is defined by concrete equations, we can obtain rather
 explicit expression for the obstruction.

In Subsection \ref{subsec:Kuranishi}, we give general remarks on 
 Kuranishi maps.
In Subsections \ref{subsec:convention} to \ref{subsec:lifthigh},
 we give preparatory construction needed for later calculation.
In Subsection \ref{subsec:step3}, we give a definition of the obstruction 
 (Definition \ref{def:cechobst}) using the standard 
 construction based on \v{C}ech cohomology. 
The vanishing of it is 
 necessary and sufficient for the existence of the smoothing.
Then we perform the analytic continuation of functions on degenerate
 algebraic curves in Subsection \ref{subsec:step4}.
Based on it,
 we introduce the pre-obstruction using the data of 
 poles of analytic continuations of suitable functions
 in Subsection \ref{subsec:step5}
 (Definitions \ref{def:preob} and \ref{def:preobcontri}).
The pre-obstruction turns out to be equivalent to the original obstruction
 (Corollary \ref{cor:anal-ob}).
 
The vanishing of these obstructions can be understood as a
 Cousin type problem (that is, the problem of finding rational functions
 with prescribed singularities on a given Riemannian surface)
 through this calculation.
By the explicit expression of the pre-obstruction mentioned above, we can 
 solve this problem, resulting in
 the necessary and sufficient condition for
 an immersive superabundant tropical curve genus one
 $(\Gamma, h)$ under which there  
 is a pre-log curve associated to $(\Gamma, h)$ allowing a 
 smoothing (Theorem \ref{thm:immersive}).
The result turns out to be the extension of 
 the {well-spacedness condition} 
 of Speyer \cite{S} (see Definition 
 \ref{well-spaced}).
We note that when there are edges with weights larger than one, 
 the well-spacedness condition should be modified 
 (this point seems to be missing in the existing literature,
 see Remark \ref{rem:generalweight}).

An important point is that we consider tropical curves
 which need not be immersions, so that the images may have
 higher valent vertices.
We generalize the study of obstructions in Section \ref{sec:I}
 to these cases in Section \ref{sec:II}.
In contrast to the genus zero (or more generally regular)
 cases, such tropical curves generically appear in the parameter space of
 smoothable 
 (in the sense of Definition \ref{def:smoothable}) tropical curves.

There are two major contributions of the existence of higher valent
 vertices
 to the study of obstructions.
One is the case where the higher valent vertices are not on the 
 loop of the tropical curve, and the other is those on the loop. 
These play rather different roles.

In Subsection \ref{subsec:general_genus_one}, we study the former
 case, which is the direct 
 extension of the calculation in Section \ref{sec:I}.
Again, degeneration technique is crucial in that 
 although the curves corresponding to tropical curves with higher valent
 vertices are quite complicated, in the degenerate situation we can still
 reduce the calculation to the case of lines in projective spaces.
The latter case (i.e., higher valent vertices on the loop) 
 was studied in \cite{N1}, 
 which gives a new mechanism for the vanishing of the obstruction
 compared to the known well-spacedness in \cite{S}
 and Definition 
 \ref{well-spaced}.
We recall it
 in Subsection \ref{subsec:4-vloop} (Theorem \ref{thm:highvalloop}).

These results culminate in the generalized version of 
 the well-spacedness condition.
Roughly speaking, it goes like the following (see Definition 
 \ref{well-spaced2} for the precise statement).
\begin{defn}\label{def:intro}
A tropical curve $h\colon\Gamma\to \Bbb R^n$
 of genus one is \emph{well-spaced} if
 the following condition is satisfied.
Namely, let $A$ be the minimal affine subspace of $\Bbb R^n$
 which contains all the edges emanating from the vertices 
 contained in the loop of $h(\Gamma)$.
Let $B$ be any affine hyperplane of $\Bbb R^n$ containing $A$.
Let $\{v_1, \dots, v_k\}$ be the set of vertices on the subgraph
 $h(\Gamma)\cap B$ from which edges not contained in $B$ emanate.
Then,
\begin{itemize}
\item the set of distances from $v_i$ to the loop of $h(\Gamma)$
 contains at least two minimum, or
\item there is only one minimum, but from that vertex at least three edges
 emanates which are not contained in $B$.
\item Moreover, a certain polytope constructed from the vectors in the directions  
 of the edges emanating from the vertices in the loop of $h(\Gamma)$
 contains lattice points (see Theorem \ref{thm:119} for the precise definition).
\end{itemize}
\end{defn} 
The first two conditions concern the vanishing of the obstruction and the 
 third one concerns the existence of degenerate curves associated to the given tropical curve.
See \cite{T} for a related definition from purely combinatorial view point.
The following is our main result in this paper 
 (see Theorem \ref{thm:general} for the precise statement).
\begin{thm}
A tropical curve $(\Gamma, h)$ of genus one corresponds to
 a classical algebraic curve if and only if it is well-spaced
 in the sense of Definition \ref{def:intro}.
\end{thm}
Although this is an existence theorem, 
 the calculation of the obstruction contains sufficient information
 for more detailed study (for example, construction of enumerative
 invariants).
Concerning this point, we give a few results in 
 Subsection \ref{subsec:resol}.
In particular, we prove that tropical curves which correspond to 
 classical curves have a nice deformation property 
 (Proposition \ref{prop:nonassumption}).
Moreover, by the results of Torchiani \cite{T},
 the space of such tropical curves
 is of expected dimension and also of pure dimension,
 and in fact has a natural 
 structure of a tropical variety.
More detailed study of
 this space and its applications to classical and tropical curves
 will be done in the forthcoming paper.
 
Although we mainly concentrate on curves of genus one in this paper, 
 the techniques developed here and in \cite{N1}
 will almost suffice for the study of tropical curves of higher genus, 
 at least when the genus is not very large (when the genus is very large, 
 the combinatorial complexity would become formidable which prevents us
 from the detailed understanding of the deformation theory
 as in this paper, 
 although the mechanism of the appearance of the obstruction will be the same).
In fact, if the loop part of the tropical curve is a disjoint union of single loops, 
 then the method in this paper works 
 without substantial modifications even if the genus is arbitrarily large.
In general, the loops will be combined in a complicated way, and higher valent vertices 
 play important roles here, too 
 (see \cite[Example 85]{N1} for the calculation of one of such situations).
One direction which is lacking for the study of these general curves up to the present time
 is the detailed study of the case where 
 the map $h\colon\Gamma\to \Bbb R^n$ contracts the loop.
A related problem is that whether there exists a smoothable tropical curve
 which cannot deformed into an immersive one (see \cite[Question 72]{N1}). 
Proposition \ref{prop:nonassumption} claims that this has the positive answer 
 in the case of curves of genus one, and in this case
 essentially one do not need to consider the case
 where the loop is contracted.
However, in the case of curves of higher genus, such a study will be more important
 (see Example \ref{ex:high}).\\


\noindent
{\bf Acknowledgments.}
Needless to say, I was inspired by D. Speyer's paper \cite{S}.
My study of tropical geometry began from the joint work \cite{NS} with
 B. Siebert, and many ideas from it appear in this paper, too.
It is a great pleasure for me to express my gratitude to him.
Carolin Torchiani made an important observation 
 through Example \ref{ex:loop} which I
 missed in the first draft.
I also express my deep gratitude to her for it and also for sending me 
 her thesis \cite{T}.
This work was supported by  JSPS KAKENHI Grant Number 26400061.

\section{Preliminaries}\label{sec:pre}
In this section, we recall and define some notations and notions
 which are used in this paper.

\subsection{Tropical curves}\label{subsec:pre}
First we recall some definitions about tropical curves, see \cite{M, NS}
 for more information.
Let $\overline \Gamma$ be a weighted, connected finite graph.
Its sets of vertices and edges are denoted by $\overline \Gamma^{[0]}$,
 $\overline \Gamma^{[1]}$, respectively.
We write by $w_{\overline \Gamma} \colon 
  \overline \Gamma^{[1]} \to \Bbb N \setminus \{ 0 \}$
 the weight function.
An edge $E \in \overline \Gamma^{[1]}$ has adjacent vertices
 $\partial E = \{ V_1, V_2 \}$.
Let $\overline \Gamma^{[0]}_{\infty} \subset \overline \Gamma^{[0]}$
 be the set of all 1-valent vertices.
We write $\Gamma = \overline \Gamma \setminus \overline\Gamma^{[0]}_{\infty}$.
Noncompact edges of $\Gamma$ are called \emph{unbounded edges}.
Let $\Gamma^{[1]}_{\infty}$ be the set of all unbounded edges.
Let $\Gamma^{[0]}, \Gamma^{[1]}, w_{\Gamma}$
 be the sets of vertices and edges of $\Gamma$ and the weight function
 of $\Gamma$ (induced from $w_{\overline\Gamma}$ in 
 the obvious way),
 respectively.
Let $N$ be a free abelian group of rank $n\geq 1$
 and we write $N_{\Bbb K} = N\otimes_{\Bbb Z}\Bbb K$, where 
 $\Bbb K = \Bbb Q, \Bbb R, \Bbb C$.

\begin{defn}
We call a continuous map $h\colon \Gamma\to N_{\Bbb R}$ 
 a \emph{semi-affine map} if the following conditions hold.
\begin{enumerate}
\item $h$ is a proper map.
In particular, the closure of the image of an unbounded edge is non-compact.
\item For every edge $E \in \Gamma^{[1]}$, the restriction $h \big|_E$
is either an embedding with the image $h(E)$ 
contained in an affine line, or 
a contraction so that $h(E)$ is a point.
\end{enumerate}
\end{defn}

Let $h$ be a semi-affine map.
Let $\mathcal S\subset h(\Gamma)$ be the set of points
 with the property that $p\in \mathcal S$ if and only if
 for any neighborhood $O_p$ of $p$ in $N_{\Bbb R}$, 
 the intersection $O_p\cap h(\Gamma)$ is not
 homeomorphic to an open interval.
The following is easy to see.
\begin{lem}\label{lem:savert}
Let $h$ be a semi-affine map as above.
Then the following statements hold.
\begin{enumerate}
\item $\mathcal S$ is a finite set.
\item By adding finite number of vertices to $\Gamma$, 
 we can assume that the inverse image $h^{-1}(p)$ 
 of each $p\in\mathcal S$
 consists of closed subgraphs of $\Gamma$.\qed
\end{enumerate}
\end{lem} 

Under the condition of Lemma \ref{lem:savert} (2), the image $h(\Gamma)$
 has a natural structure of a graph as follows.

\begin{defn}\label{def:imagevert}
Let $h\colon \Gamma\to N_{\Bbb R}$ be a semi-affine map satisfying 
 the condition of Lemma \ref{lem:savert} (2).
Then a vertex of the image $h(\Gamma)$ is the image of some vertex of $\Gamma$.
Similarly, an edge of $h(\Gamma)$ is the image of some edge of $\Gamma$
 which is not contracted.
\end{defn}

Note that under the condition of Lemma \ref{lem:savert} (2), 
 the inverse image $h^{-1}(v)$ of a vertex of $h(\Gamma)$ is a union of 
 closed subgraphs.
Similarly, the inverse image $h^{-1}({\mathfrak E}^{\circ})$ of 
 the open part ${\mathfrak E}^{\circ}$ of an edge 
 $\mathfrak E$ of $h(\Gamma)$ is the union of the open part of the edges of $\Gamma$
 each of which is mapped to $\mathfrak E$ homeomorphically.

Henceforth, we always assume 
Lemma \ref{lem:savert} (2) holds.
Let $p\in\mathcal S$ and $\Gamma_1$ be one of connected components of
$h^{-1}(p)$.
Then $\Gamma_1$ contains several 1-valent vertices $q_1, \dots, q_a$.
Let $E_1, \dots, E_b$ be the edges of $\Gamma\setminus\Gamma_1$
emanating from some of 
$q_1, \dots, q_a$.
\begin{defn}[{\cite[Definition 2.2]{M}}]\label{def:param-trop}
A \emph{parametrized tropical curve} in $N_{\Bbb R}$ is a semi-affine map
 $h \colon \Gamma \to N_{\Bbb R}$ satisfying the following conditions.
\begin{enumerate}
\item[(i)] For every edge $E \in \Gamma^{[1]}$, 
 the image $h(E)$ is either contained
 in an affine line with a rational slope, or 
 a point.
\item[(ii)] For every vertex $V \in \Gamma^{[0]}$, $h(V)\in N_{\Bbb Q}$.
\item[(iii)] 
 The following \emph{balancing
 condition} holds.
Namely, for each $p\in\mathcal S$ and a connected component
 $\Gamma_1$ of $h^{-1}(p)$, 
 the equality
\begin{equation}
\sum_{j=1}^b w_{\Gamma}(E_j)u_j = 0
\end{equation}
 holds, using the notation in the above paragraph.
Here $u_j$ is the primitive integral vector of $N$ 
 in the direction of the edge $h(E_j)$ emanating from $p$.
\end{enumerate}
\end{defn}
\begin{rem}
In \cite{NS}, $h|_E$ is assumed to be an embedding
 (see \cite[Definition 1.1]{NS}) for every edge $E$.
The reason that we take the above definition is that such cases appear naturally
 when we consider superabundant tropical curves.
\end{rem}
An isomorphism of parametrized tropical curves 
 $h \colon \Gamma \to N_{\Bbb R}$ and 
 $h' \colon \Gamma' \to N_{\Bbb R}$ 
 is a homeomorphism $\Phi \colon \Gamma \to \Gamma'$
 respecting the weights such that $h = h' \circ \Phi$.
\begin{defn}\label{def:tropical curve}
	A \emph{tropical curve} is an isomorphism class of parametrized 
	tropical curves.
	A tropical curve is \emph{3-valent} if any vertex of $\Gamma$ is at  most
	3-valent.
	The \emph{genus} of a tropical curve is the first Betti number of $\Gamma$.
	The set of \emph{flags} of $\Gamma$ is 
	\[
	F\Gamma = \{(V, E) \big|
	V \in \partial E \}.
	\]
\end{defn}
Note that in our definition, 
there can be 2-valent vertices in the abstract graph $\Gamma$
of a 3-valent tropical curve.
By (i) of Definition \ref{def:param-trop}, we have a map
$u \colon F\Gamma \to N$ sending a flag $(V, E)$
to the primitive integral vector $u_{(V, E)} \in N$
emanating from $h(V)$ in the direction of $h(E)$ or to the zero vector. 
\begin{rem}\label{rem:0edge}
We allow that even when $u_{(V, E)}\neq 0$, the edge $E$ is contracted
 by the map $h$
 $($but in such a case the direction $u_{(V, E)}$ does not contribute to the 
 balancing condition$)$.
We think of such a tropical curve as the limit of a
 family of tropical curves $(\Gamma, h_s)$, 
 $s\in [0, 1)\cap\Bbb Q$, where the image $h_s(E)$ has the direction $u_{(V, E)}$
 for each $s$ 
 and the length goes to 0 as $s\to 1$
 (see Definition \ref{def:moduli} and Remark \ref{rem:type}).
\end{rem}
\begin{defn}\label{def:type}
The (unmarked) \emph{combinatorial type}
 of a tropical curve $(\Gamma, h)$
 is the graph $\Gamma$ 
 together with the map $u \colon F\Gamma \to N$.
We write this by the pair $(\Gamma, u)$.
\end{defn}

\begin{defn}\label{def:degree}
The \emph{degree} of a combinatorial type $(\Gamma, u)$
 is the function $\Delta(\Gamma, u)=\Delta\colon N \setminus \{ 0 \}
  \to \Bbb N$
 with finite support defined by 
\begin{equation*}
 \Delta(\Gamma, u)(v):= \sharp \{ (V, E) \in F\Gamma\; |\;
    E \in \Gamma^{[1]}_{\infty}, w(E)u_{(V, E)} = v \} 
\end{equation*}
Let $e=|\Delta| = \sum_{v\in N\setminus\{0\}}\Delta(v)$.
This is the same as the number of unbounded edges 
 of $\Gamma$ (not necessarily of $h(\Gamma)$
 since some of the edges may have the same image).
\end{defn}
\begin{defn}\label{immersive}
	We call a tropical curve $(\Gamma, h)$ \emph{immersive} 
	if for any $E\in\Gamma^{[1]}$, the restriction of 
	$h$ to $E$ is an embedding.
	Note that even if $(\Gamma, h)$ is immersive, some of the edges of $\Gamma$
	can have the same image. 
\end{defn}
\begin{prop}[{\cite[Proposition 2.13]{M}}]\label{prop:trop_moduli}
The space parameterizing immersive 3-valent tropical
 curves of a given combinatorial type is, if it is
 non-empty, an open
 convex polyhedral domain in the real affine $k$-dimensional space,
 where $k \geq e+(n-3)(1-g)$. 
Here $e$ is the number of unbounded edges of $\Gamma$ as
 in Definition \ref{def:degree}, $n$ is the dimension of the target space $N_{\Bbb R}$,
 and $g$ is the genus of $\Gamma$.\qed
\end{prop}
\begin{defn}\label{def:moduli}
Fix a combinatorial type of 3-valent tropical curves whose parameter space of 
 immersive curves is non-empty.
Then we define
 the \emph{parameter space} of all 3-valent tropical
 curves of the given combinatorial type as the closure
 of the parameter space of 
 immersive curves.
\end{defn}
\begin{rem}\label{rem:type}
In other words, an element of the parameter space of 
 tropical curves of a given combinatorial type is a tropical curve
 which can be deformed into an immersive tropical curve
 of that combinatorial type.
If there is no immersive tropical curve of the given combinatorial type, 
 then the corresponding parameter space is defined to be empty
 in this paper.
This does not give any restriction 
 in the case of genus one due to Proposition \ref{prop:nonassumption}.
\end{rem}
\begin{defn}[{\cite[Definition 2.22]{M}}]\label{superabundant}
A 3-valent tropical curve is called \emph{superabundant}
 if the parameter space is of dimension larger than 
 $e+(n-3)(1-g)$.
\end{defn}
In the case of genus zero (\cite{NS}), or more generally, 
 if a tropical curve is regular (see Definition \ref{def:nonsuperabundant2} below),
 one can see that if $(\Gamma, h)$ is a tropical curve
 which satisfies general incidence conditions
 ${\bf A} = (A_1, \dots, A_l)$
 of codimension ${\bf d} = (d_1, \dots, d_l)$
 (see \cite[Definition 2.3]{NS})
  with $|{\bf d}|=\sum_i d_i = (n-3)(1-g)+e$,
  then $h$ is an immersion (embedding if $n$ is greater than two).

However, in superabundant cases, the situation 
 where the map $h$ contracts some of the
 edges of $\Gamma$ naturally appears
 in view of correspondence with algebraic curves
 (see for example Example \ref{ex:loop}).
Consequently, the image 
 $h(\Gamma)$ may have vertices of higher $(>3)$
 valence, and some of the edges
 of $\Gamma$ may have the same image.
On the other hand, 
 we do not need to allow all kinds of maps $h$ for the enumerative results
 (for example, maps which contract a loop are not generic 
 in the space of tropical curves of genus one which correspond to algebraic curves, 
 see Proposition \ref{prop:nonassumption}).
In this paper, 
 we assume that a tropical curve $(\Gamma, h)$ satisfies
 Assumption A below, which is a reasonable condition for the study of 
 superabundant tropical curves in view of correspondence with algebraic curves
 (see also Remark \ref{rem:assumption}).

Before stating Assumption A, we prepare some terminologies.
Let $\Gamma$ be a finite non-compact graph as above.
\begin{defn}\label{def:loops}
\begin{enumerate}[(i)]
\item An edge $E\in\Gamma^{[1]}$ is said to be a 
 \emph{part of a loop} of $\Gamma$
 if the graph given by $\Gamma\setminus E^{\circ}$ has the smaller first Betti number
 than $\Gamma$.
Here $E^{\circ}$ is the interior of $E$ (that is, 
 $E^{\circ} = E\setminus\partial E$).
\item The \emph{loop part} of $\Gamma$ is the subgraph of $\Gamma$ composed of
 the union of  parts of a loop of $\Gamma$.
\item A \emph{bouquet} of $\Gamma$ is 
 a connected component of the loop part of $\Gamma$.
A subset of a bouquet which is homeomorphic to
 a circle is called a \emph{loop}.
\end{enumerate}
\end{defn}
In particular, a bouquet or a loop does not contain unbounded edges.
Now we state Assumption A. \\

\noindent
{\bf Assumption A.}
\begin{enumerate}[(i)]
\item A vertex of the abstract graph $\Gamma$ is at most 3-valent.
Therefore, $(\Gamma, h)$ is always a 3-valent tropical curve in view of
 Definition \ref{def:tropical curve},
 although some vertices of the image $h(\Gamma)$ may not be 2- or 3-valent.
\item The map $h$ may contract some of the bounded edges of $\Gamma$.
However, $h$ does not contract a loop to a point.
\item $(\Gamma, h)$ can be deformed into an immersive tropical curve.
Thus, the parameter space containing $(\Gamma, h)$ is not empty.
\item The inverse image of any vertex of $h(\Gamma)$ is a 
 disjoint union of closed subgraphs of $\Gamma$ (by (ii), 
 each component must be a tree or a vertex).
\end{enumerate}
The essential part of Assumption A is the part (iii). 
The part (ii) is also related to the part (iii) (see \cite[Example 19]{N1}.
See also Remark \ref{rem:assumption}).
The part (iv) is (2) of Lemma \ref{lem:savert}, 
 and it can always be realized by taking
 a suitable refinement of  the graph
 $\Gamma$.
In fact, this Assumption A does not give any restriction 
 if one is interested only in tropical curves of genus one which 
 correspond to some classical curves 
 (see Proposition \ref{prop:nonassumption}).

Note that immersive 3-valent tropical curves satisfy Assumption A
 possibly after a suitable refinement of the graph $\Gamma$.
When $(\Gamma, h)$ satisfies Assumption A, we
 define a part of a loop, the loop part, a bouquet and a loop 
 of the image 
 $h(\Gamma)$ as the images of those of $\Gamma$.

In order to distinguish the valence and the weights of $\Gamma$ 
 from those of $h(\Gamma)$, we introduce the following definition. 
\begin{defn}\label{def:trop}
We assume that a tropical curve $(\Gamma, h)$ satisfies Assumption A.
\begin{enumerate}[(i)]
\item Let $v\in h(\Gamma)$
  be a vertex and let $\gamma_1, \dots, \gamma_s$
 be the connected components of $h^{-1}(v)$.
The \emph{valence} of $v$, $val(v)$ is defined as follows.
Namely, by Assumption A, $\gamma_1, \dots, \gamma_s$
 are closed subgraphs of genus zero in
 $\Gamma$.
Contracting $\gamma_i$ to a point 
 produces a new graph with a vertex $W_i$ which is the 
 image of $\gamma_i$ under this contraction.
Let $a_i$ be the valence of this vertex in the new graph.
Then, the valence of the vertex $v$ is the unordered set of positive integers
\[
val(v) = \{a_1, \dots, a_s\}.
\]
\item Let $\mathfrak E\in h(\Gamma)$ be an edge.
Let $E_1, \dots, E_s\in \Gamma^{[1]}$ be the edges of $\Gamma$ such that 
 $h(E_i) = \mathfrak E$
 (in particular, the edges $E_1, \dots, E_s$ are not contracted by $h$).
Then the \emph{weight} of $ \mathfrak E$, $w(\mathfrak E)$,
 is the unordered set of positive integers
 $\{w_1, \dots, w_s\}$, here $w_i$ is the weight of $E_i$ in $\Gamma$.
The sum $w_s(\mathfrak E) = 
 \sum_{i=1}^{s}w_i$ is called the \emph{total additive weight} of 
 $\mathfrak E$.
The product $w_m(\mathfrak E) = 
 \prod_{i=1}^sw_i$ is called the \emph{total weight} of $\mathfrak E$.
\end{enumerate}
\end{defn}

Finally, although our main concern in this paper is 
 properties of superabundant tropical
 curves and associated algebraic curves, we recall the definition of regular tropical curves 
 for readers' convenience since
 it will be referred to several times.

Let $\Gamma$ be a weighted abstract graph as in the beginning of 
Subsection \ref{subsec:pre}, but now we allow
that $\Gamma$ is not necessarily 3-valent.
Let $h\colon \Gamma\to N_{\Bbb R}$ be an embedding giving $\Gamma$
a structure of a tropical curve.
We identify the graph $\Gamma$ with its image.
Let $L$ be the loop part of $\Gamma$ 
and $L^{[1]} = \{E_1, \dots, E_l\}$ be the set of edges in $L$
seen as 1-chains of the simplicial homology group
(in particular, we choose an orientation for each $E_i$).
Let $u_i$ be the primitive integral vector in the direction of $E_i$.
\begin{defn}[{\cite[Definition 4.1]{CFPU}, see also \cite[Section 2.6]{M}
		and \cite[Section 1]{Katz}}]\label{def:nonsuperabundant2}
	The tropical curve $(\Gamma, h)$ is \emph{regular}
	if the following abundancy map
	\[
	\Phi_{(\Gamma, h)}\colon
	\Bbb R^{L^{[1]}}\to Hom(H_1(\Gamma), N_{\Bbb R}),
	\]
	\[
	(\ell_{E_i})\mapsto \left(\sum_i a_i[E_i]\mapsto \sum_i\ell_{E_i}a_iu_i
	\right)
	\]
	is surjective.
\end{defn}
\begin{rem}
	In \cite{CFPU}, a tropical curve satisfying this condition is called non-superabundant.
	We call it differently, since some people already call these tropical curves regular, 
	and we used the word non-superabundant in \cite[Definition 17]{N1}
	to refer to different objects.
Also, from the point of view of Theorem \ref{thm:regsm} below, the use of the
 word regular for these tropical curves
 seems reasonable.
\end{rem}
\begin{example}
Immersive 3-valent plane tropical curves or immersive tropical curves 
 of genus zero in 
 any dimensional ambient space
 are regular.
\end{example}

Regular tropical curves behave very nicely in view of correspondence 
 to algebraic curves.
Namely, in \cite{N1}, we proved the following.
\begin{thm}\label{thm:regsm}
	Any regular tropical curve is smoothable
	in the sense of Definition \ref{def:smoothable} below.\qed
\end{thm}

\subsection{Toric varieties associated 
 to tropical curves and pre-log curves on them}\label{subsec:toric}
\begin{defn}\label{toric}
A toric variety $X$ defined by a fan $\Sigma$ 
 is called \emph{associated to a tropical curve $(\Gamma, h)$}
 if the set of the 
 rays of $\Sigma$ contains the set of the rays spanned by the vectors in $N$
 which are in the 
 support of the degree map $\Delta\colon N\setminus\{0\}\to \Bbb N$
 of $(\Gamma, h)$.

If $ \mathfrak E$
 is an unbounded edge of $h(\Gamma)$, there is 
 the obvious unique divisor of $X$
 corresponding to it.
We write it as $D_{\mathfrak E}$ and call it the \emph{divisor associated to the edge 
 $\mathfrak E$}.
\end{defn}
\begin{defn}\label{def:degeneration}
Given a tropical curve $(\Gamma, h)$ in $N_{\Bbb R}$
 defined over $\Bbb Q$ (that is, the vertices have rational coordinates),
 we can construct a polyhedral
 decomposition $\mathscr P$ of $N_{\Bbb R}$ defined over $\Bbb Q$
 such that $h(\Gamma)$ is contained in the 1-skeleton of $\mathscr P$
 (\cite[Proposition 3.9]{NS}).
Given such $\mathscr P$, 
 we construct a degenerating family $\mathfrak X\to \Bbb C$
 of
 a toric variety $X$ associated to $(\Gamma, h)$ (\cite[Section 3]{NS}).
We call such a family a \emph{degeneration of $X$ defined respecting $(\Gamma, h)$}.
Let $X_0$ be the central fiber.
It is a union $X_0 = \cup_{v\in\mathscr P^{[0]}}X_{0,v}$
 of toric varieties intersecting along toric strata.
Here $\mathscr P^{[0]}$ is the set of the vertices of $\mathscr P$.
\end{defn}
\begin{rem}\label{rem:basechange}
	To define $\mathfrak X$, in general we need to multiply $(\Gamma, h)$ by some constant
	so that it becomes defined over $\Bbb Z$.
	There are choices of this multiplication factor, which correspond to 
	base changes on the algebraic geometry side.
	One requirement to this choice is that the integral length of each bounded edge is a multiple of
	its weight (see \cite[Proposition 7.1]{NS}).
If this condition is satisfied, the difference of the multiplication factor has essentially no
 effect to the results (see also Remark \ref{rem:2-val}).
In the following arguments,
	we just talk about degenerations 
	defined respecting $(\Gamma, h)$ which is defined over $\Bbb Q$, 
	assuming a suitable choice of the multiplication factor is made.
It is convenient to do so, since we sometimes refine the graph $\Gamma$
 by adding 2-valent vertices, and mentioning the base change each time will be
 cumbersome.
\end{rem}
\begin{defn}[{\cite[Definition 4.1]{NS}}]\label{torically transverse}
Let $X$ be a toric variety.
An algebraic curve $C\subset X$ is \emph{torically transverse}
 if it is disjoint from all toric strata of codimension greater than 1.
A stable map $\phi\colon
 C\to X$ is torically transverse if $\phi^{-1}(int X)\subset C$
 is dense and $\phi(C)\subset X$ is a torically transverse curve. 
Here $int X$ is the complement of the union of toric divisors.
\end{defn}
\begin{defn}\label{def:pre-log}
Let $C_0$ be a prestable curve.
A \emph{pre-log curve} on $X_0$ is a stable map 
 $\varphi_0\colon C_0\to X_0$
 with the following properties.
\begin{enumerate}
\item[(i)] For any $v\in\mathscr P^{[0]}$,
 the restriction $C_0\times_{X_0}X_{0,v}\to X_{0,v}$
 is a torically transverse stable map.
\item[(ii)] Let $P\in C_0$ be a point which maps to the singular locus of $X_0$.
Then $C_0$ has a node at $P$, and $\varphi_0$ maps the two branches
 $(C_0', P), (C_0'', P)$ of $C_0$ at $P$ to different irreducible components 
 $X_{0,v'}, X_{0,v''}\subset X_0$.
Moreover, if $w'$ is the intersection index 
 of the restriction $(C_0', P)\to (X_{0,v'}, D')$ with the toric divisor
 $D'\subset X_{0,v'}$, 
 and $w''$ accordingly for $(C_0'', P)\to (X_{0,v''}, D'')$,
 then $w' = w''$.
\end{enumerate}
\end{defn}

Let $X$ be a toric variety and $D$ be the union of toric divisors.
In \cite[Definition 5.1]{NS}, a non-constant
 torically transverse map $\phi\colon \Bbb P^1\to X$ is called a \emph{line}
 if 
\[
\sharp\phi^{-1}(D)\leq 3.
\]
In this case, the image of $\phi$ is contained in the closure of 
 the orbit of a subtorus of dimension at most two 
 of the big torus acting on $X$ (\cite[Lemma 5.2]{NS}).
Because we consider 
 more general tropical curves, we have to extend this notion.

Let $\Gamma$ be a weighted 3-valent tree 
 and $h\colon \Gamma\to N_{\Bbb R}$
 be a map which gives $\Gamma$ a structure of a tropical curve,
 and assume that the image $h(\Gamma)$ has only one vertex $v$.
Let $\mathfrak E_1, \dots, \mathfrak E_s$ be the edges of $h(\Gamma)$
 (some of 
 these can coincide, and
 these correspond to the unbounded edges of $\Gamma$ in a natural way).
Let $X$ be a toric variety associated to $(\Gamma, h)$.
\begin{defn}\label{def:typev}
A non-constant torically transverse map
 $\phi\colon \Bbb P^1\to X$ is called 
 \emph{of type $(\Gamma, h)$}, or  \emph{of type $v$}
 when $h$ is clear from the context,
 if $\phi$ satisfies the following property:
\begin{itemize}
\item Let $\mathfrak E_{i}$ be an edge of $h(\Gamma)$ and
 let $w_i$ be the weight of $\mathfrak E_{i}$.
Then $\phi(\Bbb P^1)$ has an intersection
 with the divisor $D_{\mathfrak E_i}$
 with 
  intersection multiplicity 
 $w_i$, and there is no intersection between $\phi(\Bbb P^1)$
 and toric divisors other than these.
\end{itemize} 
\end{defn}
Note that when some edges of $\mathfrak E_1, \dots, \mathfrak E_s$
 coincide, then there are several intersections between $\phi(\Bbb P^1)$
 and the corresponding toric divisor.
 
\begin{figure}[h]
	\includegraphics[height=4cm]{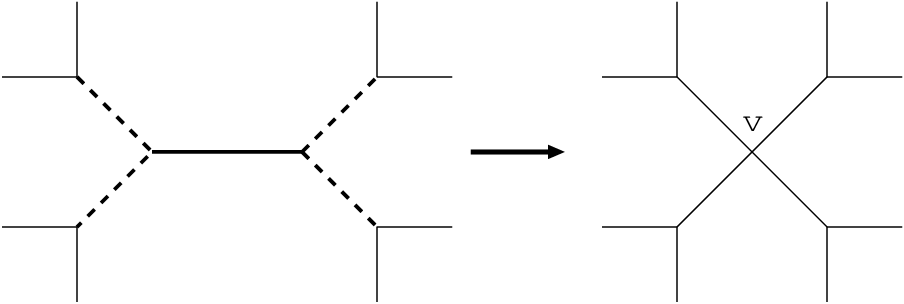}
	\caption{The abstract graph $\Gamma$ (the picture on the left)
		is mapped to a graph in $\Bbb R^n$ (the picture on the right).
		The bold line segment in $\Gamma$ is the inverse image of the vertex $V$
		and the union of the bold line segment and the dotted line segments
		in $\Gamma$ is the subgraph $\Gamma_V$.
		In this example, the graph $\Gamma_{\varphi_0}$
		is the same as the picture on the right as an abstract graph.}\label{fig:invim}
\end{figure}

Let $(\Gamma, h)$ be a tropical curve satisfying Assumption A.
Let $X$ be a toric variety associated to $(\Gamma, h)$ and  
 $\mathfrak X\to\Bbb C$ be a degeneration of $X$ defined respecting $(\Gamma, h)$.
Let $X_0$ be the central fiber.
\begin{defn}\label{def:oftype}
A pre-log curve $\varphi_0\colon C_0\to X_0$ is called 
 \emph{of type $(\Gamma, h)$}
 if for any $v\in h(\Gamma)^{[0]}$,  
 the restriction of $C_0\times_{X_0}X_{0, v}\to 
  X_{0, v}$ to each of the components $\{C_{0, v(i)}\}$
  of $C_0$
  mapped to $X_{0, v}$
  is a rational torically transverse curve of type 
  $(\Gamma_{v(i)}, h|_{\Gamma_{v(i)}})$.
 Here
  $\{\Gamma_{v(i)}\}$ is the set of open
  subgraphs of $\Gamma$ each element of which is 
  the union of 
 \begin{itemize}
 \item a connected component of $h^{-1}(v)$, and 
 \item the open
  part of the edges emanating from the vertices of it.
 \end{itemize}
See Figure \ref{fig:invim} for an example.
\end{defn}

\begin{defn}\label{def:intergraph}
When $\varphi_0\colon C_0\to X_0$ is a pre-log curve of type
 $(\Gamma, h)$, let $\Gamma_{\varphi_0}$ be the graph obtained from 
 $\Gamma$ by contracting 
 all the bounded edges of each $\Gamma_{v(i)}$
 in the notation of Definition \ref{def:oftype}.
The graph $\Gamma_{\varphi_0}$ is the dual intersection graph of $C_0$
 (containing the set of unbounded edges, which is
 the information of intersection of $\varphi_0(C_0)$
 and come of the toric divisors of $X_0$ seen as marked points of $C_0$)
 and it is clear that the map $h\colon \Gamma\to \Bbb R^n$
 factors through $\Gamma_{\varphi_0}$.
\end{defn}

By the condition for pre-log curves, the components $\{C_{0, v(i)}\}$ 
 are disjoint and the subgraphs $\{\Gamma_{v(i)}\}$
 are also disjoint.
Note that each $C_{0, v(i)}$ is an irreducible rational curve.

\begin{rem}\label{rem:type0}
If $(\Gamma, h)$ is immersive, 
 then a pre-log curve of type $(\Gamma, h)$ is 
 the \emph{maximally degenerate curve} of \cite[Definition 5.6]{NS}.
\end{rem}
\begin{defn}\label{def:smoothable}
A 3-valent 
 tropical curve $(\Gamma, h)$  is \emph{smoothable}
 if there is a pre-log curve $\varphi_0\colon C_0\to X_0$ 
 of type $(\Gamma, h)$ with the following property.
Namely, there exists 
 a family of stable maps over $\Bbb C$
\[
\Phi\colon \mathfrak C/\Bbb C\to \mathfrak X/\Bbb C
\]
 such that $\mathfrak C/\Bbb C$ is a flat family of pre-stable 
 curves whose fiber over $0$ is isomorphic to $C_0$,
 and the restriction of $\Phi$ to $C_0$ is a stable map
 equivalent to $\varphi_0$.
We also call such a pre-log curve \emph{smoothable}.
\end{defn}
\begin{rem}
The smoothability of a tropical curve does not depend on the choice of
 a toric variety
 $X$ associated to it or a degeneration of $X$ defined respecting the tropical curve.
\end{rem}
See \cite[Section 5]{NS}, for more information about lines and maximally
  degenerate pre-log curves.
Given an immersive 3-valent tropical curve
 $(\Gamma, h)$ which is non-superabundant,
 we can construct maximally degenerate 
 pre-log
 curves of type $(\Gamma, h)$
 (\cite[Proposition 5.7]{NS}), and vice versa (\cite[Construction 4.4]{NS}).
More generally, in \cite[Proposition 77]{N1}, we showed that if $(\Gamma, h)$
 is a regular tropical curve, then there is a pre-log curve of type 
 $(\Gamma, h)$ (it is not difficult to see that immersive non-superabundant tropical curves
 are regular).

However, when the tropical curve $(\Gamma, h)$ is superabundant, 
 there are examples to which pre-log curves of type $(\Gamma, h)$
 do not exist, even when $(\Gamma, h)$ is an embedded tropical curve
 (\cite[Examples 81, 82]{N1}).
On the other hand, in the case of genus one, there are always pre-log curves
 of the given type, when $(\Gamma, h)$ is an embedding (or immersion).
We deal with the deformations of degenerate curves in these cases
 in Section \ref{sec:I}.
When $(\Gamma, h)$ is not an immersion, the situation 
 can be quite different.
We deal with the deformation theory in these cases in 
 Section \ref{sec:II}.
  
The smoothings of the maximally degenerate curves or pre-log curves
 of type $(\Gamma, h)$ are studied
 by log-smooth deformation theory \cite{KF, KK}.
For informations about 
 log structures relevant to our situation, see \cite[Section 7]{NS}.
We do not repeat it here, 
 because nothing new about log structures is required here, 
 other than those given in \cite{NS}.

\subsection{Description of the dual obstruction}\label{subsec:abundance}
In this subsection, we recall the description of the dual of the
 cohomology group
 where the obstructions to the smoothing of pre-log curves lie.
First we recall some terminology concerning log deformation theory of 
 pre-log curves.
See \cite{NS, N1} for more details.

Let $\mathfrak X\to \Bbb C$ be a degenerating family of toric varieties 
 defined respecting a given tropical curve $(\Gamma, h)$, in the sense
 of Definition \ref{def:degeneration}.
Let $\varphi_0\colon C_0\to X_0$ be a pre-log curve of type
 $(\Gamma, h)$ (see Definition \ref{def:type}),
 where $X_0$ is the central fiber of $\mathfrak X$.

The total space of the degeneration $\mathfrak X$ is a toric variety
 and has a natural log structure given by the toric divisors.
The central fiber $X_0$ has a log structure by restriction.
The curve $C_0$ also has natural log structures
 so that the map $\varphi_0$ extends to a log smooth map
 (see \cite[Section 7]{NS}).
We write the log smooth map by the same letter $\varphi_0$.
Let $\Theta_{X_0}$ and $\Theta_{C_0}$
 be the log tangent sheaves of $X_0$
 and $C_0$, respectively.
Let $\mathcal N_{C_0/X_0} = \varphi_0^*\Theta_{X_0}/\Theta_{C_0}$
 be the log normal sheaf.
Then the obstruction classes lie in the cohomology group
 $H^1(C_0, \mathcal N_{C_0/X_0})$.
Let $H = H^1(C_0, \mathcal N_{C_0/X_0})^{\vee}$ be the dual space.

By the Serre duality for nodal curves, it is isomorphic to 
 $H^0(C_0, \mathcal N_{C_0/X_0})^{\vee}\otimes\omega_{C_0})$, 
 where $\omega_{C_0}$ is the dualizing sheaf of the nodal curve $C_0$.
It is the sheaf of 1-forms where logarithmic poles are allowed at the nodes of $C_0$.
Note that the sheaf $\mathcal N_{C_0/X_0}^{\vee}$ is a subsheaf of 
 $(\varphi_0^*\Theta_{X_0})^{\vee}$ which is naturally isomorphic to 
 $N_{\Bbb C}^{\vee}\otimes\mathcal O_{C_0}$.

Let $(\Gamma, h)$ be a 3-valent embedded tropical curve.
We identify the edges of $\Gamma$ and their image by $h$.
Let $L$ be the loop part of $\Gamma$ (see Definition \ref{def:loops}).
If $L_i\subset L$ is a bouquet of $\Gamma$, 
 note that it is a subgraph of $\Gamma$ with 2-valent and 3-valent vertices.
Let $\{ v_i\}$ be the set of 3-valent vertices of $L$.
Cutting $L$ at each $v_i$, we obtain a set of piecewise linear segments 
 $\{l_m\}$.
Let $U_m$ be the linear subspace of $N_{\Bbb C}$ spanned by the
 direction vectors of the segments of $l_m$.

In \cite[Theorem 44]{N1}, we obtained the following description of the
 space $H$.

\begin{thm}\label{thm:obstruction}
Let $(\Gamma, h)$ be a 3-valent embedded
 tropical curve.
Elements of the space $H$ are described by the following procedure.
\begin{enumerate}[(I)]
\item Give the value zero to all the flags not contained in 
 the loop part of $\Gamma$.
\item Give a value $u_m$
 in $(U_m)^{\perp}\subset (N_{\Bbb C})^{\vee}$
 to each of the flags  associated to the edges of $\{h(l_m)\}$.
\item The data $\{u_m\}$ give an element of $H$ if and only if the following 
 conditions are satisfied.
\begin{enumerate}
\item At each vertex
 $v$ of $h(\Gamma)$, 
\[
u_1+u_2+u_3=0
\]
 holds as an element of $(N_{\Bbb C})^{\vee}$.
Here $u_1, u_2, u_3$ are the data attached to the three flags in
 $h(\Gamma)$
 which have $v$ as the vertex.
\item The data $\{u_m\}$
 is compatible on each edge of $h(l_m)$, 
 in the sense that the sum of the values attached to the two flags
 of any edge of $h(l_m)$ is zero. \qed
\end{enumerate}
\end{enumerate} 
\end{thm}
In the case where $(\Gamma, h)$ is a superabundant curve of genus one, 
 the set of generators of the space $H$ is easily written down.
Namely, in this case the space $H$ is canonically isomorphic to 
 the space $(\bar A_{\Bbb C})^{\perp}\subset (N_{\Bbb C})^{\vee}$,
 where $A$ is the minimal affine subspace of 
 $N_{\Bbb R}$ which contains the loop of $h(\Gamma)$, 
 $\bar A$ is the linear subspace of $N_{\Bbb R}$ parallel to
 $\bar A$, and $\bar A_{\Bbb C}$ is the complexification of $\bar A$.
Let $f_1, \dots, f_r$ be a basis of $(\bar A_{\Bbb C})^{\perp}$.

Let $C_{0,1}, C_{0,2}, \dots, C_{0,s}, C_{0,s+1} = C_{0,1}$
 be the loop of rational curves in 
 $C_0$ corresponding to the loop of $\Gamma$.
Let $p_i = C_{0,i}\cap C_{0,i+1}$ be the node ($i = 1, \dots, s$).
Let $z_i$ be an affine coordinate on $l_i$ which is 
 $0$ at $p_{i-1}$ and $\infty$ at $p_i$.
Then the logarithmic 1-form on $C_0$ which is 
 $\frac{dz_i}{z_i}$ on $C_{0,i}$ is the generator of the 
 space of sections $H^0(C_0, \omega_{C_0})$
 of the dualizing sheaf $\omega_{C_0}$.
Let us write it by $\xi$.
\begin{cor}\label{cor:obstbasis}
The set of sections 
\[
f_1\otimes \xi, \dots, f_r\otimes \xi
\] 
 gives a basis of the space $H$.
\end{cor}
\proof 
This is an immediate consequence of Theorem \ref{thm:obstruction}.
See \cite{N1} for more details.\qed\\
 
Although Theorem \ref{thm:obstruction} assumes that
 $(\Gamma, h)$ is an embedded curve, 
 the case when $(\Gamma, h)$ is immersive requires 
 essentially no change.
When $(\Gamma, h)$ is not immersive, we also have
 a general result (\cite[Theorem 57]{N1}).
We do not cite the full statement of it here, but we note the following.
Namely, in such general cases too,
 the dual obstruction space is described by the data attaching suitable vectors
 to the flags contained in 
 the loop part of $(\Gamma, h)$.
However, contrary to Theorem \ref{thm:obstruction},
 the space $H$ cannot be determined from the combinatorial data alone, 
 but depends on the configuration of corresponding pre-log curves. 
In general it is difficult to accurately determine the space $H$, 
 but in our genus one case, it is still possible.
We study this in Subsection \ref{subsec:general_genus_one}
 after studying immersive cases using Theorem \ref{thm:obstruction}. 

In view of Theorem \ref{thm:obstruction}
 and the above paragraph, we give the following definition.
\begin{defn}\label{def:abundancysupport}
	Let $(\Gamma, h)$ be a tropical curve satisfying Assumption A.	
	Then the \emph{support of superabundancy} of $(\Gamma, h)$ is
	the closed subgraph $\Gamma_s$ of $\Gamma$ such that 
	for any edge $E$ of $\Gamma_s$, there is an element of 
	$H$ such that the value of it on the flags associated to $E$
	is not zero.
\end{defn}
Part of the results which are valid for superabundant tropical curves of genus one
 can be extended to those tropical curves whose support of superabundancy
 is a loop.
See Remark \ref{rem:mainext}.

\section{Correspondence theorem for superabundant tropical curves  of genus one I:  Existence of smoothings for immersive  tropical curves}\label{sec:I}
Having described the dual space $H$ of obstructions,
 we now try to find a criterion describing the condition under which
 these obstructions 
 are effective or not.
When a tropical curve is immersive,
 the deformation of it is governed by 
 $H^0(C_{0}, \varphi^*_{0}\Theta_{\mathfrak X/\Bbb C}
  /\Theta_{C_{0}/O_{0}})$ regardless of whether it is superabundant or not
  (see  
    \cite[Subsection 5.2]{N1}),
   and there is no need to consider the obstruction.
Here $\varphi_0\colon C_0\to X_0$ is a pre-log curve of type
 $(\Gamma, h)$.
However,  
 the corresponding algebraic curves (or pre-log curves)
 actually have obstructions to a smoothing in general, and the smoothability cannot 
 be determined just from the calculation of the cohomology group $H$.
We have to calculate the \emph{Kuranishi map},
 and this is what we do in the rest of the paper.
First, we give some general remarks on Kuranishi maps, to motivate the
 calculation in the following arguments.

\subsection{Remarks on Kuranishi maps}\label{subsec:Kuranishi}
It is now common to formulate abstract deformation theory in terms of
 differential graded Lie algebras (dgLa).
Deformation theory of a complex manifold or a pair of them 
 consisting of an ambient space and a submanifold
 is a typical example.
The Kuranishi map, which is a nonlinear map between suitable
 vector spaces whose zero set is the moduli space of the deformations, 
 can also be formulated in terms of dgLa (see for example \cite{Ma}).
However, such an abstract formulation contains
 the use of the Green's operator, 
 a solution to non-linear partial differential equations, 
 or some other object which plays the 
 role equivalent to the Green's operator.
This is usually very difficult to write down
 in concrete problems.

Our suggestion in this paper is that degeneration technique is useful in 
 calculating the Kuranishi map.
For this purpose, we adopt two points of view on obstructions.
Let us explain it in a simplified situation.
Let $X$ be a nonsingular complex variety and 
\[
i\colon Y\hookrightarrow X
\]
 be an inclusion
 of a nonsingular subvariety,
 and suppose we want to study the deformation of $i$ with fixed $X$ and $Y$.

The first point of view is the standard definition through \v{C}ech cohomology.
A first order deformation of $i$ is naturally identified with a section of 
 $H^0(Y, i^*TX)$.
Take coordinate neighborhoods of $X$ and $Y$ in analytic category
 so that the image of each neighborhood of $Y$
 is contained in some neighborhood of $X$ (and fix this assignment).
The image of $i$ is locally presented by a vector valued function on $Y$.
A section of $H^0(Y, i^*TX)$ gives a perturbation of 
 these functions on each coordinate
 neighborhood, with the infinitesimal parameter $t$.

Up to the first order of
 $t$, these perturbations on coordinate neighborhoods are
 compatible under the coordinate change
 (because we start from a section of the pull back of the
 tangent bundle).
However, for the second and higher orders, 
 these local perturbations do not necessarily coincide with each other on
 the intersections of the coordinate neighborhoods.
These gaps, inductively with respect to the order of $t$, 
 compose the cohomology classes of degree one.
In particular, we call the map from $H^0(Y, i^*TX)$ to the cohomology
 group $H^1(Y, i^*TX)$ 
 which is given by the terms of  
 order two with respect to $t$,
\[
\kappa^{(1)}\colon H^0(Y, i^*TX)\to H^1(Y, i^*TX),
\]
 the \emph{Kuranishi map of order one}.

If we take a section of $H^0(Y, i^*TX)$ 
 for which the image of $\kappa^{(1)}$ is trivial
 (in other words, at a point in the zero locus of $\kappa^{(1)}$), we can 
 perturb the first order deformation of $Y$
 (more precisely, that of the map $i$) in the second order of $t$
 on each coordinate neighborhood, so that
 they are now compatible up to the second order of $t$.

Let $i_1$ be a first order deformation of $i$ which allows a
 second order deformation.
Then the set of second order deformations
 of $i_1$ is (non canonically) parametrised by $H^0(Y, i^*TX)$. 
Using these second order deformations, 
 we again calculate the first cohomology,
 which is now induced from the terms of order three
 with respect to $t$.
Then we can define the Kuranishi map of order two for the map $i_1$:
\[
\kappa_{i_1}^{(2)}\colon H^0(Y, i^*TX) \to H^1(Y, i^*TX).
\]

Note that although we are now considering the lift $i_1$ of $i$, 
 the Kuranishi map for $i_1$ still has the same domain and target as those for $i$. 
Continuing this process, one has a sequence $\{i_j\}_{0\leq j\leq k}$
 of deformations of $i = i_0$ such that:
\begin{itemize}
\item The map $i_j$ is defined over $\Bbb C[t]/t^{j+1}$.
\item The map $i_{j+1}$ reduces to $i_{j}$ over $\Bbb C[t]/t^{j+1}$.
\end{itemize}
Then, to the map $i_{k-1}$, we can define the Kuranishi map of 
 order $k$:
\[
\kappa_{i_{k-1}}^{(k)}\colon H^0(Y, i^*TX) \to H^1(Y, i^*TX).
\]
If this map contains the zero in its image, then there is a $k$-th order deformation of $i_{k-1}$
 which can be further deformed to the $(k+1)$-th order 
 (it is not necessarily the given $i_k$).
 
If we can continue this process infinitely many times, it defines
 a formal deformation of $i$.
Then in a favorable situation (when $X$ and $Y$ are projective, for example)
 an implicit function theorem \cite{A}
 assures the existence of an actual deformation.

We can put this construction in somewhat different way
 from another point of view.
Namely, let us take a section of $H^0(Y, i^*TX)$ again.
On a coordinate neighborhood, it gives a deformation of $i$ as before,
 which is represented by a vector valued function.
Then we analytically continue this function to $Y$.
In the first order of $t$, it restores the section of $H^0(Y, i^*TX)$.
In higher orders, in general it will not be possible to define the 
 analytic continuation on the whole $Y$, 
 resulting in a divergence or multivaluedness.
When we start from a section of $H^0(Y, i^*TX)$ 
 which is a solution of $\kappa^{(1)}=0$,
 it is possible to perturb the deformation on the coordinate neighborhood,
 so that it can be analytically
 continued to $Y$ up to the second order of $t$.
In view of this construction, 
 we may be able to calculate the obstruction 
 by looking at how the analytic continuation of
 a given local deformation fails to be global.

Both points of view are useful.
In fact, we use the \v{C}ech type representation of the obstruction classes
 in Subsection \ref{subsec:step3} and the 
 analytic continuation type representation of them in Subsections \ref{subsec:step4}
 and \ref{subsec:step5}.
The precise equivalence of these points of view is proved in Subsection \ref{subsec:pairing}.

In general, it is difficult to perform either of these constructions. 
However, in our setting, using degeneration, it becomes quite manageable.
Namely:
\begin{itemize}
\item We can find a suitable covering on which local lifts can be 
 calculated. 
\item We can also calculate the analytic continuation explicitly.
\item As a result, we can identify the obstruction for deforming degenerate 
 curves with the data of
 poles appearing in the analytic continuations of suitable functions.
This allows us to define the Kuranishi map in an explicit way.
\item In the case of genus one, 
 the Kuranishi map satisfies suitable transversality, 
 and we can qualitatively describe the solution space.
\end{itemize}
We note that the actual 
 setting is a little different from the one explained above.
First, we consider a family of maps
\[
\mathfrak C\to \mathfrak X
\]
 over $Spf\Bbb C[[t]]$, where $\mathfrak C$ and $\mathfrak X$
 are not necessarily smooth (but log smooth).
Here $\mathfrak C$ is a family of prestable curves and $\mathfrak X$ is 
 a toric degeneration of a toric variety $X$ (precisely, the completion of it
 along the central fiber).
Second, the domain of the Kuranishi map is not the 
 zeroth cohomology of the tangent sheaf, but the parameter space of 
 pre-log curves.
However, this parameter space is locally modelled on the 
 zeroth cohomology of a sheaf similar to $i^*TX$ 
 (precisely, the log normal sheaf of the map on the central fiber.
Here, in contrast to the above discussion,
 the tangent sheaf is replaced by the normal sheaf 
 to take the deformation of the domain curve into account).  
Therefore, essentially the argument proceeds parallel to the above simple setting.

In this section, the analytic continuations of local lifts are calculated.
This culminates in the calculation of the pre-obstruction,
 which is the function theoretic version of the obstruction class
 (Proposition
 \ref{prop:obstruction1}).
It is proved to be equivalent to the standard obstruction class 
 defined through \v{C}ech cohomological
 construction (Corollary \ref{cor:anal-ob}).
Therefore, it is essentially the same as the Kuranishi map.
Using this, we define the Kuranishi map
 (Definition \ref{def:preob} and Remark \ref{rem:kuranishi}).
By the transversality property
 of the Kuranishi map, the necessary and sufficient condition 
 for the smoothability of a tropical curve can be deduced (Theorems \ref{thm:immersive},
 \ref{thm:general}).
More detailed study of the Kuranishi map 
 will be done in the subsequent paper.
\subsection{Conventions on tropical curves}\label{subsec:convention}
Let $(\Gamma, h)$ be a 3-valent 
 superabundant tropical curve of genus one in $N_{\Bbb R}\cong \Bbb R^n$
 satisfying Assumption A.
Moreover, we assume $(\Gamma, h)$ is immersive in this section.
We impose the following natural conditions on
 all the tropical curves from now on:
\begin{itemize}
\item The direction vectors of the edges of $\Gamma$ span $\Bbb R^n$.
\item $(\Gamma, h)$ is defined over $\Bbb Z$.
\end{itemize}
If the first condition is not satisfied, then we can
 take an affine subspace of $N_{\Bbb R}$ so that this condition is
  satisfied.
Tropical curves which are related to toric varieties should be defined over
 rational numbers.
By multiplying suitable integers, we can 
 make such tropical curves to satisfy the
 second condition (see also Remark \ref{rem:basechange}).
 
Let $X$ be a toric variety associated to $(\Gamma, h)$
 (see Definition \ref{toric})
 and $\mathfrak X\to\Bbb C$ a degeneration of $X$ defined 
 respecting $(\Gamma, h)$ (see Definition \ref{def:degeneration}).
Let $\varphi_0\colon C_0\to X_0$ be a pre-log curve 
 of type $(\Gamma, h)$. 
The space $\mathfrak X$ and the base $\Bbb C$ have
 natural structures of log schemes induced from toric charts,
 see \cite{KF, KK}.
We also put a log structure on $C_0$ 
 so that $\varphi_0$ can be equipped with
 a structure of a morphism between log schemes
 with the following properties.
Namely, its composition with the projection 
 $\mathfrak X\to\Bbb C$ is a log smooth map
 and it is strict on the complement of the special points of $C_0$
 (that is, the inverse image of the toric divisors of components of $X_0$).
We write the log morphism also by $\varphi_0$.
Such log structures on $C_0$ are classified in \cite[Proposition 7.1]{NS}.

Let $L$ be the loop of $\Gamma$.
Let $A$ be the minimal dimensional affine subspace of $\Bbb R^n$ which contains
 $h(L)$, and let $\bar A$ be the subspace parallel to $A$.
The subset $A\cap h(\Gamma)$ of $h(\Gamma)$
 may have several connected components, 
 and let $\mathcal G$ be the unique component containing $h(L)$.
Let
 $\overline\Gamma'$  
 be the connected component of $h^{-1}(\mathcal G)$
 which contains the loop of $\Gamma$.

Because $(\Gamma, h)$ is superabundant, $\mathcal G$ necessarily has
 1-valent vertices when $h(\Gamma)$ does not have higher ($\geq 4$)
 valent vertices (those cases where there are
 such vertices are treated in Section \ref{sec:II}).
Let $\{ \alpha_i\}$ be the set of these 1-valent vertices.
These are in 
 one to one correspondence with the  1-valent vertices of $\overline{\Gamma}'$
 when $h$ is an immersion.
We write by $\Gamma'$ the open subset of
 $\overline\Gamma'$ with these 1-valent vertices deleted. 
 (see Figure \ref{fig:Gamma'}).
When $h$ is an embedding, 
 then $\Gamma'$ can be identified with $\mathcal G\setminus \{\alpha_i\}$. 

\begin{rem}\label{rem:2-val}
In general the graph $\Gamma$ must contain 2-valent vertices
 to satisfy (iv) of Assumption A.
However, since the components of pre-log curves corresponding to 
 these vertices play rather minor role (see the proof of \cite[Proposition 7.1]{NS}, 
 for example), we ignore these vertices
 in many cases to simplify the exposition.
\end{rem}
\subsection{Tropical curve with two vertices}\label{subsec:weight one}
\emph{In the rest of Section \ref{sec:I}, we assume that $h$ is an embedding.}
The case when $h$ is an immersion requires essentially no modification, 
 but we work with embeddings for notational simplicity.
General cases where we allow contractions of edges
 will be studied in Subsections \ref{subsec:general_genus_one}
 and \ref{subsec:4-vloop}.
Therefore, we identify $\Gamma$ and the image $h(\Gamma)$.
We can also identify $\mathcal G$ and $\overline\Gamma'$.
In particular, all the vertices of $h(\Gamma)$ are 3-valent
 and the weight of each edge of $h(\Gamma)$ is a single integer
 (see Definition \ref{def:trop}).

Intuitively, using the notation of Subsection \ref{subsec:convention},
 if we remove $\{ \alpha_i\}$ from $\mathcal G$ and extend the open edges 
 to infinity, we have 
 a tropical curve in the affine plane $A$ which is regular.
Thus, there is no obstruction to the smoothing
 of pre-log curves corresponding to it.
Therefore, roughly speaking,
 the first possible
 obstructions appear when we try to extend the local smoothing of
 the node corresponding to the edge of $\mathcal G$
 attached to some $\{ \alpha_i\}$ to the whole curve,
 from the point of view of analytic continuation 
 discussed in Section \ref{subsec:Kuranishi}.
In the following, we will make this observation precise.

\begin{figure}
\includegraphics{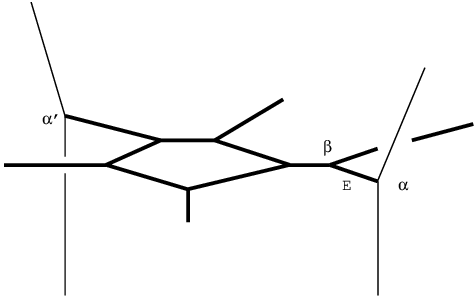}
\caption{The part drawn by bold lines is $h(\Gamma')$}\label{fig:Gamma'}
\end{figure}

Now, take a vertex $\alpha$ from  $\{ \alpha_i\}$
 and let $E$ be the edge of $\overline\Gamma'$ attached to $\alpha$.
Note that there is a unique path from $\alpha$ to the loop $L$.
Assume that the integral length of this path is the shortest among 
 the 1-valent vertices of $\overline\Gamma'$
 (more precisely, we need to use Definition \ref{def:path length}
 for the length of the path).
Let $\beta$ be the other vertex of $E$, see Figure \ref{fig:Gamma'}.

Let $t$ be the pull back to $\mathfrak X$ of the standard coordinate on the
 base space of the family $\mathfrak X\to \Bbb C$. 
Let $\varphi_0\colon C_0\to X_0$ be a pre-log curve of type $(\Gamma, h)$.
We note that such a curve always exists.
\begin{lem}\label{lem:immpre-logexist}
Let $(\Gamma, h)$ be an immersive tropical curve of genus one.
Let $\mathfrak X$ be a toric degeneration defined respecting $(\Gamma ,h)$
 and let $X_0$ be its central fiber.
Then there is a pre-log curve $\varphi_0\colon C_0\to X_0$ of type $(\Gamma, h)$. 
\end{lem}
\proof
We use the notation in Subsection \ref{subsec:convention}.
The open subgraph $h(\Gamma')$ of $h(\Gamma)$ 
 can be seen as
 a tropical curve in $\bar A$
 by extending the open edges to infinity. 
The resulting tropical curve is regular in the sense of 
 Definition \ref{def:nonsuperabundant2}.
 It is not difficult to see that for regular tropical curves, there always exist
 corresponding pre-log curves (in fact, it is a part of the claim of 
 Theorem \ref{thm:regsm}).
 Since the complement of $h(\Gamma')$ in $h(\Gamma)$ is a tree, 
 it is easy to extend such a pre-log curve corresponding to 
  $h(\Gamma')$ to a pre-log curve corresponding to $h(\Gamma)$.\qed\\

In the argument below, tropical curves with two vertices play
 an important role.
The reason for this is that the most part of 
 the corresponding pre-log curve is contained
 in one coordinate neighborhood of a toric degeneration $\mathfrak X$
 defined respecting the tropical curve.
Such a property is important for the calculation of the obstruction because
 as we mentioned in Subsection \ref{subsec:Kuranishi}, 
 for explicit calculation of the obstruction, representation of the 
 subvarieties (curves in our case) in terms of coordinates is very useful.
In view of this observation, 
 we introduce the following basic condition 
 for tropical curves with only two vertices.
\begin{defn}\label{def:twovertex}
Let $\Gamma_s$ be a 3-valent graph with two vertices 
 $\alpha_0$ and $\beta_0$
 such that $\alpha_0$ and $\beta_0$ are connected by an edge $E_0$, and
 the remaining edges are noncompact edges.
Let $(\Gamma_s, h_s)$ be a tropical curve in 
 $N_{\Bbb R}\cong \Bbb R^n$, $n\geq 3$, 
 which is an immersion and satisfies the following conditions.
\begin{itemize}
\item There is a basis of $N$ such that 
 the directions of the
 edges emanating from $\beta_0$ are spanned by the vectors
\[
(1, 0, 0, 0, \dots, 0), \;\; (0, 1, 0, 0, \dots, 0),\;\; (-1, -1, 0, 0, \dots, 0),
\]
 where the direction of the
 edge $E_0$ is spanned by $(1, 0, 0, 0, \dots, 0)$.
Similarly, the directions of the 
 edges emanating from $\alpha_0$ are spanned by the vectors
\[
(-1, 0, 0, 0, \dots, 0),\;\; (0, 0, 1, 0, \dots, 0),\;\; (1, 0, -1, 0, \dots, 0).
\] 
\item All the edge weights are 1.
\end{itemize}
We call such a tropical curve \emph{standard}.
\end{defn}
When we fix the rank of the lattice $N$, 
 standard tropical curves are parametrized by the integral length
 of the edge $h_s(E_0)$ up to isomorphisms.
Although standard tropical curves are rather special among the set of 
 tropical curves with two vertices, 
 they are of fundamental importance due to the following observations.
\begin{lem}\label{lem:twovertex}
Let $(\Xi, h_{\Xi})$ be
 an embedded 3-valent tropical curve in $N_{\Bbb R}$ with two vertices
 such that the integral length of the unique bounded edge is an integer
 multiple of the weight of that edge.
Then $(\Xi, h_{\Xi})$
 is represented as the image of a 
 standard tropical curve $\Gamma_s$
 by an integral affine endomorphism $\mathcal L$ of $N_{\Bbb R}$.
\end{lem} 
\proof
Let $\alpha$ and $\beta$ be the vertices of $\Xi$,
 and $E$ be the edge connecting $\alpha$ and $\beta$.
Let $E_1, E_2$ be the other edges emanating from $\beta$
 and $E_1', E_2'$ be the other edges emanating from $\alpha$.
We write by $w_E$ the weight of $E$ and by $r_E$ the integral length
 of the image $h_{\Xi}(E)$.
Similarly, let $w_{E_1}, w_{E_2}$ and $w_{E_1'}, w_{E_2'}$ be the
 weights of the corresponding edges.
  
We take $(\Gamma_s, h_s)$ so that the integral length of 
 the image $h_s(E_0)$ is $\frac{r_E}{w_E}$.
Let 
 $u$, $v_1$ and $v_2$ be the primitive integral vectors which
 generate the directions of the edges 
 $h_{\Xi}(E), h_{\Xi}(E_1)$ and $h_{\Xi}(E_2)$
 emanating from $h_{\Xi}(\beta)$.
Similarly, let $-u$, $v_1'$ and $v_2'$ be the 
 primitive integral vectors which
 generate the directions of the edges 
 $h_{\Xi}(E), h_{\Xi}(E_1')$ and $h_{\Xi}(E_2')$
 emanating from $h_{\Xi}(\alpha)$.

Let $\bar A$ be the subspace of $N_{\Bbb R}$ parallel to 
 the minimal affine subspace containing $h(\Xi)$.
The dimension of $\bar A$ is two or three.
We fix a basis 
 $\{f_1, \cdots, f_n\}$ of $N$ so that the set of vectors
 $\{f_1, f_2, f_3\}$ or $\{f_1, f_2\}$ is a basis of $A\cap N$.
 
Define the integral affine map $\mathcal L\colon N_{\Bbb R}\to N_{\Bbb R}$ by
 the following properties.
\begin{itemize}
\item $\mathcal L$ maps the vertex $h_{s}(\alpha_0)$
 to $h_{\Xi}(\alpha)$, and
 $\mathcal L$ is linear when $h_{s}(\alpha_0)$ and 
 $h_{\Xi}(\alpha)$ are regarded as the origins.
Let $\overline{\mathcal L}$ be the linear transformation obtained in this way.
\item $\overline{\mathcal L}$ transforms the vector $(0, 1, 0, \dots, 0)$ to 
 $w_{E_1}v_1$, 
 $(0, 0, 1, 0, \dots, 0)$ to $w_{E_1'}v_1'$ and
 $(1, 0, \dots, 0)$ to $w_Eu$.
\item $\overline{\mathcal L}$ transforms the other standard generators
 $e_i$, $i = 4, \cdots, n$ of $N$ to $f_i$.
\end{itemize}
Then by the balancing condition, we see that the vectors 
 $(-1, -1, 0, \dots, 0)$ and $(1, 0, -1, 0, \dots, 0)$ 
 are mapped to the vectors
 $w_{E_2}v_2$ and $w_{E_2'}v_2'$, respectively.
Moreover, the integral length of the image of the edge $h_s(E_0)$
 is $r_E$.\qed\\

Another property of standard tropical curves $(\Gamma_s, h_s)$
 concerns with the description of pre-log curves of type
 $(\Gamma_s, h_s)$.
Let $(\Gamma, h)$ be a superabundant embedded tropical curve of
 genus one and 
 take the vertices $\alpha$ and $\beta$ and the edge $E$ connecting them
 given at the beginning part 
 of this subsection (see Figure \ref{fig:Gamma'}).
Assume the restriction of $h$ to the union of the
 edges emanating from $\alpha$
 and $\beta$ gives a standard tropical curve (after extending the 
 edges other than $h(E)$ to infinity in the obvious way).

Let $\varphi_0\colon C_0\to \mathfrak X$ be a general pre-log curve of 
 type $(\Gamma, h)$.
Let $C_{0, \alpha}$ and $C_{0, \beta}$ be
 the components of $C_0$ corresponding
 to the vertices $\alpha$ and $\beta$.
Let $p = C_{0, \alpha}\cap C_{0, \beta}$ be the node corresponding
 to the edge $E$.
With this notation, it is easy to see the following.
\begin{lem}\label{lem:coordinate}
We can take a set of functions
 $\{ x, y, z, w_1, \dots, w_{n-2}\}$ on 
 $\mathfrak X$ around $\varphi_0(p)$ 
 and coordinates $s$ and $\sigma$ on the branches $C_{0, \alpha}$
 and $C_{0, \beta}$ of $C_{0}$  
 around the node $p$ with the following properties.
\begin{enumerate}[(i)]
\item Each of the functions in
 $\{ x, y, z, w_1, \dots, w_{n-2}\}$ is a character of the big torus acting
 on $\mathfrak X$.
\item The equation $xy = t^r$ holds, here $r$ is the integral length of the edge $h(E)$.
Each of the equations $x= 0$, $y=0$
 determines an irreducible component of $X_0$ around $\varphi_0(p)$,
 corresponding to the vertices $\alpha$ and $\beta$. 
Let us write the reduced structure of the variety corresponding to
 $\alpha$
 by $X_{0, \alpha}$ and assume it is (set theoretically) locally given by $y=0$.
Similarly, write by 
 $X_{0, \beta}$ the component of $X_0$ corresponding
 to the vertex $\beta$ and assume it is locally given by $x=0$.
\item The set of functions $\{x, z, w_1, \dots\}$ is
 a local coordinate system of $X_{0, \alpha}$
 around $\varphi_0(p)$.
Similarly, the set of functions $\{y, z, w_1, \dots\}$ is
 a local coordinate system of $X_{0, \beta}$ around $\varphi_0(p)$.
\item $\varphi_0^*(x) = s$, $\varphi_0^*(y) = \sigma$.
\item The component $C_{0, \alpha}$
 is mapped to $X_{0, \alpha}$ by $\varphi_0$.
The defining equations of the image are given by
\[
kx+lz+m= 0,\;\; y = 0, \;\;w_1 = -\frac{\mu}{\lambda},
 \;\;w_2 = a_2, \dots, w_{n-2}  =a_{n-2}.
\]
Here $k, l, m, \mu, \lambda, a_i$
 are generic complex numbers (in particular, nonzero).
\item Similarly, the component
 $C_{0,\beta}$ is mapped to $X_{0, \beta}$ and
 the defining equations of the image are given by
\[
\kappa y+\lambda w_1+\mu
  = 0, \;\;x = 0, \;\; z = -\frac{m}{l}, \;\; w_2 = a_2, \dots, w_{n-2} = a_{n-2}.
\]
Here $\kappa$ is a generic complex number.\qed
\end{enumerate}
\end{lem}
\begin{rem}\label{rem:coord}
\begin{enumerate}
\item
If $r$ is larger than one, then the set of functions 
 $\{x, y, z, w_1, \dots, w_{n-2}\}$ in the lemma is not a 
 coordinate system because $\mathfrak X$ has singularity
 around $\varphi_0(p)$.
In fact, if one wants a coordinate system, then one can 
 add 2-valent vertices to
 edges of $\Gamma$ so that every edge of $h(\Gamma)$ has 
 integral length one.
Then the set of functions 
 $\{x, y, z, w_1, \dots, w_{n-2}\}$ will be a coordinate system
 (note that our relevant edge $E$ has weight one by the assumption
 that $h(\Gamma)$ is standard around a neighborhood of the edge $E$).
 
This quite technical remark does not play an essential role in the following argument, 
 but it affects the normalization of log tangent vectors.
Namely, as we see below, there is a relation $\frac{dx}{x}+\frac{dy}{y} = r\frac{dt}{t}$
 between log cotangent vectors.
By this relation, when we define a basis of a log tangent space by taking 
 the dual 
 of the basis of a log cotangent space
 constructed from the functions
 $\{x, y, z, w_1, \dots, w_{n-2}\}$ above,
 the values of the bases vectors evaluated by the covector
 $\frac{dt}{t}$ will depend on $r$.
When we set $r = 1$, or more generally set $r$ equal to the weight of the edge, 
 such an ambiguity vanishes.   
This will be convenient for our later calculation, and we will usually assume this normalization
 afterwards.
\item
The coefficients $k, l, m$ etc. depend on the choice of
 the functions $x, y, z$ etc..
We will write obstructions using these coefficients, but eventually 
 the dependence on the choice vanishes.
See Remarks \ref{rem:coeff} and \ref{rem:dependence}.
\end{enumerate}
\end{rem}
Around $\varphi_0(p)$, the log cotangent sheaf $\Omega_{\mathfrak X}$
 (which is free of rank equal to $\dim\mathfrak X=n+1$)
 of $\mathfrak X$ is spanned by
\[
\frac{dx}{x},\;\; \frac{dy}{y},\;\; \frac{dz}{z},\;\; \frac{dw_1}{w_1}, \dots,\;\; \frac{dw_{n-2}}{w_{n-2}}.
\]
Let
\[
x\partial_x, \;\;y\partial_y, \;\;z\partial_z, \;\;w_1\partial_{w_1}, \dots, 
 w_{n-2}\partial_{w_{n-2}}
\]
 be its dual basis of the log tangent sheaf $\Theta_{\mathfrak X}$.
On the other hand, we have the identity
\[
\frac{dx}{x}+\frac{dy}{y} = r\frac{dt}{t}
\]
 of log cotangent vectors
 which restricts to  
\[
\frac{dx}{x}+\frac{dy}{y} = 0
\] 
 on $X_0$.
Following Remark \ref{rem:coord} (1), we assume the integral length of the $E$ 
 is one, that is, $r = 1$.

On the image $\varphi_0(C_{0, \alpha})$, we have identities
\[
kdx+ldz = 0,\;\; dy = dw_1 = \cdots = dw_{n-2} = 0.
\]
Similarly, on the image $\varphi_0(C_{0,\beta})$, we have
\[
\kappa dy+\lambda dw_1 = 0,\;\; dx = dz = dw_2 = \cdots = dw_{n-2}  = 0.
\]
From these, we see that, 
 around $\varphi_0(p)$, the push forward by $\varphi_0$ of the
 fibers of the log 
 tangent sheaf of $C_{0, \alpha}$
 (recall that $\varphi_0$ is a morphism of
 log schemes)
 are spanned by the vectors
\[
x\partial_x-y\partial_y-\frac{kx}{lz}\cdot z\partial_z
\]
 as a subsheaf of the log tangent sheaf $\Theta_{\mathfrak X}$
 of $\mathfrak X$. 
Similarly, the fibers of the image of 
 the log tangent sheaf of the union $C_{0,\alpha}\cup
 C_{0,\beta}$
 are spanned by
\[
x\partial_x-y\partial_y-\frac{kx}{lz}\cdot z\partial_z
 +\frac{\kappa y}{\lambda w_1}\cdot w_1\partial_{w_1}
\]
 around $\varphi_0(p)$.

What is relevant to us is the log normal sheaf, which
 is the quotient of the sheaf $\varphi_0^*\Theta_{\mathfrak X}$
 by the log tangent sheaf of the curve described above.
We will see the explicit local correspondence between the 
 sections of the log normal sheaf and the first order lifts
 of the map $\varphi_0$.
Before doing this, we study what changes are needed when 
 the local structure of $(\Gamma, h)$
 around the edge $E$ is not standard.
\subsubsection{Local calculation for non-standard tropical curves
 with two vertices.}\label{subsec:general_edge}
\begin{figure}[h]
\includegraphics{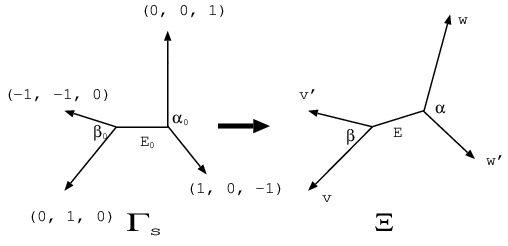}
\caption{}\label{fig:transform}
\end{figure}

As we saw in Lemma \ref{lem:twovertex},
 as long as $(\Gamma, h)$ is an embedding, 
 a neighborhood of the edge $E$ of Figure \ref{fig:Gamma'}
 is obtained as a result of 
 an integral affine map applied to
 a standard tropical curve,
 see Figure \ref{fig:transform} (here we identify the graphs and
 their images).
We write by $(\Xi, h_{\Xi})$
 the tropical curve with two vertices which is isomorphic to a
  neighborhood of the edge $E$ of $\Gamma$
 (with the edges other than $E$ extended to infinity).
First we assume the image $h(\Gamma)$ is not contained in 
 a two dimensional affine subspace of $N_{\Bbb R}$. 
The case when $h(\Gamma)$ is contained in 
 a two dimensional affine subspace is easier and we mention it later.

For notational simplicity, we assume $N\cong \Bbb Z^3$.
In Figure \ref{fig:transform}, the integral affine map 
 $\mathcal L\colon N\to N$ is 
 uniquely determined by the following properties:
\begin{itemize}
\item $\mathcal L$ maps the vertex $\alpha_0$ to $\alpha$, and
 $\mathcal L$ is linear when $\alpha_0$
 and $\alpha$ are regarded as the origins.
Let $\overline{\mathcal L}$ be the linear transformation obtained in this way.
\item $\overline{\mathcal L}$ transforms the vector $(0, 1, 0)$ to $v$, 
 $(0, 0, 1)$ to $w$ and
 $(1, 0, 0)$ to the vector $u$ spanning the edge $E$
 in the direction from $\beta$ to $\alpha$.
\end{itemize}
The linear map $\overline{\mathcal L}$ gives a map
\[
\Phi_{\mathcal L}\colon \Bbb P^3 \to \Bbb P,
\]
 where $\Bbb P$ is the toric variety defined by 
 the complete fan whose rays are spanned by
 the vectors $v, w, v'$ and $w'$.
Here $v' = \overline{\mathcal L}((-1,-1,0))$ and
 $w' = \overline{\mathcal L}((1,0,-1))$
 are vectors along one of the edges emanating from the vertices $\alpha$
 and $\beta$, respectively (see Figure \ref{fig:transform}).
 
Note that the vectors $u$, $v$ and $w$ (as well as $v'$, $w'$) may not be primitive.
They are multiples of primitive vectors by their edge weights.
The integral length of the edge $E$ is the weight of it times
 the integral length of the edge $E_0$.
\begin{rem}\label{coordinate}
In particular, the length of the edge $E$ is an integer multiple of 
 its weight, which is the condition required 
 in \cite{NS} when we apply the log deformation theory.
\end{rem}
Let $N^{\vee}$ be the dual of $N$ and 
\[
f_1, f_2, f_3
\]
 be the dual basis of the basis
 $(1, 0, 0), (0, 1, 0), (0, 0, 1)$ of $N$.
The standard tropical curve $\Gamma_s$ corresponds to 
 a linear curve in $\Bbb P^3$.
Consider a toric degeneration $\mathfrak X_s\to\Bbb C$ of 
 $\Bbb P^3$ defined respecting $\Gamma_s$.
As we described in the Lemma \ref{lem:coordinate}, there are components 
 $X_{s, 0, \alpha_0}, X_{s, 0, \beta_0}$ 
 of the central fiber $X_{s, 0}$ of $\mathfrak X_s$,
 corresponding to the vertices $\alpha_0, \beta_0$.
Moreover, there are coordinate systems
 $\{x, w, z\}$ on $X_{s, 0, \alpha_0}$ and 
 $\{y, w, z\}$ on $X_{s, 0, \beta_0}$
 corresponding to the vectors $-f_1, f_2, f_3$ and $f_1, f_2, f_3$
 by which a general pre-log curve of type $\Gamma_s$
 can be written
 in the form shown in the 
 properties (iv) and (v) of Lemma \ref{lem:coordinate}.
\begin{rem}
Precisely, the functions
 $x, y, z, w$ are restrictions of functions on $\mathfrak X$
 to $X_0$.
Therefore, these functions correspond to vectors in 
 the dual of $N\oplus \Bbb Z$, 
 where the generator of the dual of the summand $\Bbb Z$ 
 corresponds to the function
 $t$, the pull back of the coordinate on the base space of the degeneration.
The vectors corresponding to $x, w, z$ 
 project to $-f_1, f_2, f_3$ respectively by the natural projection
 $(N\oplus \Bbb Z)^{\vee} = N^{\vee}\oplus \Bbb Z^{\vee}\to N^{\vee}$.
The case of the functions $y, w, z$ is similar.
See also Remark \ref{rem:function}.
\end{rem}

On the other hand, let 
\[
g_1, g_2, g_3
\]
 be the dual basis in $N_{\Bbb Q}^{\vee}$
 of the vectors $u, v, w$ in $N$.
These are not in the lattice $N^{\vee}$ in general, but their suitable integral
 multiples are in $N^{\vee}$.
Let  
\[
\tilde g_1, \tilde g_2, \tilde g_3
\]
 be the positive integer multiples of $g_1, g_2, g_3$, which are primitive
 in $N^{\vee}$.
However, these may not compose a basis of $N^{\vee}$ in general.
A basis of $N^{\vee}$ is obtained from $\tilde g_3$ by choosing
 two other suitable vectors.
Let us choose one of them from $u^{\perp}\cap N^{\vee}$
 and write it by $\tilde g_2'$.
It is a linear combination of $\tilde g_2$ and $\tilde g_3$ over $\Bbb Q$.
Moreover, choose $\tilde g_1'$ so that the pairing
 $(\tilde g_1', u)$ is positive, and that the triplet
\[
\tilde g_1', \tilde g_2', \tilde g_3
\]
 becomes a basis of $N^{\vee}$.

The tropical curve $(\Xi, h_{\Xi})$ corresponds to a curve 
 in the toric variety $X = \Bbb P$
 described above.
Consider a toric degeneration $\mathfrak X\to\Bbb C$ of 
 $\Bbb P$ defined respecting $(\Xi, h_{\Xi})$
 and take a pre-log curve of type $\Xi$ on the central fiber $X_0$
 of $\mathfrak X$.
As in the case of $\mathfrak X_s$, 
 there are components $X_{0, \alpha}, X_{0, \beta}$
 of the central fiber
 $X_0$ corresponding to the vertices $\alpha, \beta$ of $\Xi$.
The component of the pre-log curve
 corresponding to the vertex $\beta$
 is contained in the closure of an orbit of the subtorus generated by
 the vectors $u$ and $v$.
The function corresponding to $\tilde g_3$ is constant on that orbit.


There are functions $\{Y, W, Z\}$ on $\mathfrak X$
 corresponding to the vectors $\tilde g_1', \tilde g_2', \tilde g_3$
 which restrict to a coordinate system on $X_{0, \beta}$
 in a neighborhood of the interior of $X_{0,\beta}\cap X_{0,\alpha}$
 (see Remark \ref{rem:function}).
Here interior means the complement of the
 toric divisors.
Also, we can take similar coordinates
 $\{X, W, Z\}$ on $X_{0, \alpha}$ 
 so that the functions $X$ and $Y$ satisfy the relation
\[
XY = t^{r}, 
\]
 where $r$ is the integral length of the edge $E$ of $\Xi$.
We will assume that $r$ equals to the weight of $E$ in view of Remark \ref{rem:coord} (1).

As we noted above, the linear map $\overline{\mathcal L}$
 induces a map $\Phi_{\mathcal L}$ from 
 $\Bbb P^3$ to $\Bbb P$, which is a branched covering map.
Any torically transverse
 curve corresponding to $\Xi$ in $\Bbb P$ is obtained as the image of
 a linear curve in $\Bbb P^3$ by $\Phi_{\mathcal L}$.
In fact, there are $\det\overline{\mathcal L}$
 different linear curves in $\Bbb P^3$ which map to 
 such a curve in $\Bbb P$.

It is easy to see that we can take the degeneration $\mathfrak X$
 so that there is a natural map from $\mathfrak X_s$
 to $\mathfrak X$ over the base $\Bbb C$.
In fact, the fan corresponding to $\mathfrak X_s$ 
 lies in 
 the space $\Bbb R^3\oplus\Bbb R$, and applying the 
 linear transformation $\overline{\mathcal L}\oplus Id_{\Bbb R}$,
 we obtain the fan corresponding to $\mathfrak X$.
We write this map by $\widetilde\Phi_{\mathcal L}$.
Moreover, any pre-log curve of type $\Xi$ in $X_0$
 is obtained as the image of a pre-log curve of type $\Gamma_s$
 in $X_{s, 0}$
 by the induced morphism.
In other words, a pre-log curve $\varphi_0\colon C_0\to X_0$ of type $\Xi$
 factors via a pre-log curve $\varphi_{s, 0}\colon C_0\to X_{s, 0}$
 of type $\Gamma_s$,
 so that $\varphi_0 = \widetilde\Phi_{\mathcal L}\circ\varphi_{s, 0}$.

The upshot is that although it is not easy in general to find 
 defining equations of the curve $\varphi_0(C_0)$ in $X_0$, 
 we can pull it back to $X_{s, 0}$ 
 where it is easy to find defining equations
 so that explicit calculation becomes possible. 

Let
\[
kx+lz+m = 0,\;\; y = 0,\;\; w = -\frac{\mu}{\lambda},
\]
and
\[
\kappa y+\lambda w+\mu = 0,\;\; x = 0,\;\; z = -\frac{m}{l}
\]
 be the defining 
 equations of the images of the components $C_{0, \alpha}$
 and $C_{0, \beta}$ of
 a pre-log curve of type $\Gamma_s$, respectively (see Lemma 
 \ref{lem:coordinate}).
As we calculated above, the
 log tangent sheaf of the image is spanned by the vector
\[
x\partial_x-y\partial_y-\frac{kx}{lz}\cdot z\partial_z
 +\frac{\kappa y}{\lambda w}\cdot w\partial_{w}.
\]
 near the node corresponding to the edge $E_0$.
The log tangent sheaf
 of a pre-log curve $\varphi_0$ of type $\Xi$
 is spanned by the push-forward of these vectors by 
 $\widetilde\Phi_{\mathcal L}$.

Important for us will be 
 the functions $X, Y, Z, W$ (particularly the function $Z$),
 and they are pulled back by 
 $\widetilde\Phi_{\mathcal L}$ to monomials of $x, y, z, w$
 (see Case (a) and (b) in Subsection \ref{subsec:step4} below).
They depend on how we choose $\tilde g_1'$ and $\tilde g_2'$,
 but as for the function $Z$, 
 since it corresponds to $\tilde g_3$,
 it is pulled back to a power of $z$:
\[
\widetilde\Phi_{\mathcal L}^*(Z) = z^c
\]
 for some integer $c$, up to a nonzero constant multiple.

Finally, when the image $h(\Gamma)$ is contained in a one or two dimensional 
 affine subspace, 
 the above construction works almost
 unchanged so that the curve $\varphi_0(C_0)$
 is obtained as the image of a standard curve which can be explicitly described
 by defining equations.
In this case too, the calculation of the restriction of functions on $\mathfrak X$
 to $\varphi_0(C_0)$ (or its deformations) will be important later.
This can be done by pulling them back to the standard curve again.


\noindent
\subsection{Calculation of obstructions in embedded cases. Step 1: Local sections of the log normal sheaf}\label{subsec:step1}
As we noted in Subsection \ref{subsec:Kuranishi}, 
 if $\varphi\colon Y\to X$ is an embedding between smooth varieties, 
 the first order deformations of $\varphi$ naturally correspond to 
 the sections of the pull back 
 $\varphi^*TX$ of the tangent bundle of $X$
 (when we allow the deformations of $Y$, then we use 
 the sections of the normal bundle $\varphi^*TX/TY$ instead).
Though our varieties are singular,
 since we deal with the log smooth situation, this principle is still valid.
Namely, the first order deformations correspond to 
 the sections of the log normal sheaf 
 $\mathcal N_{C_0/\mathfrak X} = 
  \varphi_0^*(\Theta_{\mathfrak X})/\Theta_{C_0}$.
In this subsection we describe sections of this sheaf.
In the following, 
 for notational simplicity
 we use the notation in Lemma \ref{lem:coordinate}, 
 which assumes that $(\Gamma, h)$ is standard in the neighborhood of
 the edge $E$ (see Definition \ref{def:twovertex}), and the integral length $r$ of the edge
 $E$ is one.
However, the arguments in this subsection require little change when 
 this assumption is removed (but we still assume that 
 $(\Gamma, h)$ is an embedding), see Remark \ref{rem:nonstdmod} below.

First, we describe the sheaf $\mathcal N_{C_0/\mathfrak X}$.
As we noted in Subsection \ref{subsec:abundance},  
 since $\mathfrak X$ is a toric variety, it has a natural log structure and
 the log tangent sheaf
 is free of rank $n+1$.
Using the functions in Lemma \ref{lem:coordinate},
 the generators are given by
\[
x\partial_x, \;\;y\partial_y, \;\;z\partial_z, \;\;w_1\partial_{w_1}, \dots, 
 w_{n-2}\partial_{w_{n-2}}.
\]
Also note that the log tangent sheaf of the central fiber $X_0$ of $\mathfrak X$ 
 is generated by
\[
x\partial_x-y\partial_y, \;\;z\partial_z,\;\; w_1\partial_{w_1}, \dots, w_{n-2}\partial_{w_{n-2}}.
\]

On the other hand, 
 the image $\varphi_0(C_{0, \alpha})$ of the component $C_{0, \alpha}$
 of $C_0$ corresponding to the vertex $\alpha$ of $\Gamma$
 is given by the defining equations
\[
kx+lz+m= 0,\;\; y = 0, \;\;w_1 = -\frac{\mu}{\lambda},
 \;\;w_2 = a_2, \dots, w_{n-2}  =a_{n-2}.
\]
As we calculated before,
 fibers of 
 the log tangent sheaf of $C_{0, \alpha}$ 
 near the node $\varphi_0(p)$ corresponding to the edge connecting the
 vertices $\alpha$ and $\beta$
 are generated by the vectors
\[
lz\cdot(x\partial_x-y\partial_y) -kx\cdot z\partial_z
\]
 as a subsheaf of $\varphi_0^*(\Theta_{\mathfrak X})$.
 
Thus, the restriction of the log normal sheaf 
 $\mathcal N_{C_0/\mathfrak X}$ to the 
 (suitable affine open subset of) component $C_{0, \alpha}$
 is given by
\[
\mathcal O_{C_{0, \alpha}}
 \langle x\partial_x, \;\;y\partial_y, 
 \;\;z\partial_z, \;\;w_1\partial_{w_1}, \dots, 
 w_{n-2}\partial_{w_{n-2}}\rangle
  /(lz\cdot(x\partial_x-y\partial_y) -kx\cdot z\partial_z, )
\]
 near $\varphi_0(p)$
 (precisely, the functions $lz$ and $kx$ should be pulled back
 to $C_{0, \alpha}$, but we use the same letters for notational simplicity).

Moreover, the restriction of the
 sheaf $\mathcal N_{C_0/\mathfrak X}$
 to $C_{0, \alpha}$ is isomorphic to
\[
\mathcal O_{C_{0, \alpha}}(1)\oplus\mathcal O_{C_{0,\alpha}}^{\oplus n-1}, 
\] 
 and the space of sections is generated by the images of the sections
\[
x\partial_x, \;\;y\partial_y, 
 \;\;z\partial_z, \;\;w_1\partial_{w_1}, \dots, 
 w_{n-2}\partial_{w_{n-2}}
\]
 of $\varphi_0^*\Theta_{\mathfrak X}$.
 
As we 
 noted above, the first order lifts of $\varphi_0$ naturally correspond to
 the sections of $\mathcal N_{C_0/\mathfrak X}$.
However, since 
 we consider smoothings over the base space $\Bbb C$ (precisely, 
 the infinitesimal neighborhood of the origin in $\Bbb C$), 
 the sections should be evaluated to one by the covector $\frac{dt}{t}$.

Recalling the relation
\[
\frac{dx}{x}+\frac{dy}{y} =\frac{dt}{t},
\]
 the following is an easy conclusion of this description.
\begin{lem}\label{lem:lift1.5}
The sections of the sheaf $\mathcal N_{C_0/\mathfrak X}$
 restricted to $C_{0, \alpha}$, 
 which are evaluated to one by the covector
 $\frac{dt}{t}$,  can be written in the form
\[\begin{array}{l}
\mathfrak n_{0, \alpha} = 
 y\partial_y + c(x\partial_x-y\partial_y) + c_0z\partial_z + 
 c_1w_1\partial_{w_1}+\cdots +c_{n-2}w_{n-2}\partial_{w_{n-2}}\\
 \hspace{3in}({\rm mod}\;\; 
  lz(x\partial_x-y\partial_y)-kx\cdot z\partial_z)
 \end{array}
\]
 \text{\rm around $\varphi_0(p)$}.
Here $c, c_0, c_1, \dots, c_{n-2}$ are arbitrary constants.
In particular, these $\mathfrak n_{0, \alpha}$
 compose an $n$ dimensional affine space.\qed
\end{lem}

Similarly, on $C_{0, \beta}$, whose image by $\varphi_0$ is given by
\[
\kappa y+\lambda w_1+\mu
  = 0, \;\;x = 0, \;\; z = -\frac{m}{l}, \;\; w_2 = a_2, \dots, w_{n-2} = a_{n-2},
\]
 such sections are represented by
\[\begin{array}{l}
\mathfrak n_{0, \beta} = 
 x\partial_x + c'(x\partial_x-y\partial_y) + c'_0z\partial_z + 
 c'_1w_1\partial_{w_1}+\cdots +c'_{n-2}w_{n-2}\partial_{w_{n-2}}\\
  \hspace{2.7in}({\rm mod}\;\; \lambda w_1(x\partial_x-y\partial_y)
   + \kappa y\cdot w_1\partial_{w_1})
 \end{array}
\]
 \text{\rm around $\varphi_0(p)$},
 with $c', c'_0, c'_1, \dots, c'_{n-2}$ arbitrary constants.

For later use, we also give the version of these descriptions
 for the union $C_{0, \alpha}\cup C_{0, \beta}$.
\begin{lem}\label{lem:lift2}
The set of sections of the sheaf
 $\mathcal N_{C_0/\mathfrak X}$ restricted to
 $C_{0, \alpha}\cup C_{0, \beta}$ 
 which are evaluated to one by the 
 covector $\frac{dt}{t}$ is given by the set
\[
\{(\mathfrak n_{0, \alpha}, \mathfrak n_{0, \beta})\,|\, c_i = c'_i, \, i = 0, 1, \dots, n-2\},
\] 
 where $\mathfrak n_{0, \alpha}$ and $\mathfrak n_{0, \beta}$ are the sections of 
 $\mathcal N_{C_0/\mathfrak X}$ restricted to $C_{0,\alpha}$
 and $C_{0,\beta}$ respectively, and $c_i$ and $c_i'$ are their coefficients
 in the above expression (note that $c$ and $c'$ can be different).
In particular, these sections compose an $n+1$ dimensional affine space.
\end{lem}
\proof
The sections $\mathfrak n_{0, \alpha}$ and 
 $\mathfrak n_{0, \beta}$ on $C_{0, \alpha}$ and $C_{0, \beta}$
 glue and give a section on $C_{0,\alpha}\cup C_{0,\beta}$ if and only if
 they coincide at the node $C_{0,\alpha}\cap C_{0,\beta}$.
Clearly it is necessary that $c_i$ equals to $c_i'$ for $i = 0, 1, \dots, n-2$.
On the other hand, at the node, we have
\[
x\partial_x-y\partial_y = \frac{kx}{lz}\cdot z\partial_z
 -\frac{\kappa y}{\lambda w_1}\cdot w_1\partial_{w_1} = 0
\]
 since $x = y = 0$ there.
Thus, there is no condition imposed on the coefficients of 
 the term $x\partial_x-y\partial_y$ of $\mathfrak n_{0, \alpha}$
  and $\mathfrak n_{0, \beta}$. \qed

\begin{rem}\label{rem:zeroorder}
From the view point of
 the correspondence between the set of sections of the log normal sheaf 
 and the set of lifts of $\varphi_0$ studied in the next subsection, 
 one should consider the section $\mathfrak n_{0, \alpha}$
 as
\[
\mathfrak n_{0, \alpha} = 
(y\partial_y -\frac{\kappa y}{\lambda w_1}\cdot w_1\partial_{w_1})
 + c(x\partial_x-y\partial_y) + c_0z\partial_z + 
c_1w_1\partial_{w_1}+\cdots +c_{n-2}w_{n-2}\partial_{w_{n-2}}
\]
 and similarly, 
\[
\mathfrak n_{0, \beta} = 
(x\partial_x -\frac{kx}{lz}\cdot z\partial_{z})
+ c'(x\partial_x-y\partial_y) + c_0'z\partial_z + 
c_1'w_1\partial_{w_1}+\cdots +c_{n-2}'w_{n-2}\partial_{w_{n-2}},
\]
 although the terms $\frac{\kappa y}{\lambda w_1}\cdot w_1\partial_{w_1}$ and
 $\frac{kx}{lz}\cdot z\partial_{z}$ vanish on $C_0$
 (so that $\mathfrak n_{0, \alpha}$ and $\mathfrak n_{0, \beta}$
 reduce to the expression in Lemma \ref{lem:lift2} above). 
See Remark \ref{rem:secorder2}.
\end{rem}

Finally, we note the modifications required when the tropical curve 
 is not locally isomorphic to the standard one.
\begin{rem}\label{rem:nonstdmod}
When the directions of the edges emanating from the vertices
 $\alpha$ and $\beta$ span a three dimensional subspace, 
 then we can obtain sections of the normal sheaf
 by pushing forward the sections of the normal sheaf in the standard case described above
 by $(\widetilde{\Phi}_{\mathcal L})_*$
 constructed in Subsection \ref{subsec:general_edge}.
When this subspace is two dimensional, we can again push forward
 by $(\widetilde{\Phi}_{\mathcal L})_*$, but
 the map $(\widetilde{\Phi}_{\mathcal L})_*$ has one dimensional cokernel,
 and it gives a trivial summand of $\mathcal N_{C_0/\mathfrak X}$
 restricted to $C_{0, \alpha}\cup C_{0, \beta}$.
\end{rem}

\subsection{Calculation of obstructions in embedded cases.  Step 2: The local lifts, first order and higher.}\label{subsec:lifthigh}
Consider a tropical curve $(\Gamma, h)$ as before, 
 and we use the same notation $\Gamma', \alpha, \beta, E$ etc.
  as in Subsections \ref{subsec:convention} and \ref{subsec:weight one}.
In the previous subsection, 
 we described the sections of the log normal sheaf restricted to one or two 
 components of $C_0$.
As we mentioned at the beginning of Subsection \ref{subsec:step1},
 there is a natural correspondence between
 the set of the first order deformations of
 $\varphi_0$ and the set of sections of the log normal sheaf.
In this subsection, we explicitly describe the local lifts corresponding to 
 local sections of the log normal sheaf.
In fact, when one restricts attention to a local coordinate neighborhood, 
 there is a one to one correspondence between local lifts and local sections
 even in higher order.
For notational simplicity, we again assume that the restriction of 
 $(\Gamma, h)$ to a neighbourhood of the edge $E$ is standard
 (see Definition \ref{def:twovertex}) and the integral length of the edge $E$
 is one.
The argument in this subsection can be applied to general immersive tropical curves
 by reducing the problem to the standard one by pulling back the curve
 as in Subsection \ref{subsec:general_edge}.

Let $\mathscr D$ be the set of all the nodes of the curve $C_0$.
\begin{defn}\label{def:C^*}
Let us write $C_{0, \alpha}^* = C_{0,\alpha}\setminus \mathscr D$
 with the induced log scheme structure
 and similarly define $C_{0, \beta}^*$.
These are rational curves with at most three points removed.
If $C_k$ is a $k$-th order deformation of $C_0$, then 
 we write by $C_{k, \alpha}^*$ the log scheme structure
 on the topological space underlying $C_{0, \alpha}^*$
 induced by the restriction.
\end{defn}
We first describe the spaces of lifts of $\varphi_0$
 restricted to $C_{0, \alpha}^*$ and 
 $C_{0, \beta}^*$.
Recall that the curve $\varphi_0(C_{0, \alpha})$
 is presented by the equations
\[
(\ast) \;\;\;\; kx+lz+m= 0, \;\; y = 0,  \;\;w_1 = -\frac{\mu}{\lambda},\;\;
 w_2 = a_2, \;\; \dots, \;\; w_{n-2}  =a_{n-2}.
\]
\begin{defn}\label{def:D_alpha}
Let $\mathcal D_{\alpha}$ be the set of isomorphism classes of 
 stable maps from 
 $C_{0, \alpha}^*\times Spec\Bbb C[t]/t^2$ to $\mathfrak X$
 (precisely, to the restriction of  $\mathfrak X$ over $\Bbb C[t]/t^2$.
 We neglect this point hereafter)
 log smooth over $\Bbb C[t]/t^2$ which specialize to the restriction
 $\varphi_0|_{C_{0,\alpha}^*}$ of $\varphi_0$
 on the central fiber.
\end{defn}
Let $\Gamma(C_{0, \alpha}^*, \mathcal N_{C_0/\mathfrak X})$
 be the space of sections of $\mathcal N_{C_0/\mathfrak X}$ restricted to 
 $C_{0, \alpha}^*$ which are evaluated to one by $\frac{dt}{t}$.
A section $\mathfrak n_{0, \alpha}$ 
 in $\Gamma(C_{0, \alpha}^*, \mathcal N_{C_0/\mathfrak X})$
 can be described in the form similar to the one in 
 Lemma \ref{lem:lift1.5},
\[
\mathfrak n_{0, \alpha} =
 y\partial_y + c(s)(x\partial_x-y\partial_y) + c_0(s)z\partial_z + 
 c_1(s)w_1\partial_{w_1}+\cdots +c_{n-2}(s)w_{n-2}\partial_{w_{n-2}},
\]
 here $s$ is a parameter on $C_{0, \alpha}^*$, and 
 now the coefficients $c, c_0, \dots, c_{n-2}$ are rational functions of $s$
 which can have poles of any order at the punctures of $C_{0, \alpha}^*$.

Now the curve $\varphi_0(C_{0, \alpha})$ is parametrized as
\[
(x, y, z, w_1, w_2, \dots, w_{n-2})
= (s, 0, -\frac{k}{l}s-\frac{m}{l}, -\frac{\mu}{\lambda}, a_2, \dots, a_{n-2}),
\]
 using a parameter $s$ on $C_{0, \alpha}$which is zero at
 the intersection with $x = 0$ (we may take $s$ to be the pullback of $x$
 by $\varphi_0$).
Note that for a smoothing $C_k$ of $C_0$, the restriction
of $C_k$ to the topological space underlying 
$C_{0, \alpha}^*$ is isomorphic to the product
$C_{0, \alpha}^*\times Spec\Bbb C[t]/t^{k+1}$, since $C_{0, \alpha}^*$ is
 a rational curve with at most three punctures, and there is no nontrivial 
 flat deformations.
Thus, if $s$ is a parameter of $C_{0, \alpha}^*$ 
$($defined over $\Bbb C[t]/t$$)$,
we can regard $s$ also as a parameter of $C_{k, \alpha}^*$ 
$($defined over $\Bbb C[t]/t^{k+1}$$)$.
The following is a consequence of direct calculation.
\begin{lem}\label{lem:lift3}
There is a natural one to one correspondence 
 between the space $\Gamma(C_{0, \alpha}^*, \mathcal N_{C_0/\mathfrak X})$
 and the set $\mathcal D_{\alpha}$.
Namely, 
 a section $\mathfrak n_{0, \alpha}= 
 y\partial_y + c(s)(x\partial_x-y\partial_y) + c_0(s)z\partial_z + 
 c_1(s)w_1\partial_{w_1}+\cdots +c_{n-2}(s)w_{n-2}\partial_{w_{n-2}}$ 
 corresponds to a lift of $\varphi_0$ parametrized as
\[\begin{array}{l}
(x, y, z, w_1, w_2, \dots, w_{n-2})  \\
 \hspace{.1in}= ((1+c(s))s, \frac{t}{s}, (1+c_0(s))(-\frac{k}{l}s-\frac{m}{l}),
(1+c_1(s))(-\frac{\mu}{\lambda}), (1+c_2(s))a_2, \dots, (1+c_{n-2}(s))a_{n-2}).
\end{array}
\]
\qed
\end{lem}
\begin{rem}\label{rem:secorder}
Note that the parameterization
 of the first order lift of $\varphi_0|_{C_{0, \alpha}}$
 given in Lemma \ref{lem:lift3} should reduce to the parameterization of 
 $\varphi_0|_{C_{0, \alpha}}$ over $\Bbb C[t]/t$.
Therefore, the coefficients $c(s), c_0(s), \dots, c_{n-2}(s)$
 should have values in $t\cdot\mathcal O_{C_{0, \alpha}^*}$.
Namely, when we consider the natural 
 correspondence between the set of sections
 of the log normal sheaf and the set of lifts of $\varphi_0$ as in Lemma \ref{lem:lift3},
 it is more legitimate to think that
 the section $\mathfrak n_{0, \alpha}$ belongs not to 
 $\Gamma(C_{0, \alpha}^*, \mathcal N_{C_0/\mathfrak X})$, but to
 the space $y\partial_y +
  \Gamma(C_{0, \alpha}^*, \mathcal N_{C_0/X_0}) \otimes_{\Bbb C} t\Bbb C[t]/t^2$.
This is a quite technical point, but plays a role in the calculation of 
 obstructions, see Subsection \ref{subsec:resobs}.
\end{rem}
In general, it is not easy to describe the lifted curve by defining equations.
However, there is one distinguished lift.
Namely, since the image $\varphi_0(C_{0, \alpha}^*)$ is contained in 
 a single coordinate neighborhood of $\mathfrak X$, 
 the defining equations ($\ast$) of $\varphi_0(C_{0, \alpha}^*)$ above
 (restricted to the locus $x\neq 0, z\neq 0$),
 which are equations defined over $\Bbb C[t]/t$, can be regarded as 
 equations defined over $\Bbb C[t]/t^2$, or even over $\Bbb C[t]/t^{k+1}$
 for an arbitrary positive integer $k$, 
 only by replacing $xy = 0$ with $xy = t$.
These equations define the lift of $\varphi_0$ parametrized as
\[
(x, y, z, w_1, w_2, \dots, w_{n-2})  \\
\hspace{.1in}= (s, \frac{t}{s}, -\frac{k}{l}s-\frac{m}{l},
-\frac{\mu}{\lambda}, a_2, \dots, a_{n-2}).
\]
In other words, this is the $k$-th order
 lift corresponding to the section $\mathfrak n_{0, \alpha} = y\partial_y$.
\begin{defn}\label{def:trivlift1}
We call the $k$-th order lift corresponding to $\mathfrak n_{0, \alpha} = y\partial_y$
 the \emph{trivial lift} of $\varphi_0|_{C_{0, \alpha}^*}$.
\end{defn}
Using this, we can generalize the local correspondence between 
 the set of sections of the log normal sheaf
 and the set of lifts of $\varphi_0$ given in Lemma \ref{lem:lift3} to higher order
 as follows.
 
Let $\mathcal D_{k, \alpha}$ be the set of isomorphism classes of stable maps
 from $C_{0, \alpha}^*\times Spec\Bbb C[t]/t^{k+2}$ to $\mathfrak X$ log smooth over 
 $\Bbb C[t]/t^{k+2}$ which specialize to the restriction $\varphi_0|_{C_{0, \alpha}^*}$
 of $\varphi_0$ on the central fiber.

Also, consider the space 
 $y\partial_y+
 \Gamma(C_{0, \alpha}^*, \mathcal N_{C_0/X_0}) \otimes_{\Bbb C} t\Bbb C[t]/t^{k+1}$.
\begin{lem}\label{lem:lift3h}
There is a natural one to one correspondence between the set $\mathcal D_{k, \alpha}$
 and the set of sections in
 $y\partial_y+
  \Gamma(C_{0, \alpha}^*, \mathcal N_{C_0/X_0}) \otimes_{\Bbb C} t\Bbb C[t]/t^{k+2}$.
Explicitly, the lift corresponding to a section
 $y\partial_y + c(s)(x\partial_x-y\partial_y) + c_0(s)z\partial_z + 
  c_1(s)w_1\partial_{w_1}+\cdots +c_{n-2}(s)w_{n-2}\partial_{w_{n-2}}$ 
  is parametrized as
 \[\begin{array}{l}
 (x, y, z, w_1, w_2, \dots, w_{n-2})  \\
  \hspace{.1in}= ((1+c(s))s, \frac{1}{1+\bar c(s)}\frac{t}{s}, (1+c_0(s))(-\frac{k}{l}s-\frac{m}{l}),
 (1+c_1(s))(-\frac{\mu}{\lambda}), (1+c_2(s))a_2, \dots, (1+c_{n-2}(s))a_{n-2}),
 \end{array}
 \]
 where $c(s), c_0(s), \dots, c_{n-2}(s)\in 
  \Gamma(C_{0, \alpha}^*, \mathcal N_{C_0/X_0}) \otimes_{\Bbb C} t\Bbb C[t]/t^{k+2}$.
Also, $\bar c(s)$ is the part of $c(s)$
 whose coefficient is of order less than $k+1$ with respect to the exponent of $t$.
\end{lem}
\proof
This is proved by induction with respect to the order $k$.
See the proof of Lemma \ref{lem:lifthigh} below for details.\qed\\

Similarly, there is a natural one to one correspondence 
 between the set of isomorphism classes 
 of stable maps over $\Bbb C[t]/t^{k+2}$ which specialize to 
 $\varphi_0|_{C_{0,\beta}^*}$
 and the space of sections of
  $x\partial_x+
    \Gamma(C_{0, \beta}^*, \mathcal N_{C_0/X_0}) \otimes_{\Bbb C} t\Bbb C[t]/t^{k+2}$.\\

We can perform the same construction for larger subsets of $C_0$.
Let $\alpha$ and $\beta$ be the adjacent vertices of $\Gamma$ 
 as before.
\begin{defn}\label{def:twin}
Let $(C_{0,\alpha}\cup C_{0,\beta})^*$ be the open subset 
 of $C_0$ whose underlying topological space is
 $C_{0, \alpha}^*\cup C^*_{0,\beta}\cup
 (C_{0,\alpha}\cap C_{0,\beta})$ 
 with the log scheme structure induced from $C_0$.
If $C_k$ is a $k$-th order deformation of $C_0$, then we write by 
 $(C_{k,\alpha}\cup C_{k,\beta})^*$ the log scheme structure on 
 the topological space underlying 
 $(C_{0,\alpha}\cup C_{0,\beta})^*$ induced by restriction.
\end{defn}
The importance of these subsets lies in the following two properties:
\begin{itemize}
\item They give a covering of the curve $C_k$.
\item Each of them is mostly contained in a single coordinate neighborhood of 
 $\mathfrak X$ (when there is 
 a point corresponding to an unbounded edge
 of $\Gamma$, then $(C_{0,\alpha}\cup C_{0,\beta})^*$ is not
 contained in a single coordinate neighborhood of $\mathfrak X$.
But this does not affect our argument below, since the image of such a point is
 uniquely determined
 by its complement by the valuative criterion).
\end{itemize}
When the edge length $r$ of the edge between the vertices $\alpha$ and 
 $\beta$ is larger than one, then the single neighborhood mentioned in the second
 property above is not a coordinate neighborhood in the strict sense 
 (since $\mathfrak X$ has singularity). 
In this case, we add 2-valent vertices so that the refined
 edges have length one (see Remark \ref{rem:coord}).
These two facts allow us to calculate the obstruction cohomology class 
 using these open subsets. 
\begin{defn}\label{def:twoverticeslift}
Let $\mathcal D_{\alpha\cup\beta}$ be the set of 
 isomorphism classes of stable maps log smooth over 
 $\Bbb C[t]/t^2$ which specialize to the restriction
 $\varphi_0|_{(C_{0,\alpha}\cup C_{0,\beta})^*}$ on the central fiber.
\end{defn}
\begin{rem}\label{rem:twocomppara}
Note that $(C_{0, \alpha}\cup C_{0, \beta})^*$ is an affine log scheme
 log smooth over the standard log point whose underlying affine scheme is 
 given by the ring $\Bbb C[s, u, \frac{1}{s-1}, \frac{1}{u-1}]/(su)$.
Then as in the case of usual affine smooth schemes, 
 any two log deformations of $(C_{0, \alpha}\cup C_{0, \beta})^*$
 are mutually isomorphic.
In particular, over the base $\Bbb C[t]/t^{k+1}$, 
 the underlying scheme is given by the ring
 $\Bbb C[s, u, \frac{1}{s-1}, \frac{1}{u-1}, t]/(su-t, t^{k+1})$.
 
Therefore, for log deformations $(C_{k, \alpha}\cup C_{k, \beta})^*$
 of $(C_{0, \alpha}\cup C_{0, \beta})^*$
 of any order, we can use $s$ and $u$ as parameters, as in the irreducible 
 case (see the paragraph before Lemma \ref{lem:lift3}). 
\end{rem}

Let $\Gamma((C_{0, \alpha}\cup C_{0, \beta})^*, \mathcal N_{C_0/\mathfrak X})$
 be the space of sections of $\mathcal N_{C_0/\mathfrak X}$ restricted to 
 $(C_{0, \alpha}\cup C_{0, \beta})^*$ which are evaluated to one by $\frac{dt}{t}$.
A section $(\mathfrak n_{0, \alpha}, \mathfrak n_{0, \beta})$ in 
 $\Gamma((C_{0, \alpha}\cup C_{0, \beta})^*, \mathcal N_{C_0/\mathfrak X})$
 can be described as in Lemma \ref{lem:lift2}.
Precisely, let $s$ be a parameter on $C_{0, \alpha}$ which is zero at 
 the intersection with $x = 0$ as before, and let $u$ be a parameter
 on $C_{0, \beta}$ which is zero at the intersection with $y = 0$.
Then the coefficients $c, c_0, \dots, c_{n-2}$ defining 
 $\mathfrak n_{0,\alpha}$ are rational functions
 of $s$ which may have poles at the punctures of 
 $(C_{0,\alpha}\cup C_{0,\beta})^*\cap C_{0, \alpha}$.
Similarly, the coefficients $c', c'_0, \dots, c'_{n-2}$
 are rational functions of $u$ which may have poles at the punctures of 
 $(C_{0,\alpha}\cup C_{0,\beta})^*\cap C_{0, \beta}$.
Moreover, the constant terms of the pairs 
$(c_0, c_0'), \dots, (c_{n-1}, c_{n-2}')$ (that is, their values at the node 
 $C_{0, \alpha}\cap C_{0, \beta})$ must be the same.

The following is the analogue of Lemma \ref{lem:lift3}.
\begin{lem}\label{lem:lift4}
There is a natural one to one correspondence between the set
 $\mathcal D_{\alpha\cup\beta}$ and 
 the space 
 $\Gamma((C_{0,\alpha}\cup C_{0,\beta})^*,
  \mathcal N_{C_0/\mathfrak X})$. 
See Remark \ref{rem:secorder2} below for more explicit description of the correspondence. \qed
\end{lem}
Note that the image of the union $(C_{0, \alpha}\cup C_{0, \beta})^*$ by $\varphi_0$
 is defined by the equations 
\[
\begin{array}{l}
kx+lz+m = 0,\\
\kappa y+\lambda w_1+\mu = 0,\\
xy =  0,\\
 w_2 = a_2,\;\; \dots,\;\; w_{n-2} = a_{n-2}.\\
\end{array}
\]
Since the image of $(C_{0,\alpha}\cup C_{0,\beta})^*$ is mostly
 contained 
 in a single coordinate neighborhood of $\mathfrak X$
 (see the paragraph after Definition \ref{def:twin}), 
 there is a particular lift
 as in Definition \ref{def:trivlift1}, 
 whose image is given by the same set of equations 
 as $\varphi_0((C_{0, \alpha}\cup C_{0, \beta})^*)$
 except replacing the equation $xy = 0$ 
 by $xy = t$.
This corresponds to the section of $\mathcal N_{C_0/\mathfrak X}$
 on 
 $(C_{0, \alpha}\cup C_{0, \beta})^*$
 given by the pair
 $(\mathfrak n_{0, \alpha}, \mathfrak n_{0, \beta})$ where 
 all the coefficients $c, c', c_i = c'_i$, $i = 0, 1, \dots, n-2$ are zero
 in the notation of Lemma \ref{lem:lift2}.

\begin{rem}\label{rem:secorder2}
\begin{enumerate}
\item More explicitly, the lift in the above paragraph 
 corresponds to the pair
 $(\mathfrak n_{0, \alpha}, \mathfrak n_{0, \beta})$ with
\[
\mathfrak n_{0, \alpha} = y\partial_y-\frac{\kappa y}{\lambda w_1}\cdot w_1\partial_{w_1},
\;\;
\mathfrak n_{0, \beta} = x\partial_x-\frac{k x}{l z}\cdot z\partial_{z},
\]
 where we wrote the sections $\mathfrak n_{0, \alpha}$ and
 $\mathfrak n_{0, \beta}$ following Remark \ref{rem:zeroorder}.
The terms $\frac{\kappa y}{\lambda w_1}\cdot w_1\partial_{w_1}$
 and $\frac{k x}{l z}\cdot z\partial_{z}$ are zero over $\Bbb C[t]/t$, 
 but are written as  $\frac{\kappa t}{\lambda w_1x}\cdot w_1\partial_{w_1}$
 and $\frac{k t}{l zy}\cdot z\partial_{z}$ using the relation $xy = t$,
 and in fact relevant to the lift.
Namely, the lift is parametrized as
\[
(x, y, z, w_1, \dots, w_{n-2}) =
 (s, \frac{t}{s}, -\frac{k}{l}s-\frac{m}{l}, 
 -\frac{\kappa}{\lambda}\frac{t}{s}-\frac{\mu}{\lambda},
 a_2, \dots, a_{n-2}),
\]
 using the parameter $s$ extended to the lift of 
 $(C_{0, \alpha}\cup C_{0, \beta})^*$, see Remark \ref{rem:twocomppara}.
This can also be written by using the parameter $u$ by
 replacing $s$ with $\frac{t}{u}$.

\item As in Remark \ref{rem:secorder}, under the correspondence 
 in Lemma \ref{lem:lift4}, we should regard the components of
 the section 
 $(\mathfrak n_{\alpha, 0}, \mathfrak n_{0, \beta})$ to be in the set
\[
\mathfrak n_{0, \alpha}\in (y\partial_y-\frac{\kappa y}{\lambda w_1}\cdot w_1\partial_{w_1})
   + t\Bbb C[t]/t^{2}\otimes_{\Bbb C}
 \Gamma({(C_{0,\alpha}\cup C_{0,\beta})^*\cap C_{0, \alpha}}, \mathcal N_{C_0/X_0})
 \]
 and 
\[
\mathfrak n_{0, \beta}\in
 (x\partial_x-\frac{k x}{l z}\cdot z\partial_{z})
  + t\Bbb C[t]/t^{2}\otimes_{\Bbb C}
 \Gamma({(C_{0,\alpha}\cup C_{0,\beta})^*\cap C_{0, \beta}}, \mathcal N_{C_0/X_0})
\]
 whose constant terms satisfying the conditions in Lemma \ref{lem:lift2},
  though there is a one to one correspondence between the set of such sections and 
  the set $\Gamma((C_{0, \alpha}\cup C_{0, \beta})^*, \mathcal N_{C_0/\mathfrak X})$.
See Subsection \ref{subsec:resobs} for the role played by this remark.
See also Remark \ref{rem:globalsect} for a related note.
\item In the lemma, we assumed that the integral length of the edge $E$
 in Figure \ref{fig:Gamma'} is one, as we noted at the beginning of this subsection. 
When the integral length of $E$ is larger than one, then we add 2-valent vertices on $E$
 so that all the resulting refined edges have length one.
In this case, the image of the curve $C_{0, \beta}$ coincides with the closure of an orbit of a
 one dimensional torus acting on the component $X_{0, \beta}$, and the defining equation
 of the curve becomes
\[
\begin{array}{l}
kx+lz+m = 0,\\
xy =  0,\\
w_1 = -\frac{\mu}{\lambda},\;\; w_2 = a_2,\;\; \dots,\;\; w_{n-2} = a_{n-2}.\\
\end{array}
\]
The calculation for this case in the following arguments 
 is simpler than the case in the lemma, and 
 we usually omit further details for this case in the later arguments.
\end{enumerate}
\end{rem}

\begin{defn}\label{def:trivlift2}
We call the lift of $\varphi_0|_{(C_{0, \alpha}\cup C_{0, \beta})^*}$
 corresponding to the section 
 $(\mathfrak n_{0, \alpha}, \mathfrak n_{0, \beta})
  = (y\partial_y-\frac{\kappa y}{\lambda w_1}\cdot w_1\partial_{w_1},
  \; x\partial_x-\frac{k x}{l z}\cdot z\partial_{z})$ the \emph{trivial lift}.
\end{defn}

As before, this also gives a lift over $\Bbb C[t]/t^{k+1}$
 for an arbitrary positive integer $k$, and using this we can generalize the correspondence 
 in Lemma \ref{lem:lift4} to higher order.

Namely, let $\mathcal D_{k, \alpha\cup \beta}$ be the set of isomorphism
 classes of stable maps log smooth over $\Bbb C[t]/t^{k+2}$
 which specialize to the restriction $\varphi_0|_{(C_{0, \alpha}\cup C_{0, \beta})^*}$
 on the central fiber.

Also, consider the space of sections $(\mathfrak n_{k, \alpha}, \mathfrak n_{k, \beta})$
 where 
\[\begin{array}{ll}
\mathfrak n_{k, \alpha}
  & = (y\partial_y-\frac{\kappa y}{\lambda w_1}\cdot w_1\partial_{w_1})
      + c(s)(x\partial_x-y\partial_y) + c_0(s)z\partial_z + 
        c_1(s)w_1\partial_{w_1}+\cdots +c_{n-2}(s)w_{n-2}\partial_{w_{n-2}}
      \\
     & \in 
 (y\partial_y-\frac{\kappa y}{\lambda w_1}\cdot w_1\partial_{w_1})
    + t\Bbb C[t]/t^{k+2}\otimes_{\Bbb C}
  \Gamma({(C_{0,\alpha}\cup C_{0,\beta})^*\cap C_{0, \alpha}}, \mathcal N_{C_0/X_0})
 \end{array}
 \]
  and 
 \[\begin{array}{ll}
 \mathfrak n_{k, \beta}
   & = (x\partial_x-\frac{k x}{l z}\cdot z\partial_{z})+
   c'(u)(y\partial_y-x\partial_x) + c'_0(u)z\partial_z + 
     c'_1(u)w_1\partial_{w_1}+\cdots +c'_{n-2}(u)w_{n-2}\partial_{w_{n-2}}\\
    &\in 
   (x\partial_x-\frac{k x}{l z}\cdot z\partial_{z})
     + t\Bbb C[t]/t^{k+2}\otimes_{\Bbb C}
    \Gamma({(C_{0,\alpha}\cup C_{0,\beta})^*\cap C_{0, \beta}}, \mathcal N_{C_0/X_0}).
  \end{array}
  \] 
Here $s$ and $u$ are the parameters as in Remark \ref{rem:secorder2}.
Also, $c, c_0, c_1, \dots, c_{n-2}$ are rational functions of $s$ 
 which may have poles at the
 punctures of $(C_{0,\alpha}\cup C_{0,\beta})^*\cap C_{0, \alpha}$.
Similarly, $c', c'_0, \dots, c'_{n-2}$
  are rational functions of $u$ which may have poles at the punctures of 
  $(C_{0,\alpha}\cup C_{0,\beta})^*\cap C_{0, \beta}$.
Moreover, their constant terms satisfy
 $c_i(0) = c'_i(0)$ as elements of $t\Bbb C[t]/t^{k+2}$ for $0\leq i\leq n-2$
 (note that $s = 0$ and $u = 0$ corresponds to the node of 
 $(C_{0, \alpha}\cup C_{0, \beta})^*$ and $c, c_0, \dots, c_{n-2}$ as well as 
 $c', c_0', \dots, c'_{n-2}$ do not have poles at this node). 
We write this constant as $\gamma_i$.
\begin{lem}\label{lem:lifthigh}
There is a natural one to one correspondence between the set $\mathcal D_{k, \alpha\cup\beta}$
 and the set of sections $(\mathfrak n_{k, \alpha}, \mathfrak n_{k, \beta})$
 described above.
Explicitly, the lift corresponding to a section
 $(\mathfrak n_{k, \alpha}, \mathfrak n_{k, \beta})$ in the paragraph above 
  is parametrized as
 \[\begin{array}{l}
 (x, y, z, w_1, w_2, \dots, w_{n-2}) \\
= (\frac{1+c(s)}{1+\bar c'(u)}s, \frac{1+c'(u)}{1+\bar c(s)}u, 
 -\frac{k}{l}(s+\frac{m}{k})(1+c_0(s)+c_0'(u)-\gamma_0),
  -\frac{\kappa}{\lambda}(u
  + \frac{\mu}{\kappa})(1+c_1(s)+c_1'(u)-\gamma_1), \\
  \hspace{1.6in}
   (1+c_2(s)+c_2'(u)-\gamma_2)a_2, \dots, 
    (1+c_{n-2}(s)+c_{n-2}'(u)-\gamma_{n-2})a_{n-2}),
 \end{array}
 \]
 mod $t^{k+2}$, see Remark \ref{lem:paramtwo} below.
Here $\bar c(s)$ is the part of $c(s)\in t\Bbb C[t]/t^{k+2}\otimes
  \mathcal O_{(C_{0, \alpha}\cup C_{0, \beta})^*\cap C_{0, \alpha}}$
  whose coefficient is of order less than $k+1$ with respect to the exponent of $t$.
The function $\bar c'(u)$ is defined similarly.
\end{lem}
\begin{rem}\label{lem:paramtwo}
The precise meaning of the
 parameterization in the lemma is as follows.
Namely, on the component $C_{{k+1}, \alpha}^*$, the lift corresponding to
 $(\mathfrak n_{k, \alpha}, \mathfrak n_{k, \beta})$
 is parametrized by the above expression where $u$ is replaced by $\frac{t}{s}$
 and then terms with order higher than $k+2$ with respect to 
 $t$ are discarded.
Similarly, on $C_{{k+1}, \beta}^*$, the parameter $s$
 is replaced by $\frac{t}{u}$ and 
 then terms with order higher than $k+2$ with respect to 
 $t$ are discarded.
\end{rem}
\proof
We prove this by induction.
The case $k = 0$ corresponds to first order lifts and it corresponds to the
 statement of Lemma \ref{lem:lift4}.
They are obtained from the trivial lift given in Remark \ref{rem:secorder2} (1),
 by perturbations corresponding to sections of 
 $t\Bbb C[t]/t^{2}\otimes_{\Bbb C}
   \mathcal N_{C_0/X_0}|_{(C_{0,\alpha}\cup C_{0,\beta})^*\cap C_{0, \alpha}}$
 and 
 $t\Bbb C[t]/t^{2}\otimes_{\Bbb C}
     \mathcal N_{C_0/X_0}|_{(C_{0,\alpha}\cup C_{0,\beta})^*\cap C_{0, \beta}}$.
These sections can be regarded as 
 Lie algebra (of the torus acting on $\mathfrak X$, with coefficients in $t\Bbb C[t]/t^2$)
 valued functions (modulo the tangent vector), and act on the trivial lift by
 exponentiation.
Explicitly, for sections $c(s)(x\partial_x-y\partial_y) + c_0(s)z\partial_z + 
        c_1(s)w_1\partial_{w_1}+\cdots +c_{n-2}(s)w_{n-2}\partial_{w_{n-2}}$
        in $\Gamma({(C_{0,\alpha}\cup C_{0,\beta})^*\cap C_{0, \alpha}},
         t\Bbb C[t]/t^{2}\otimes_{\Bbb C}
        \mathcal N_{C_0/X_0})$
    and $c'(u)(y\partial_y-x\partial_x) + c'_0(u)z\partial_z + 
         c'_1(u)w_1\partial_{w_1}+\cdots +c'_{n-2}(u)w_{n-2}\partial_{w_{n-2}}$,
         in $\Gamma({(C_{0,\alpha}\cup C_{0,\beta})^*\cap C_{0, \beta}},
         t\Bbb C[t]/t^{2}\otimes_{\Bbb C}
         \mathcal N_{C_0/X_0})$,
   the result is
\[
\begin{array}{l}
 (x, y, z, w_1, w_2, \dots, w_{n-2}) \\
= ((1+c(s))s, (1+c'(u))u, 
 -\frac{k}{l}(s+\frac{m}{k})(1+c_0(s)+c_0'(u)-\gamma_0),
  -\frac{\kappa}{\lambda}(u
  + \frac{\mu}{\kappa})(1+c_1(s)+c_1'(u)-\gamma_1), \\
  \hspace{2in}
   (1+c_2(s)+c_2'(u)-\gamma_2)a_2, \dots, 
   (1+c_{n-2}(s)+c_{n-2}'(u)-\gamma_{n-2})a_{n-2}).
 \end{array}
\]

Now assume the lemma is proved up to a non-negative integer $k-1$.
Then we have a curve parametrized as in the expression in the statement of the lemma
 (defined over $\Bbb C[t]/t^{k+1}$),
 and it can be regarded as a curve over $\Bbb C[t]/t^{k+2}$ in an obvious way
 (namely, regard all the coefficients $c(s), s_0(s), \dots, c_{n-2}(s)$
 and $c'(u), c_0'(u), \dots, c_{n-2}'(u)$ over $\Bbb C[t]/t^{k+1}$
 as defined over $\Bbb C[t]/t^{n+2}$). 
The curve obtained in this way corresponds to the section 
 $(\mathfrak n_{k, \alpha}, \mathfrak n_{k, \beta})$
 without terms whose coefficients are in $t^{k+1}\Bbb C[t]/t^{k+2}$.
Other lifts of order $k+1$ (which restricts to the original curve over $\Bbb C[t]/t^{k+1}$)
 are obtained from this by perturbing by the action corresponding to sections
 $c(s)(x\partial_x-y\partial_y) + c_0(s)z\partial_z + 
         c_1(s)w_1\partial_{w_1}+\cdots +c_{n-2}(s)w_{n-2}\partial_{w_{n-2}}$
     and $c'(u)(y\partial_y-x\partial_x) + c'_0(u)z\partial_z + 
          c'_1(u)w_1\partial_{w_1}+\cdots +c'_{n-2}(u)w_{n-2}\partial_{w_{n-2}}$
  where all the coefficients $c, c_1, \dots, c'_{n-2}$ are of order $t^{k+1}$.
They act on the given curve of order $k+1$ by exponentiation as above, and
 the result coincides with the one given in the statement of the lemma.\qed

\begin{rem}\label{rem:param}
\begin{enumerate}
\item The correspondence given in the lemma does not depend on the coordinates of 
 $\mathfrak X$ as long as they satisfy the properties of Lemma \ref{lem:coordinate}
 (with $r = 1$).
This is because the correspondence is described using sections of 
 $\mathcal N_{C_0/X_0}$ and actions of them on curves through the 
 canonical torus action on $\mathfrak X$, and these do not depend on 
 the choice of coordinates.
\item Although the sections $\mathfrak n_{k, \alpha}$ and $\mathfrak n_{k, \beta}$
 belong to the spaces $(y\partial_y-\frac{\kappa y}{\lambda w_1}\cdot w_1\partial_{w_1})
 + t\Bbb C[t]/t^{k+2}\otimes_{\Bbb C}
 \Gamma({(C_{0,\alpha}\cup C_{0,\beta})^*\cap C_{0, \alpha}}, \mathcal N_{C_0/X_0})$
 and $(x\partial_x-\frac{k x}{l z}\cdot z\partial_{z})
 + t\Bbb C[t]/t^{k+2}\otimes_{\Bbb C}
 \Gamma({(C_{0,\alpha}\cup C_{0,\beta})^*\cap C_{0, \beta}}, \mathcal N_{C_0/X_0})$,
 respectively (therefore the parameters $s$ and $u$ are those for the curve
 $(C_{0, \alpha}\cup C_{0, \beta})^*$), 
 in the expression in Lemma \ref{lem:lifthigh},
 the parameters $s$ and $u$ are those for the deformed curve
 $(C_{k+1, \alpha}\cup C_{k+1, \beta})^*$ through Remark \ref{rem:twocomppara},
 since the expression in Lemma \ref{lem:lifthigh} is the parameterization of 
 the image of a map from $(C_{k+1, \alpha}\cup C_{k+1, \beta})^*$
 to $\mathfrak X$.
 
In other words, we can regard $\mathfrak n_{k, \alpha}$ and $\mathfrak n_{k, \beta}$
 as sections of the normal sheaf of the $(k+1)$-th order trivial lift of 
 $\varphi_0|_{(C_{0, \alpha}\cup C_{0, \beta})^*}$ so that the parameters
 $s$ and $u$ satisfy the relation $su = t$.
However, the lift corresponding to $(\mathfrak n_{k, \alpha}, \mathfrak n_{k, \beta})$
 in the lemma is not the one which is obtained from the $(k+1)$-th order trivial lift
 of $\varphi_0|_{(C_{0, \alpha}\cup C_{0, \beta})^*}$ by the action of the 
 exponential of $(\mathfrak n_{k, \alpha}, \mathfrak n_{k, \beta})$.
Instead, it is obtained by repeating the trivial lift 
 and perturbing by the action order by order.
\item In the expression of Lemma \ref{lem:lifthigh},
 the parameter $s$ is the one introduced in Remark \ref{rem:twocomppara}.
On the other hand, it is also convenient to use 
 $S= \frac{1+c(s)}{1+\bar c'(\frac{t}{s})}s$ as a parameter of the curve. 
This is the pull back of the function $x$ on $\mathfrak X$.
Note that $s$ can be expressed by $S$ in the form 
 $s = S(1 + t h(t, S))$, where $h(t, S)$ is a Laurent series in $S$ with coefficients 
 in $\Bbb C[t]/t^{k+2}$.
\end{enumerate}
\end{rem}

Before studying a global lift of $\varphi_0$, we make some remarks.
\begin{rem}\label{rem:noncan}
The trivial lift of $\varphi_0|_{C_{0, \alpha}^*}$
 or $\varphi_0|_{C_{0, \beta}^*}$ does not coincide with the 
 the restrictions of the trivial lift of $\varphi_0|_ {(C_{0, \alpha}\cup C_{0, \beta})^*}$
 to $C_{0, \alpha}^*$ or to $C_{0, \beta}^*$.
The notion of trivial lifts depends on the choice of a coordinate neighborhood.
\end{rem}
We note the following result, which is
 in contrast to the later calculation of obstruction classes.
We do not need the result itself later, 
 but we give this argument here since we use its ideas and notations 
 in the following subsections.
\begin{lem}\label{lem:first_order_deform}
	The map $\varphi_0\colon C_0\to \mathfrak X$ has a first order lift
	$\varphi_1\colon C_1\to \mathfrak X$.
	The set of these lifts is a torsor over the vector space of 
	global sections $\Gamma(C_0, \mathcal N_{C_0/X_0})$.
\end{lem}
\proof
Once we prove the existence of a lift, then the latter claim of the 
lemma is a standard fact in deformation theory.
For the existence of a lift, it suffices to prove that there is a global
section of the sheaf $\mathcal N_{C_0/\mathfrak X}$
which is evaluated to one by the covector $\frac{dt}{t}$
 (precisely, such a section which locally has
 the form given in Remark \ref{rem:secorder2}, 
 see also Remark \ref{rem:globalsect} below).
We can construct such a section as follows.

Recall that we write by $L$ the loop of $\Gamma$.
Let $C_{0, L}$ be the corresponding union of components of $C_0$.
Let $A$ be the minimal affine subspace in $N_{\Bbb R}$
 which contains the image $h(L)$ of $L$.
Let $\bar A$ be the vector subspace parallel to $A$.
Choose a basis 
\[
\{\mathfrak v_1, \dots, \mathfrak v_{n-r}, 
 \mathfrak w_{1},\dots, \mathfrak w_r\}
\]
 of $N$
 so that $\{\mathfrak v_1, \dots, \mathfrak v_{n-r}\}$
 is a basis of $\bar A\cap N$.
On $C_{0, L}$, the vectors 
 $\mathfrak w_{1},\dots, \mathfrak w_r$
 give rise to globally defined sections
 $s_{1}, \dots, s_r$ of 
 $\varphi_0|_{C_{0, L}}^*(\Theta_{\mathfrak X})$
 which do not have a zero locus.
Note that these sections are evaluated to zero by the covector
 $\frac{dt}{t}$.

We have correspondingly the decomposition of the log normal
 sheaf $\mathcal N_{C_{0, L}/\mathfrak X}$:
\[
\mathcal N_{C_{0, L}/\mathfrak X}
 \cong \mathcal N_{C_{0, L}/\mathfrak X}'
 \oplus \mathcal O_{C_{0, L}}\langle
  s_{1}, \dots, s_r\rangle.
\]
Here $\mathcal N_{C_{0, L}/\mathfrak X}'$ is isomorphic to 
 the log normal sheaf of the maximally degenerate curve
 associated to the tropical curve
\[
h_{\bar A}\colon \tilde L\to \bar A.
\]
Here $\tilde L$ is the open subgraph of $\Gamma$
 which is the union of $L$ and the adjacent (open) edges.
The map $h_{\bar A}$ is the natural map obtained from 
 the restriction $h|_{\tilde L}$ by extending the images of the
 open edges to infinity.
The map $h_{\bar A}$ clearly gives 
 the structure of a tropical curve to the graph $\tilde L$.
The associated
 maximally degenerate curve has the target space different from 
 $\mathfrak X$ and there is no canonical choice of the target space, 
 but the isomorphism class of its log normal sheaf 
 does not depend on the choice.
We write this sheaf by 
 $\mathcal N'_{C_{0, L}/\mathfrak X}$.

On the other hand, since the tropical curve 
 $h_{\bar A}\colon \tilde L\to \bar A$ is regular, 
 associated maximally degenerate curves admit smoothings
 (\cite{CFPU, N1, S, Tyo}).
In particular, 
 there is a global section of $\mathcal N_{C_{0, L}/\mathfrak X}'$
 which is evaluated to one by the covector $\frac{dt}{t}$.
This can be regarded as a section of $\mathcal N_{C_{0, L}/\mathfrak X}$
 by the above splitting.
 
The complement $\Gamma\setminus L$ is a union of trees.
Then by the description of local sections of the log normal sheaf 
 restricted to neighboring 
 two components of $C_0$ in Lemma \ref{lem:lift2}
 (and a similar result for more general cases, see Remark \ref{rem:nonstdmod}),
 it is easy to see that a given local section on the part 
 $C_{0, L}$ can be extended to the parts of $C_0$ corresponding to 
 these trees.
Thus, there is a global section of $\mathcal N_{C_0/\mathfrak X}$
 which is evaluated to one by the covector $\frac{dt}{t}$.
Therefore, there is a first order lift of $\varphi_0$, too. \qed \\

\begin{rem}\label{rem:globalsect}
Recall that the total space $\mathfrak X$ is a toric variety
 defined by the fan contained in $\Bbb R^n\times \Bbb R_{\geq 0}$.
In particular, the log tangent sheaf $\Theta_{\mathfrak X}$ 
 decomposes reflecting the product structure of
 $\Bbb R^n\times \Bbb R_{\geq 0}$.
When $\Theta_{\mathfrak X}$ is pulled back to the central fiber $X_0$
 by the inclusion, the part corresponding to the summand 
 $\Bbb R^n$ gives $\Theta_{X_0}$.
We write the decomposition as
\[
\Theta_{\mathfrak X}|_{X_0} \cong \Theta_{X_0}\oplus t\partial_t\cdot \mathcal O_{X_0},
\]
 here $t\partial_t$ is the log tangent vector corresponding to the 
 vector $(0, 1)\in \Bbb R^n\times \Bbb R_{\geq 0}$.
It gives a global section of $\Theta_{\mathfrak X}|_{X_0}$ which is evaluated to one
 by the covector $\frac{dt}{t}$.
Therefore, when it is restricted to the image of $C_0$ by $\varphi_0$, it gives a
 section of the log normal sheaf which is evaluated to one by $\frac{dt}{t}$,
 and thus gives a deformation of $\varphi_0$.
 
However, this section locally has the form $x\partial_x + \eta$, in terms of the
 expression as in Remark \ref{rem:secorder2}.
Here $\eta$ is a local section of $\mathcal N_{C_0/X_0}$ over $\Bbb C[t]/t$,
 not over $t\Bbb C[t]/t^2$.
Thus, when we reduce the deformed curve over $\Bbb C[t]/t^2$ to over
 $\Bbb C[t]/t$, then it does not in general recover the original pre-log curve, 
 but gives some other curve.
\end{rem}

As we can see in the proof of Lemma \ref{lem:first_order_deform},
 a first order global lift of $\varphi_0$
 on the whole $C_0$ is obtained by a simple gluing operation on 
 local lifts which are easier to construct.
However, we cannot unconditionally assume that such a gluing 
 process is always possible for higher order lifts.
In general, as we discussed in Subsection \ref{subsec:Kuranishi},
 local lifts given on suitable open coverings of $C_0$
 give rise to a cohomology class
 in $H^1(C_0, \mathcal N_{C_0/\mathfrak X})$,
 whose vanishing guarantees the 
 existence of a global lift.

If the tropical curve $(\Gamma, h)$ is regular, the 
 cohomology group $H^1(C_0, \mathcal N_{C_0/\mathfrak X})$
 is trivial.
Thus, the existence of a global lift is automatic 
 (at least when $h$ is an embedding).
However, in our superabundant case, the group 
 $H^1(C_0, \mathcal N_{C_0/\mathfrak X})$ is not trivial.
Therefore, we have to explicitly calculate the cohomology class.
In the following, 
 we carry out the calculation of the obstruction class. 
Since elements of 
 the dual of the group $H^1(C_0, \mathcal N_{C_0/\mathfrak X})$
 are supported on the loop (see Corollary \ref{cor:obstbasis}),
 essentially we only need to care what happens at the loop.
We formulate this point precisely in the next subsection.

\subsection{Calculation of obstructions in embedded cases.  Step 3: Decomposition of the obstruction.}\label{subsec:step3}
Assume we have constructed a $k$-th order lift 
\[
\varphi_k\colon C_k\to\mathfrak X
\]
 of $\varphi_0$.
Here $C_k$ is a suitable lift of $C_0$ over $\Bbb C[t]/t^{k+1}$.
The problem is to construct a $(k+1)$-th order lift of $\varphi_k$.
When $k = 0$, by Lemma \ref{lem:lift4}, 
 locally we have a trivial lift of $\varphi_0$ restricted to
 each $(C_{0, \alpha}\cup C_{0, \beta})^*$,
 where $\alpha$ and $\beta$ are adjacent vertices of $\Gamma$.
Moreover, there is an explicit correspondence between the set of
 first order lifts and the set of sections of the log normal sheaf.

We first note that locally we also 
 have trivial lifts of higher order lifts.
In fact essentially we have seen it in the proof of 
 Lemma \ref{lem:lifthigh}.
\begin{lem}\label{lem:canlift}
The restrictions $\varphi_k|_{C_{k, \alpha}^*}$
 and $\varphi_k|_{(C_{k,\alpha}\cup C_{k, \beta})^*}$
 have the trivial lifts
 $\varphi_{k+1}|_{C_{k+1, \alpha}^*}$ and 
 $\varphi_{k+1}|_{(C_{k+1,\alpha}\cup C_{k+1, \beta})^*}$
 over $\Bbb C[t]/t^{k+2}$.
Here $C_{k+1, \alpha}^*$ and $(C_{k+1,\alpha}\cup C_{k+1, \beta})^*$
 are suitable lifts of $C_{k, \alpha}^*$ and
 $(C_{k,\alpha}\cup C_{k, \beta})^*$, respectively.
\end{lem}
\proof
The proof is the same as those of Lemmas \ref{lem:lift3} and
 \ref{lem:lift4}. 
Namely, since the images of $\varphi_{k}|_{C_{k, \alpha}^*}$ and 
 $\varphi_{k}|_{(C_{k,\alpha}\cup C_{k, \beta})^*}$
 are mostly contained in one chart of $\mathfrak X$
 (see the paragraph after Definition \ref{def:twin}), we can regard the
 defining equations for them (defined over $\Bbb C[t]/t^{k+1}$)
 as equations over $\Bbb C[t]/t^{k+2}$.
The resulting curves give the lift.\qed

\begin{rem}
As in Remark \ref{rem:twocomppara}, the isomorphism classes of the lifts
 $C_{k+1, \alpha}^*$ and $(C_{k+1,\alpha}\cup C_{k+1, \beta})^*$
 are unique.
\end{rem}

Recall that as we remarked after Definition \ref{def:twin},
 the open subsets $(C_{k,\alpha}\cup C_{k, \beta})^*$
 cover the curve $C_k$.
As we argued in Subsection \ref{subsec:Kuranishi}, 
 the differences of local lifts compose the obstruction cohomology class
 in $H^1(C_0, \mathcal N_{C_0/\mathfrak X})$.
More precisely, if we have local $(k+1)$-th order lifts on $(C_{k,\alpha}\cup C_{k, \beta})^*$
 and $(C_{k,\beta}\cup C_{k, \gamma})^*$,
 then their difference 
 is converted to a local section
 of the log normal sheaf of $C_0$ restricted to $C_{0, \beta}^*$ 
 through Lemma \ref{lem:lift3h},
 by taking the part of order $k+1$ with respect to the exponent of $t$.
Therefore, using the trivial lifts
 $\varphi_{k+1}|_{(C_{k+1,\alpha}\cup C_{k+1, \beta})^*}$
 in Lemma \ref{lem:canlift}, we can calculate the obstruction
 class.
For later reference, we record this as a definition.
\begin{defn}\label{def:cechobst}
We call the cohomology class in $H^1(C_0, \mathcal N_{C_0/\mathfrak X})$
 defined in the above paragraph the \emph{obstruction class} to lift
 $\varphi_k$ one step further, and write it as $ob(\varphi_k)$.
\end{defn}

Note that the inclusion $X_0\to \mathfrak X$ induces 
 a map 
\[
\mathcal N_{C_0/X_0}\to \mathcal N_{C_0/\mathfrak X}
\]
 of sheaves on $C_0$.
Let $H^1(C_0, \mathcal N_{C_0/X_0})\to 
 H^1(C_0, \mathcal N_{C_0/\mathfrak X})$ 
 be the map induced on cohomology groups.
\begin{lem}\label{lem:obst_0}
The obstruction class belongs to 
 the image of $H^1(C_0, \mathcal N_{C_0/X_0})$
 by the above map.
\end{lem}
\proof
This is because we only consider
 the local lifts which are evaluated to one by
 the covector $\frac{dt}{t}$.
The set of local sections of $\mathcal N_{C_0/\mathfrak X}$ satisfying this
 condition is a torsor over that of $\mathcal N_{C_0/X_0}$.\qed\\

Recall that the dual of the group of obstruction classes 
 $H^1(C_0, \mathcal N_{C_0/X_0})^{\vee}$
 is described 
 in Corollary \ref{cor:obstbasis}.
Namely, an element of it is given by a $\mathcal N^{\vee}_{C_0/X_0}$-valued
 1-form under the Serre duality for nodal curves:
 $H^1(C_0, \mathcal N_{C_0/X_0})^{\vee}\cong
 	 H^0(C_0, \mathcal N_{C_0/X_0}^{\vee}\otimes\omega_{C_0})$.
Here $\omega_{C_0}$ is the dualizing sheaf of $C_0$, which is the sheaf of 1-forms
 with logarithmic poles allowed at the nodes of $C_0$
 (see \cite{N2}, Theorem 39 for more detail).

Moreover, under the natural isomorphism 
 $\mathcal N_{C_0/X_0}^{\vee}\cong N_{\Bbb C}^{\vee}\otimes\mathcal O_{C_0}$,
 an element of $H^1(C_0, \mathcal N_{C_0/X_0})^{\vee}$
 is given by a $(\bar A)^{\perp}$-valued logarithmic 1-form supported on $C_{0, L}$, 
 the part of $C_0$ corresponding to the loop $L$ of $\Gamma$.
Here $\bar A$ is the subspace of $N_{\Bbb R}$ spanned by the direction
 vectors of the edges in the loop.\\

With this observation in mind, we can simplify the calculation of the 
 obstruction by decomposing it into the part annihilated
 by $\bar A^{\perp}$ (which does not actually contribute to the obstruction)
 and the other part.
Explicitly, we can do this as follows.

Recall that we defined 
 the connected subgraph $\overline\Gamma'$ of $\Gamma$ which
 contains the loop, and it has 
 1-valent vertices $\{\alpha_i\}$, see Subsection \ref{subsec:convention}.
Let $E_i$ be the unique 
 edge of $\overline\Gamma'$ adjacent to $\alpha_i$.
Applying a base change and adding a vertex if necessary, 
 we can assume that the end of $E_i$ other than $\alpha_i$
 is a 2-valent vertex which we write by $\alpha_i'$.
Let $E_i'$ be the other edge adjacent to $\alpha_i'$
 and let $\stackrel{\circ}{E'}_i$ be its open part
 (that is, the complement of the ends).
 
Then the subgraph of $\Gamma$ consisting of $\stackrel{\circ}{E'}_i$, $E_i$, 
 $\alpha'_i$, $\alpha_i$, 
 and the components of $\Gamma\setminus
  \overline{\Gamma}'$ which
 contain $\alpha_i$ in their closures (there are two such components)
 is a tree, and we write it by $T_i$ (see Figure 3). 

\begin{figure}[h]\label{fig:T}
\includegraphics{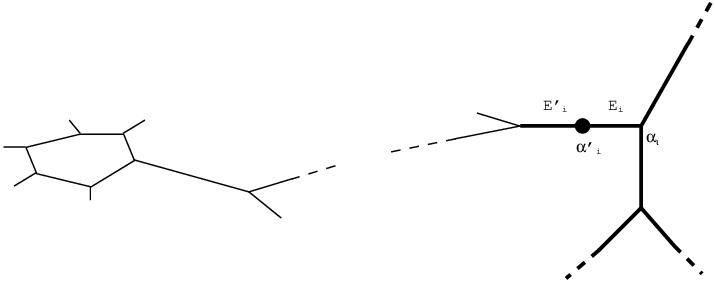}
\caption{The part drawn by bold lines is $T_i$.}
\end{figure}

Let $C_{0, T_i}$ be the chain of rational curves in $C_0$ 
 corresponding to $T_i$.
Let $p_i$ and $p_i'$
 be the nodes of $C_{0}$ corresponding to the edges $E_i$ and $E_i'$,
 respectively. 
Let $C_{0, T_i}^*$ be the open subset 
\[
C_{0, T_i}\setminus \{p_i'\}
\]
 of $C_0$.
For a $k$-th order lift $C_k$ of $C_0$, let $C_{k, T_i}^*$ be the log scheme structure
 on the topological space underlying $C_{0, T_i}^*$ induced by the restriction.
Then we have the following.
\begin{lem}\label{lem:liftT}
Let $\varphi_k\colon C_k\to \mathfrak X$ be a lift of $\varphi_0$.
Then the restriction of $\varphi_k$ to $C_{k, T_i}^*$ has a $(k+1)$-th order lift.
\end{lem}
\proof
For any adjacent vertices $\alpha$ and $\beta$ of $T_i$, 
 there is a trivial lift of $\varphi_k|_{C_{k, T_i}^*}$
 restricted to $(C_{k, \alpha}\cup C_{k, \beta})^*$
 to a map of order $k+1$ by Lemma \ref{lem:canlift}.
If $\beta, \gamma$ is another pair of adjacent vertices of $T_i$, 
 then the difference 
 of these lifts on $C_{k, \beta}^*$
 determines a section of 
 $\mathcal N_{C_0/X_0}|_{C_{0, \beta}^*}
 \otimes_{\Bbb C} t^{k+1}\Bbb C[t]/t^{k+2}\cong \mathcal N_{C_0/X_0}|_{C_{0, \beta}^*}$
 by Lemma \ref{lem:lift3h}.

On the other hand, let $\mathcal N_{C_0/X_0}|_{C_{0, T_i}^*}$
 be the sheaf of sections of $\mathcal N_{C_0/X_0}$ restricted to 
 $C_{0, T_i}^*$.
This can be thought of as the sheaf $\mathcal N_{C_0/X_0}$
 restricted to $C_{0, T_i}$ whose element allows a pole of any order
 at the point $p_i'$.
The restriction of $\mathcal N_{C_0/X_0}$
 to a component of $C_{0, T_i}$ is isomorphic to 
 $\mathcal O_{\Bbb P^1}(1)\oplus\mathcal O_{\Bbb P^1}^{\oplus n-2}$ when the component corresponds
 to a 3-valent vertex, and isomorphic to $\mathcal O_{\Bbb P^1}^{\oplus n-1}$
 when it corresponds to a 2-valent vertex. 
Then using Serre duality for nodal curves,
 it is easy to see that the cohomology group
 $H^1(C_{0, T_i}^*, \mathcal N_{C_0/X_0}|_{C_{0, T_i}^*})$ vanishes.

Now by a general principle of deformation theory, there is a lift of 
 $\varphi_k|_{C_{k, T_i}^*}$ to the next order.
This proves the lemma.\qed\\

Thus, given a $k$-th order lift 
 $\varphi_k\colon C_k\to \mathfrak X$ of $\varphi_0$, 
 we can construct a \v{C}ech cocycle using the following local lifts:
\begin{itemize}
\item On $C_{k, T_i}^*$, take some lift as in Lemma \ref{lem:liftT}.
\item For any two adjacent vertices $\sigma$ and $\tau$ in
 $\overline\Gamma'\setminus \{\alpha_i\}$, take the trivial lift of
 $\varphi_k|_{(C_{k, \sigma}\cup C_{k, \tau})^*}$ 
 as in Lemma \ref{lem:canlift}.  
\end{itemize}
It is clear that the cohomology class constructed in this way is the same as 
 $ob(\varphi_k)$ defined in Definition \ref{def:cechobst}.

As in the proof of Lemma \ref{lem:first_order_deform}, 
 on the closed subset
 $C_{0, {\Gamma'}}$
 of $C_0$ corresponding to the subgraph 
 $\Gamma'$ of $\Gamma$
 (recall that $\Gamma'$ is the graph $\overline\Gamma'$ with the
 1-valent vertices removed), 
 the log normal sheaf $\mathcal N_{C_0/X_0}$
 has a direct sum decomposition
\[
\mathcal N_{C_0/X_0}|_{C_{0, {\Gamma'}}}
 \cong \mathcal N'\oplus
  \mathcal O_{C_{0, {\Gamma'}}}^{\oplus r},
\]
 where $\mathcal N'$ is the subsheaf of 
 $\mathcal N_{C_0/X_0}|_{C_{0, {\Gamma'}}}$
 consisting of the sections annihilated by vectors in 
 $\bar A^{\perp}$,
 and $r= n-\dim \bar A$, see also Lemma \ref{lem:trivsub} below.  
\begin{rem}\label{rem:decomp}
This direct sum decomposition is not canonical.
The part $\mathcal N'$ is uniquely determined, but 
 there are choices for the other summand. 
We fix one such a choice.
This choice does not affect the following argument.
\end{rem}
\begin{prop}\label{prop:decomposition}
A representative of the obstruction class can be taken so that its support
 is contained in $C_{0, \Gamma'}$.
In particular, the obstruction class lies in the group
 $H^1(C_{0, \Gamma'}, \mathcal N_{C_0/X_0}|_{C_{0, \Gamma'}})$.
Moreover, in the decomposition
\[
H^1(C_{0, \Gamma'}, \mathcal N_{C_0/X_0}|_{C_{0, \Gamma'}})
 \cong 
 H^1(C_{0, \Gamma'}, \mathcal N')
  \oplus H^1(C_{0, \Gamma'}, \mathcal O_{C_{0, \Gamma'}}^{\oplus r}),
\]
 the component of the obstruction class in 
 $H^1(C_{0, \Gamma'}, \mathcal N')$ is zero.
\end{prop}
\begin{rem}
The inclusion $C_{0, \Gamma'}\to C_{0, \Gamma}$ induces a map 
 between the cohomology groups
 $H^1(C_{0, \Gamma}, \mathcal N_{C_0/X_0})
  \to H^1(C_{0, \Gamma'}, \mathcal N_{C_0/X_0}|_{C_{0, \Gamma'}})
  $.
By the fact that the complement $\Gamma\setminus\Gamma'$ is a union of
 trees, it is easy to see that this map is an injection
 (in fact, an isomorphism).
\end{rem}
\proof
By the construction of the above local lifts,
 the obstruction cocycle is represented by 
 a \v{C}ech 1-cocyle  
 with values in 
 $\mathcal N_{C_0/X_0}$
 whose support is contained in $C_{0, {\Gamma'}}$.
Therefore, the first statement follows.

By the splitting 
 $\mathcal N_{C_0/X_0}|_{C_{0, {\Gamma'}}}
 \cong \mathcal N'\oplus
  \mathcal O_{C_{0, {\Gamma'}}}^{\oplus r}$ over $C_{0, \Gamma'}$,
  the cocycle can be decomposed into
 the $\mathcal N'$ part and
 the $\mathcal O_{C_{0, {\Gamma'}}}^{\oplus r}$ part.
As we noted in Subsection \ref{subsec:abundance}, any element of 
 the dual space of 
 $H^1(C_0, \mathcal N_{C_0/X_0})$ is represented by 
 a logarithmic 1-form with values in $(\bar A)^{\perp}$
 supported on $C_{0, L}$.
Thus, the $\mathcal N'$ part of the obstruction cocycle 
 is annihilated by any element of 
 $H^1(C_0, \mathcal N_{C_0/X_0})^{\vee}$
 (for the explicit definition of the pairing between 
 \v{C}ech representatives of $H^1(C_0, \mathcal N_{C_0/X_0})$ and 
 1-forms of $H^0(C_0, \mathcal N_{C_0/X_0}\otimes\omega_{C_0})$, 
 see Subsection \ref{subsec:resobs}).
Consequently, this part must be zero. \qed\\

Thus, we need to consider only the 
 $\mathcal O_{C_{0, {\Gamma'}}}^{\oplus r}$
 valued part of the obstruction.

\subsection{Calculation of obstructions in embedded cases. Step 4: Analytic continuation of local lifts.}\label{subsec:step4}
Now we study the 
 $\mathcal O_{C_{0, {\Gamma'}}}^{\oplus r}$ part
 of the obstruction.
Although in the previous subsection we represented the obstruction class 
 as a \v{C}ech cocycle, now we prefer to study the
 deformation through analytic continuation.

The basic idea is the following.
Namely, suppose we are given a 
 $(k-1)$-th order lift $\varphi_{k-1}\colon C_{k-1}\to \mathfrak X$ of $\varphi_0$
 as before.
Lemmas \ref{lem:canlift} and \ref{lem:liftT}
 give local $k$-th order lifts on suitable
 open subsets of $C_{k-1}$.
If one of such local lifts on some open subset 
 can be analytically extended over $C_{k-1}$
 (this means we also choose a suitable lift $C_{k}$
 of $C_{k-1}$), then it gives a 
 $k$-th order lift of $\varphi_{k-1}$, as we argued in Subsection 
 \ref{subsec:Kuranishi}.

Though the calculation of such 
 an extension of a local lift of a map $\varphi_k$
 is difficult in general, by the argument in Subsection \ref{subsec:step3}
 (see in particular Proposition \ref{prop:decomposition}),
 we only need to 
 calculate the part which is not annihilated by $\bar A^{\perp}$. 
This reduces the calculation of the extension of a map to the 
 calculation of analytic continuation of functions.

Here an important difference from the discussion in Subsection 
 \ref{subsec:Kuranishi} is that we need to allow the deformation
 of the domain of the map.
This point makes it difficult to fully calculate the analytic
 continuation of a function.
However, since each connected component of the complement
 $C_0\setminus C_{0, L}$ (here $C_{0, L}$ is the chain of rational
 curves corresponding to the loop $L$ of $\Gamma$)
 is a tree of rational curves, the existence of lifts of $\varphi_0$ to any
 order is guaranteed on these parts as in Lemma \ref{lem:liftT}.
Using this fact, we can calculate the
 terms with the lowest exponent 
 of
 $t$ (except the constant terms)
 of 
 the analytic continuation.
Since $t$ is supposed to be small (algebraically, it is nilpotent),
 such a term dominates the others,
 and it suffices for our purpose by
 a suitable transversality result (see Lemma \ref{lem:transversality}).
In fact, the lowest order term can be calculated just 
 from the information of $\varphi_0$.
Thus we do not need to 
 worry about how the domain curve should be deformed.

Such an analytic continuation results in a pole on the loop part 
 $C_{0, L}$ (or on its suitable $k$-th order local lift).
Thus, through the calculation of analytic continuations on each of the 
 components of $C_0\setminus C_{0, L}$,
 we have poles distributed on 
 $C_{0, L}$ .
These poles define a cohomology class on $C_{0, L}$, 
 which is the obstruction to extend these analytic continuations
 over a lift of $C_{0, L}$.
Then the problem of lifting 
 the map $\varphi_k$ is reduced to the
 existence of functions on $C_{0, L}$ 
 with the prescribed pole.
See Subsections \ref{subsec:resobs}
 and \ref{subsec:pairing} for the precise formulation.

First, we describe the decomposition
\[
\mathcal N_{C_0/X_0}|_{C_{0, {\Gamma'}}}
 \cong \mathcal N'\oplus
  \mathcal O_{C_{0, {\Gamma'}}}^{\oplus r},
\]
 of the log normal sheaf 
 of $\varphi_0|_{C_{0, {\Gamma'}}}$
 a little more explicitly.
Recall that we write by $\bar A$ the subspace of $N_{\Bbb R}$
 spanned by the direction vectors of the edges in the loop $h(L)$.
Take a basis $\mathfrak v_1, \dots, \mathfrak v_{n-r}$ of $N\cap \bar A$
 and take vectors $\mathfrak w_1, \dots, 
 \mathfrak w_r$ of $N$ so that the set of vectors 
 $\{\mathfrak v_1, \dots, \mathfrak v_{n-r}, 
 \mathfrak w_1, \dots, \mathfrak w_r\}$ is a basis of $N$.
These vectors correspond to log tangent vectors on $\mathfrak X$, and 
 we have an isomorphism
\[
\Theta_{X_0}\cong \left(\mathcal O_{X_0}\cdot \mathfrak v_1\oplus\cdots
 \mathcal O_{X_0}\cdot \mathfrak v_{n-r}\right)
 \oplus\left(\mathcal O_{X_0}\cdot \mathfrak w_1\oplus\cdots
 \mathcal O_{X_0}\cdot \mathfrak w_{r}\right).
\]
This splitting induces a splitting of $\mathcal N_{C_0/X_0}$
restricted to $C_{0, \Gamma'}$ as follows.
\begin{lem}\label{lem:trivsub}
The restriction of the log normal sheaf
 $\mathcal N_{C_0/X_0}$ to $C_{0, {\Gamma'}}$
 has a splitting
\[
\mathcal N_{C_0/X_0}|_{C_{0, {\Gamma'}}}
 \cong \varphi_0|_{C_{0, \Gamma'}}^*
 \left(\mathcal O_{X_0}\cdot \mathfrak v_1\oplus\cdots
 \mathcal O_{X_0}\cdot \mathfrak v_{n-r}\right)/\Theta_{C_{0, \Gamma'}}
 \oplus\varphi_0|_{C_{0, \Gamma'}}^*
 \left(\mathcal O_{X_0}\cdot \mathfrak w_1\oplus\cdots
 \mathcal O_{X_0}\cdot \mathfrak w_{r}\right).
\]\qed
\end{lem}
\begin{rem}\label{rem:direction}
The part 
 $\varphi_0|_{C_{0, {\Gamma'}}}^*
   (\mathcal O_{X_0}\cdot \mathfrak v_1\oplus\cdots\oplus
   \mathcal O_{X_0}\cdot \mathfrak v_{n-r})/\Theta_{C_0}$
  is the summand
  $\mathcal N'$
  introduced in Subsection \ref{subsec:step3}, and
  the part 
  $\varphi_0|_{C_{0, {\Gamma'}}}^*(\mathcal O_{X_0}
    \cdot  \mathfrak w_1\oplus\cdots
   \mathcal O_{X_0}\cdot \mathfrak w_r)$
   is the summand 
   $ \mathcal O_{C_{0, {\Gamma'}}}^{\oplus r}$
   there.
As we mentioned in Remark \ref{rem:decomp}, 
 this splitting
  is not canonical.
Later it turns out that the choice of a splitting does
 not matter in our argument 
 $($see Remark \ref{rem:basis-ob}$)$.
\end{rem}

\subsubsection{Calculation on standard  curves with two components} \label{subsubsec:std}
Now we begin the calculation of the analytic continuation.
We consider the tropical curve $(\Gamma, h)$ given at the beginning 
 of Subsection \ref{subsec:weight one}, and use the same notations 
 $\Gamma', \alpha, \beta, E$ etc..
First we assume that the restriction of $h$ to a neighborhood of 
 any bounded edge is standard (see Definition \ref{def:twovertex}) for simplicity.
More general cases are treated later.

Also, recall that we chose the vertex $\alpha$ so that the integral length from it to 
 the loop $h(L)$ is the shortest among those from any 1-valent vertex
 of $\mathcal G = h(\overline{\Gamma}')$.

As we saw in Lemma \ref{lem:lift3},
 there is a direct correspondence between
 the set of local sections of the log normal sheaf and the set of local lifts.
Explicitly, a local section of the log normal sheaf
 $\mathcal N_{C_{0}/\mathfrak X}$ restricted to the
 component $C_{0, \alpha}^*$ given by the expression
\[
\mathfrak n_{0, \alpha} = 
 y\partial_y + c(s)(x\partial_x-y\partial_y) + c_0(s)z\partial_z + 
 c_1(s)w_1\partial_{w_1}+\cdots +c_{n-2}(s)w_{n-2}\partial_{w_{n-2}}
\]
 corresponds to a lift of the map $\varphi_0$ restricted 
 to $C_{0, \alpha}^*$
 whose image 
 is parametrized as
\[
(\dagger)\;\; \begin{array}{l}
(x, y, z, w_1, w_2, \dots, w_{n-2})  \\
\hspace{.1in}= ((1+c(s))s, \frac{t}{s}, (1+c_0(s))(-\frac{k}{l}s-\frac{m}{l}),
(1+c_1(s))(-\frac{\mu}{\lambda}), (1+c_2(s))a_2, \dots, (1+c_{n-2}(s))a_{n-2}).
\end{array}
\]
 using the local coordinate $\varphi_0^*(x) = s$ on $C_{0, \alpha}^*$
 as in Lemma \ref{lem:lift3} and Remark \ref{rem:secorder}
 (see also the remark before Lemma \ref{lem:lift3} 
 noting that the coordinate $s$ can be
 regarded as a coordinate on $C_{1, \alpha}^*$).
Here $c(s), c_0(s), \dots, c_{n-2}(s)$ belong to 
 $t\Bbb C[t]/t^2\otimes_{\Bbb C}\mathcal O_{C_{0, \alpha}^*}$ as in 
 Remark \ref{rem:secorder}.

By Proposition \ref{prop:decomposition}, we only need to consider
 analytic continuations of the pull back of the coordinates 
 on $\mathfrak X$ corresponding to vectors which are not annihilated by
 $\mathfrak w_1, \dots, \mathfrak w_r$.
In the parameterization $(\dagger)$ above, we can take the functions
 $(x, y, z, w_1, \dots, w_{n-2})$ so that
 the only non-constant
 coordinate in these directions is $z$ 
 (see Remark \ref{rem:functions} for more general cases).
Therefore, in this subsection
 we consider the analytic continuation of the locally defined function
\[
z = (1+c_0(s))(-\frac{k}{l}s-\frac{m}{l})
\]
 to a suitable lift of $C_0$.
The calculation containing other directions 
 is done in Subsection \ref{subsec:step6}.
For notational simplicity, we set the constant term
 to one by multiplying $-\frac{l}{m}(1-c_{0, 0}t)$ to $z$.
Here $c_{0, 0}t$ is the constant term of $c_{0}(s)$.
Thus, we consider the analytic continuation of a function of the form
\[
\frac{k}{m}s + O(t)\bar c_0(s)+1,
\]
 where $\bar c_0(s)$ is a rational function of $s$ without a constant term.

Recall we are assuming that there is a $(k-1)$-th order lift $\varphi_{k-1}$
 of $\varphi_0$ on the whole $C_{k-1}$, and 
 also a $k$-th order lift $\varphi_k$ on some subset of $C_{k-1}$
 (explicitly, on the lift of $C_0\setminus C_{0, L}$).
By Lemma \ref{lem:first_order_deform},
 we can assume $k\geq 2$.
We assume the lift $(\dag)$ is the image of the restriction of $\varphi_k$
 to $C_{k, \alpha}^*$
 (here we use Lemma \ref{lem:lift3h} and now 
 $c(s), c_0(s), \dots, c_{n-2}(s)$ are in 
 $t\Bbb C[t]/t^{k+1}\otimes_{\Bbb C}
 \Gamma(C_{0, \alpha}^*, \mathcal O_{C_{0, \alpha}^*}$)).
In particular, the function $c_0(s)$ does not
 have a pole at $s = 0$, since otherwise it violates Lemma \ref{lem:lifthigh}
 (see the paragraph before Lemma \ref{lem:lifthigh}).
There is the unique shortest path in $\Gamma$ connecting the vertex $\alpha$
 and the loop $L$.
Let 
\[
\alpha = \alpha_0\to \alpha_1\to \dots\to \alpha_N\in L
\]
 be the sequence of the vertices along the path.
Take a vertex $\alpha_i$ ($0<i\leq N$) and consider 
 the pull back of the function $z$ on $\mathfrak X$
 by the map $\varphi_k$ restricted to
  $C_{k, \alpha_i}^*$.

The function $(\varphi_k|_{C_{k, \alpha_i}^*})^*z$
 is a rational function on 
 $C_{k, \alpha_i}^*\cong C_{0, \alpha_i}^*\times\Bbb C[t]/t^{k+1}$
 allowing poles only at the punctures.
The curve $C_{k, \alpha_i}^*$ has two or three punctures, 
 and here we assume that it has three punctures
 (when it has two punctures, ignore the puncture with the coordinate $a$
 in the argument below).
Let us take a parameter $S_i$ on $C_{k, \alpha_i}^*$ such that the 
 coordinate of the punctures are $0, a$ and $\infty$,
 where $a$ is a nonzero complex number.
Then the function $(\varphi_k|_{C_{k, \alpha_i}^*})^*z$
 is of the form
\[
f_1\left(\frac{1}{S_i}\right) + f_2\left(\frac{1}{S_i-a}\right)+ f_3\left(S_i\right) + c, 
\]
 where $f_i$ are polynomials with coefficients in $\Bbb C[t]/t^{k+1}$ 
 without constant terms and $c$ is an element
 in $\Bbb C[t]/t^{k+1}$.
In particular, we note the following elementary observation.
\begin{lem}\label{lem:expand}
When we expand the function $(\varphi_k|_{C_{k, \alpha_i}^*})^*z$
 into Laurent series at the punctures, functions $f_i$ give the
 principal parts, and $(\varphi_k|_{C_{k, \alpha_i}^*})^*z$ is determined
 by these principal parts up to a constant.
On the other hand, 
 the non-principal part of the Laurent expansion is determined by the 
 principal parts at the other punctures up to the constant term. \qed
\end{lem}

Now consider the puncture corresponding to $S_i = 0$.
Let $C_{0, \beta_i}$ be the component of $C_0$ neighboring $C_{0, \alpha_i}$
 sharing the node corresponding to $S_i = 0$.
There is a parameter $T_i$ on $C_{k, \beta_i}$ satisfying the relation 
\[
S_iT_i = t,
\]
 (again we assume the length of the bounded edge of $h(\Gamma)$ 
 which connects the vertices $\alpha_i$ and $\beta_i$ to be one.
 See Remark \ref{rem:coord} (1) and Remark \ref{rem:secorder2} (3)).

Since $\varphi_k$ is globally defined on a suitable lift (which we write by $C_{k, L^c}$)
 of components of 
 $C_0\setminus C_{0, L}$,
 the function 
 $(\varphi_k|_{C_{k, \alpha_i}^*})^*z$ must be compatible with the pull back 
 $(\varphi_k|_{C_{k, \beta_i}^*})^*z$.
In particular, the principal part $f_1\left(\frac{1}{S_i}\right)$ of
 $(\varphi_k|_{C_{k, \alpha_i}^*})^*z$ at the puncture $S_i = 0$
 must be of the form that it is obtained by substituting 
 $T_i = \frac{t}{S_i}$ to a polynomial $g_1(T_i)$ of $T_i$ with coefficients 
 in $\Bbb C[t]/t^{k+1}$.

Therefore, the principal part $f_1\left(\frac{1}{S_i}\right)$
 is determined by the non-principal part of $(\varphi_k|_{C_{k, \beta_i}^*})^*z$
 at the puncture $T_i = 0$ on $C_{k, \beta_i}^*$.
By Lemma \ref{lem:expand}, 
 the non-principal part of $(\varphi_k|_{C_{k, \beta_i}^*})^*z$ at $T_i = 0$
 is determined by the principal parts of $(\varphi_k|_{C_{k, \beta_i}^*})^*z$
 at the other punctures of $C_{k, \beta_i}^*$.

\begin{figure}[h]
\includegraphics{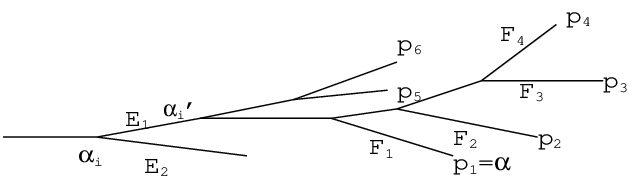}
\caption{}\label{fig:0z}
\end{figure}

Repeating this, we have the following.
Let $\overline\Gamma'_{\alpha}$ be the connected component of 
 the complement $\overline\Gamma'\setminus L$ containing the vertex $\alpha$.
For two vertices $p, q$ of $\overline\Gamma'_{\alpha}$, we write
 $p>q$ if $p$ is on the unique path from $q$ to the loop $L$. 
Take a vertex $\alpha_i$ as above.
It has two edges $E_1$ and $E_2$ in the directions opposite to the loop.
Assume $E_1$ is a bounded edge and let $\alpha_i'$ be the other end of the edge $E_1$.
Let $p_1, \dots, p_a$ be the 1-valent vertices of $\overline\Gamma'$
 satisfying $\alpha'_i>p_j$.
Let $F_j$ be the unique edge towards the loop emanating from $p_j$.
See Figure \ref{fig:0z}.
\begin{lem}\label{lem:prinpart}
The principal part of $(\varphi_k|_{C_{k, \alpha_i}^*})^*z$
 at the puncture corresponding to the edge $E_1$ can be calculated 
 from the non-principal parts of the functions 
 $(\varphi_k|_{C_{k, p_j}^*})^*z$
 at the punctures corresponding to the edges $F_j$.  \qed
\end{lem}

Now we explicitly calculate the analytic continuation which 
 contributes to the principal part of $(\varphi_k|_{C_{k, \alpha_i}^*})^*z$
 at the puncture corresponding to $E_i$.
Since the calculation is the same at any of $p_1, \dots, p_a$, 
 let us take $p_1$ and assume it coincides with the vertex 
 $\alpha$ taken in the above argument.
More precisely, the pull back $(\varphi_k|_{C_{k, \alpha_i}^*})^*z$
 is not easy to express directly.
Instead, we start from the function $(\varphi_k|_{C_{k, \alpha}^*})^*z$
 which is easy to describe, and calculate a part of $(\varphi_k|_{C_{k, \alpha_i}^*})^*z$
 relevant to us by the analytic continuation of $(\varphi_k|_{C_{k, \alpha}^*})^*z$.

The image of $\varphi_{k}$ restricted to 
 $C_{k, \alpha}^*\cong C_{0, \alpha}^*\times Spec\Bbb C[t]/t^{k+1}$
 has a parameterization similar to $(\dagger)$, where $c(s), c_0(s), \dots, 
 c_{n-2}(s)\in t\Bbb C[t]/t^{k+1}\otimes_{\Bbb C}
 \Gamma(C_{0, \alpha}^*, \mathcal O_{C_{0, \alpha}^*})$ .
Also, we use the coordinate on $C_{k, \alpha}^*$
 given by $S = \varphi_{k}^*(x)$.
This coordinate $S$ equals $s$ (considered as a function on 
 $C_{k, \alpha}^*$ by the isomorphism 
 $C_{k, \alpha}^*\cong C_{0, \alpha}^*\times Spec\Bbb C[t]/t^{k+1}$,
 see the paragraph before Lemma \ref{lem:lift3})
 modulo $t$.
Thus, the pull back of the $z$-coordinate to $C_{k, \alpha}^*$
 by the map $\varphi_k$ 
 has the form 
\[
\left(\frac{k}{m}+O(t)\right)S+ O(t)\bar\xi(S)+1+O(t),
\]
 here $\bar\xi(S)$ is a rational function of $S$ with coefficients in $\Bbb C[t]/t^{k+1}$
 (recall that we normalized the function $z$ so that its constant term
 is one modulo $t$). 
 
 
In view of Lemma \ref{lem:prinpart}, we discard the principal part of 
 the pull back of the $z$-coordinate to $C_{k, \alpha}^*$
 at the puncture corresponding to the edge $F_1$,
 and write the non-principal part by
\[
\left(\frac{k}{m}+O(t)\right)S+ O(t)S^2\xi(S)+1+O(t),
\]
 here $\xi(S)$ is a series of $S$ with coefficient in $\Bbb C[t]/t^{k+1}$
 which converges on a punctured disk around $S = 0$.

Let $C_{k, \alpha_1}^*$ be the open subset of the component of $C_k$
 corresponding to the vertex $\alpha_1$ of $\Gamma$ neighboring $\alpha$.
The image of 
$C_{0, \alpha_1}$ by the map
$\varphi_0$ 
is given by the defining equations
\[
\kappa y+\lambda w_1+\mu
= 0, \;\;x = 0, \;\; z = -\frac{m}{l}, \;\; w_2 = a_2, \dots, w_{n-2} = a_{n-2}
\]
(see Lemma \ref{lem:coordinate}).
Let us write $S_1 = \varphi_k^*(y)$.
It gives a parameter on $C_{k, \alpha_1}^*$
 and satisfies 
 \[
 SS_1 = t,
 \]
 when $S$ and $S_1$ are extended to $(C_{k, \alpha}\cup C_{k, \alpha_1})^*$.

Then the principal part of the function $(\varphi_k|_{C_{k, \alpha_1}^*})^*z$
 on the component $C_{k, \alpha_1}^*$
 at the puncture $S_1 = 0$ has the form 
\[
\left(\frac{k}{m}+O(t)\right)S+O(t)S^2\xi(S) =
 \left(\frac{k}{m}+O(t)\right)\frac{t}{S_1}+O(t)\frac{t^2}{S_1^2}
  \xi\left(\frac{t}{S_1}\right).
 \]

\begin{figure}[h]
	\includegraphics{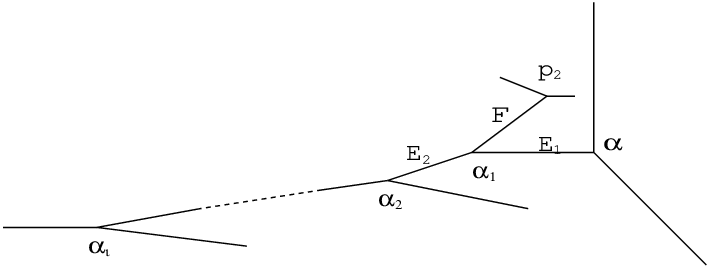}
	\caption{The principal part at the puncture corresponding to the 
		edge $E_2$ on the component $C_{k, \alpha_2}^*$ has contributions
		from the principal parts at the punctures corresponding to the edges
		$F$ and $E_1$ on the component $C_{k, \alpha_1}^*$.
		Since the vertex $\alpha$ is (one of) the closest to the loop among 1-valent vertices
		of $\overline\Gamma'$, the order of $t$
		of the contribution from the puncture corresponding to 
		$F$ is at least two (which is the 
		order of the contribution from the puncture corresponding to $E_1$).
		}\label{fig:0b}
\end{figure}

The image of $C_{k, \alpha_1}^*$ by the map $\varphi_{k}$ 
 has defining equations which reduce to the above expression 
 for the image of $C_{0, \alpha_1}^*$ over $\Bbb C[t]/t$.
Assume that the point given by $w_1 = 0$ on $C_{0, \alpha_1}$
 is the node corresponding to the edge of $\Gamma$ connecting
 the vertices $\alpha_1$ and $\alpha_2$
 (the other possibility is that the point given by $y, w_1  \to \infty$ on
 $C_{0, \alpha_1}^*$
 is this node, and the calculation in this case is similar).
If $T_1 = \varphi_{k}^*w_1$ is another coordinate on
 $C_{k, \alpha_1}^*$, then it satisfies
\[
S_1 = -\frac{\lambda}{\kappa}T_1-\frac{\mu}{\kappa} \;\;\;\; \text{mod $t$}
\]
 by the equation $\kappa y+\lambda w_1+\mu = 0$.
Using this parameter,
 the part of the $z$-coordinate
 $\left(\frac{k}{m}+O(t)\right)S+ O(t)S^2\xi(S)+1+O(t)$ above
 becomes
\[\begin{array}{ll}
\hspace{-.1in}\left(\frac{k}{m}+O(t)\right)S+O(t)S^2\xi(S)+1+O(t) & =
 \left(\frac{k}{m}+O(t)\right)\frac{t}{S_1}+O(t)\frac{t^{2}}{S_1^2}
 \xi(\frac{t}{S_1})+1+O(t)\\
& = \left(\frac{k}{m}+
   O(t)\right)\frac{t}{-\frac{\lambda}{\kappa}T_1-\frac{\mu}{\kappa}}+
   \left(\frac{t}{-\frac{\lambda}{\kappa}T_1-\frac{\mu}{\kappa}}\right)^2
   O(t)\xi_1(\frac{t}{-\frac{\lambda}{\kappa}T_1-\frac{\mu}{\kappa}})
   +1+O(t).
\end{array}
\]
Here $\xi_1$ is a polynomial with coefficients in $\Bbb C[t]/t^{k+1}$.
Let $S_2$ be the coordinate on $C_{k, \alpha_2}^*$ which satisfies
\[
S_2T_1 = t,
\]
 when $S_2$ and $T_1$ are extended to $(C_{k, \alpha_1}\cup C_{k, \alpha_2})^*$.
Then, expanding the terms including
 $\frac{t}{-\frac{\lambda}{\kappa}T_1-\frac{\mu}{\kappa}}$
 by using
\[
\frac{t}{-\frac{\lambda}{\kappa}T_1-\frac{\mu}{\kappa}}
 = \frac{t}{-\frac{\lambda}{\kappa}\frac{t}{S_2}-\frac{\mu}{\kappa}}
 = -\frac{\kappa}{\mu}t\left(1- \frac{\lambda}{\mu}\cdot\frac{t}{S_2}
 +\left(\frac{\lambda}{\mu}\right)^2\cdot\frac{t^{2}}{S_2^2}+\cdots\right),
\]
 the expression for the part of the pull back of the $z$-coordinate above
 is written as 
\[
(1+O(t))t^{2}\frac{k}{m}\cdot\frac{\kappa}{\mu}\cdot\frac{\lambda}{\mu}
 \cdot\frac{1}{S_2} + t
  O(t)\frac{t^{2}}{S_2^2}\xi_2\left(\frac{t}{S_2}\right)
 +1+O(t)
\]
  on $C_{k, \alpha_2}^*$,
 here $\xi_2$ is a polynomial.
This gives the part of the analytic continuation of the $z$-coordinate 
 on the image of
 $\varphi_{k}|_{C_{k, \alpha}^*}$ to the image of 
 $\varphi_{k}|_{C_{k, \alpha_2}^*}$.
The principal part of $(\varphi_k|_{C_{k, \alpha_2}^*})^*z$
 at the puncture $S_2 = 0$ is given by the sum of (the non-constant part of)
 this and 
 another contribution from the principal part of 
 $(\varphi_k|_{C_{k, \alpha_1}^*})^*z$
 at the puncture corresponding to $S_1 = \infty$, 
 see Figure \ref{fig:0b}. 

\begin{rem}\label{rem:leadterm}
As we noted at the beginning of Subsection \ref{subsec:step4}, 
 we have not calculated the precise form of the analytic continuation.
In fact, we have not specified what curve the lift $C_{k, L^c}$ of (a subset of) $C_0$ is.
Thus, it is impossible to calculate the precise form.
However, whatever the curve $C_k$ is, the leading term with respect
 to the order of $t$ has the same form, and
 it suffices for our purpose. 
Precisely, in the above expression, the non-constant term with the lowest
 power of $t$ is 
 $\frac{k}{m}\cdot\frac{\kappa}{\mu}\cdot\frac{\lambda}{\mu}
 \cdot\frac{1}{S_2}t^{2}$, 
 and it depends only on the data of $\varphi_0(C_0)$.
We will see that such a term is identified with a contribution to the obstruction
 to extend the partially defined $k$-th order lift of $\varphi_0$ to the whole curve. 
It follows that, given $\varphi_{k-1}$, the dominant part of its obstruction
 to lift it to the $k$-th order depends only on the data of $\varphi_0$.
\end{rem}

Repeating the above
 calculation, we can analytically continue the part of the $z$-coordinate
 at the puncture $S = 0$ on $C_{k, \alpha}^*$
 to each component $C_{k, \alpha_j}$ corresponding to the vertices on the path
 from $\alpha$ to the loop.
However, the above calculation assumes that the tropical curve
 $(\Gamma, h)$ is standard in a neighborhood 
 of the edge connecting $h(\alpha)$ and $h(\alpha_1)$ 
 as well as that of the edge connecting $h(\alpha_1)$ and $h(\alpha_2)$.
In particular, in the case when 
 the direction vectors of the edges emanating from the vertices 
 $h(\alpha_1)$ and $h(\alpha_2)$ 
 span a one or two dimensional subspace of $N_{\Bbb R}$,
 the above calculation cannot be applied
 (by definition of $\alpha$ which is a 1-valent vertex of 
 $\overline\Gamma'$, the direction vectors 
 of the edges emanating from 
 $h(\alpha)$ and $h(\alpha_1)$ span a three dimensional subspace).
In the subsection below, we study how to remove this assumption.

\subsubsection{Calculation on non-standard  curves with two components} \label{subsubsec:nonstd}
We now study how to proceed the calculation of the analytic continuation
 of the function $z$ in the cases where the local structure of 
 $(\Gamma, h)$ is not isomorphic to the standard one.
Let $\alpha', \beta'$ be a pair of adjacent vertices of 
 $\overline\Gamma'$ (see Subsection \ref{subsec:convention}). 
There are two cases to study:
\begin{enumerate}
\item[(a)] 
One of the vertices $\alpha', \beta'$ (say $\alpha'$)
 is a 1-valent vertex of $\overline\Gamma'$.
In this case, 
 the direction vectors 
 of the edges emanating from two vertices 
 $h(\alpha')$ and $h(\beta')$ span a three dimensional subspace.
\item[(b)] 
Either of the vertices $\alpha', \beta'$ is not 
 a 1-valent vertex of 
 $\overline\Gamma'$.
In this case, the direction vectors 
 of the edges
 emanating from  
 $h(\alpha')$ and $h(\beta')$ span a one, two or three dimensional subspace.
\end{enumerate}
\noindent
{\bf Case (a).}
Let $(\Xi, h_{\Xi})$ be a tropical curve with
 two vertices
 and $\mathcal L$
 be the affine linear transformation which maps a 
 standard tropical curve with two vertices $\Gamma_s$ onto $h_{\Xi}(\Xi)$
 as in Lemma \ref{lem:twovertex}.
Let $\mathfrak X$ and $\mathfrak X_s$ be toric degenerations defined respecting $(\Xi, h_{\Xi})$
 and $\Gamma_s$, respectively.
Let $X_0$ and $X_{s, 0}$ be their central fibers.
Following Remark \ref{rem:coord} (1) (see also Remark \ref{rem:secorder2} (3)),
 we assume the integral length of the unique 
 bounded edge of
 $h_{\Xi}(\Xi)$ is equal to its weight.
In this case, the integral length of the bounded edge of $\Gamma_s$ is one.

A pre-log curve $\varphi_0\colon  C_{0,\alpha'}\cup C_{0,\beta'}\to X_0$
 of type $(\Xi, h_{\Xi})$ 
 is given as
 the image of a pre-log curve $\varphi_{s, 0}\colon C_{0,\alpha'}\cup C_{0,\beta'}
 \to X_{s, 0}$ 
 of type $\Gamma_s$ by
 the branched covering map 
 $\widetilde{\Phi}_{\mathcal L}|_{X_{s,0}}\colon X_{s, 0}\to X_0$ 
 (see Subsection \ref{subsec:general_edge}).
Here $C_{0,\alpha'}\cup C_{0,\beta'}$ is
the nodal union of nonsingular rational curves.
We again assume $\dim X_0 = n = 3$ for notational simplicity, 
 as in Subsection \ref{subsec:general_edge}.
The image of 
 the pre-log curve $\varphi_{s, 0}$ of type $\Gamma_s$ is presented by the equations
\[
kx+lz+m = 0,\;\; y = 0, \;\; w = -\frac{\mu}{\lambda},
\]
\[
\kappa y+\lambda w+\mu = 0,\;\; x = 0,\;\; z = -\frac{m}{l},
\]
 using the same notation as in Lemma \ref{lem:coordinate}.
Let $\varphi_k\colon (C_{k,\alpha'}\cup C_{k,\beta'})^*\to \mathfrak X$
 be a $k$-th order lift of $\varphi_0$ restricted to the open subset
 underlying $(C_{0,\alpha'}\cup C_{0,\beta'})^*$, see Definition \ref{def:twin}.
We think of this as a restriction of a global curve $C_{k, L^c}\to \mathfrak X$ considered  
 in the first half of this Subsection \ref{subsec:step4}.

The functions $x, y$ satisfy the equation
\[
xy = t.
\]
As in Lemma \ref{lem:coordinate} (iii),
 on the pre-log curve $\varphi_{s, 0}$, 
 the pull back of the functions $x, y, z$ satisfy
\[
(\varphi_{s, 0})^*(x) = s,\;\;
 (\varphi_{s, 0})^*(y) = \sigma,\;\; (\varphi_{s, 0})^*(z) = -\frac{1}{l}(ks+m),
\]
 for appropriate coordinates $s, \sigma$ around the node of 
 $C_{0,\alpha'}\cup C_{0,\beta'}$.

Let $X_{s, 0, \alpha'}$ and $X_{s, 0, \beta'}$ be the components of 
 $X_{s, 0}$ to which $C_{0, \alpha'}$ and $C_{0, \beta'}$ are mapped,
 respectively.
Recall that in Subsection \ref{subsec:general_edge},
 we introduced 
 functions
 $X, Y, Z$ and $W$ on $\mathfrak X$ such that
 $X, Z, W$ and $Y, Z, W$ give coordinate systems on
 open subsets of the components
  $X_{0, \alpha'}$ and $X_{0, \beta'}$ of the central fiber $X_0$ of
  $\mathfrak X$, respectively.
As in the case of standard tropical curves with two vertices 
 (that is, the calculation 
 we did in Subsection \ref{subsubsec:std}), 
 it suffices to calculate the analytic continuation of the function 
 $Z$ for the calculation of the obstruction.
Precisely, we note the remark below.
\begin{rem}\label{rem:functions}
In the notation of Subsection \ref{subsec:general_edge}, in general 
 the set of functions $X, Z, W$ or $Y, Z, W$
 does not correspond to the dual basis of 
 the vectors $u, v, w$ generating the edges of the tropical curve.
In particular, the vectors in $N^{\vee}$ 
 corresponding to the functions $X$, $Y$ and $W$ do not necessarily annihilate 
 the vectors $\mathfrak w_i$ in the decomposition of
 Lemma \ref{lem:trivsub}. 

From the precise relation between analytic continuations and the obstruction
 below (see Proposition \ref{prop:anal-pair} and Corollary \ref{cor:anal-ob}), 
 this does not cause a problem because of the linearity of the obstruction with respect to 
 the directions (Lemma \ref{lem:linearity}).

From the point of view of analytic continuations, this can also be understood as follows.
Namely, the above property of the vectors corresponding to 
 the functions $X$, $Y$ and $W$ results in
 the possibility that not only the 
 pull back of the function $Z$ but also the pull back of the
 functions $X$, $Y$ and $W$ may have the obstruction to the analytic continuation.
Namely, Proposition \ref{prop:decomposition}
 assures that functions can be analytically continues globally only for 
 those corresponding to the vectors annihilated by $\mathfrak w_i$, $i = 1, \dots, r$.
However, the pull back of the function $Z$ has the form
\[
c_0 + O(t), 
\]
 on the part $C_{k, L^c}$ of $C_k$,
 as we will see below.
Here $c_0$ is a non-zero constant, and $O(t)$ is a function
 (not a constant) on $C_{k, L^c}$ which is zero modulo $t$.
Therefore, we can unambiguously fix 
 roots of $Z$.
Then, for some rational numbers $q_1$, $q_2$ and $q_3$, the functions
\[
XZ^{q_1},\; YZ^{q_2},\; WZ^{q_3}
\]
 correspond to vectors in $N^{\vee}\otimes\Bbb Q$ which annihilate the vectors
 $\mathfrak w_i$ above.
It follows that
 the analytic continuation of these functions 
 are not obstructed.
Thus, if $Z$ can be analytically continued to some global lift $C_k$ of $C_{k-1}$,
 so can be $X$, $Y$ and $W$.
Therefore, also in this case, essentially we only need to 
 study the analytic continuation of the pull back of the function  $Z$.
\end{rem}
Recall that
 the function $Z$ is pulled back by the covering map 
 $\widetilde{\Phi}_{\mathcal L}$
 to the function $z^c$ on $\mathfrak X_s$ up to a multiplicative constant, 
 where $c$ is some nonzero integer.
Thus, the function $Z$ is pulled back by 
 $\varphi_0 = \widetilde{\Phi}_{\mathcal L}\circ\varphi_{s, 0}$
 to the function
\[
\left(-\frac{1}{l}(ks+m)\right)^c
\]
 on the component $C_{0, \alpha}$.
Setting the constant term to one as we did before,
 it is expanded as
\[
1+c\frac{k}{m}s + (\text{higher order terms of $s$}),
\]
 around the node corresponding to $s = 0$.

To calculate the obstruction associated to the 
 lift $\varphi_k$ of $\varphi_0$, we need to calculate
 the analytic continuation of the function $Z$ on
 $\mathfrak X$ pulled back by $\varphi_k$
 to 
 a suitable lift $(C_{k, \alpha'}\cup C_{k, \beta'})^*$ 
 of $(C_{0, \alpha'}\cup C_{0, \beta'})^*$.
By the above argument, it is the same as pulling back the function
 $z^c$ by the suitable lift $\varphi_{s, k}$ of $\varphi_{s, 0}$
 such that $\widetilde{\Phi}_{\mathcal L}\circ \varphi_{s, k} = \varphi_k$.
Since the function $z^c$ can be represented using a local parameter
 on $C_{k, L^c}$, now we can perform the analytic continuation
 as in Subsection \ref{subsubsec:std}.
For the precise calculation of the analytic continuation,
 we need to specify the curve $C_k$ and the map $\varphi_k$,
 but as we saw in the calculation in Subsection \ref{subsubsec:std},
 the leading term of the
 analytic continuation can be described only using the data of $\varphi_0$
 (see Remark \ref{rem:leadterm}).

We note the following point which was absent in the standard case.
There are several different pre-log curves in $X_{s, 0}$
 which are mapped to a given pre-log curve of type $(\Xi, h_{\Xi})$
 in $X_0$.
In particular, when the weight of the edge $E$ in Figure \ref{fig:transform}
 in Subsection \ref{subsec:general_edge}
 is $w_E$, then any pre-log curve in $X_{s, 0}$ whose image is given by
\[
k\zeta^{-1} x+lz+m = 0,\;\; y = 0,\;\; w = -\frac{\mu}{\lambda},
\]
\[
\kappa \zeta y+\lambda w+\mu = 0,\;\; x = 0,\;\; z = -\frac{m}{l},
\]
 has the same image by $\widetilde{\Phi}_{\mathcal L}$.
Here $\zeta$ is any $w_E$-th root of unity.
Moreover, these curves naturally correspond to
 the same log structure described in 
 \cite[Proposition 7.1]{NS} (the curves defined by
\[
k\zeta^{-1} x+lz+m = 0,\;\; y = 0,\;\; w = -\frac{\mu}{\lambda},
\]
\[
\kappa \zeta' y+\lambda w+\mu = 0,\;\; x = 0, \;\;  z = -\frac{m}{l},
\] 
 here 
$
\zeta\neq \zeta',
$
 are $w_E$-th roots of unity,
 also have the same image, but correspond to different log structures).
Similarly, the deformed 
 map $\varphi_k$ also has several different lifts in $\mathfrak X_s$,
 see Remark \ref{rem:generalweight} below.

Write $(\varphi_{s, k})^*(x) = S$ and $(\varphi_{s, k})^*(y) = T$. 
Using the relation $ST = t$ 
 and another parameter 
 $(\varphi_{s, k})^*(w) = U_1$ on $C_{k, \beta'}^*$ as before,
 the calculation as in Subsection \ref{subsubsec:std} gives  
 the analytic continuation of the part  
 of the pull back of the function $Z$
 (precisely, the non-principal part 
 at the puncture of $C_{k, \alpha'}^*$ corresponding to 
 the node $C_{0, \alpha'}\cap C_{0, \beta'}$)
 in the following form:
\[
\begin{array}{l}
 \left(c\frac{k\zeta^{-1}}{m}+O(t)\right) S+S^2\xi(S)+1+O(t)\\
 \hspace{.3in} = \left(c\frac{k\zeta^{-1}}{m}+O(t)\right)\frac{t}
 {-\frac{\lambda}{\kappa\zeta}U_1-\frac{\mu}{\kappa\zeta}}
 + t^{2}U_1\xi_1(U_1)+1+O(t).
\end{array}
\]
here $\xi$ is a polynomial and $\xi_1$ is convergent series
 with coefficients in $\Bbb C[t]/t^{k+1}$.
This is expanded at the puncture 
 of $C_{k ,\beta'}^*$
 corresponding to the point $U_1 = 0$ (mod $t$) as
\[
\begin{array}{l}
(\ddagger)\;\;
(1+O(t))tc\frac{k\zeta^{-1}}{m}\cdot\frac{\kappa\zeta}{\mu}
 \cdot \frac{\lambda}{\mu}\cdot{U_1}+tU_1^2\xi_2(U_1)+ t^2U_1\xi_1(U_1)+1+O(t)\\
\hspace{1in} = (1+O(t))tc\frac{k}{m}\cdot\frac{\kappa}{\mu}
 \cdot \frac{\lambda}{\mu}\cdot{U_1} +tU_1^2\xi_2(U_1)+ t^2U_1\xi_1(U_1)+1+O(t),
\end{array}\]
 here $\xi_2$ is a convergent series of $U_1$ with coefficients in 
 $\Bbb C[t]/t^{k+1}$.

There are a few important remarks.
\begin{rem}\label{rem:generalweight}
\begin{enumerate}
\item Although in general 
 there are several lifts $\varphi_{s, k}$ of $\varphi_k$
 as we noted above, they give the same obstruction
 if they correspond to the same log structure.
This can be seen from the calculation above
 for the leading term (admitting that the analytic continuations
 determine the obstruction, the fact which is shown in Subsection \ref{subsec:pairing}).
Of course this is what is to be expected, since the map 
 $\varphi_k = \Phi_{\mathcal L}\circ\varphi_{s, k}$
 does not depend on the choice of the lift $\varphi_{s, k}$.
\item Although the integral length of the edge $E$ is equal to $w_E$,
 the weight of the edge,  
 the leading term is of order one with respect to $t$.
This point seems missing in the former study
 $($\cite{S}$)$.
See also Definition \ref{def:path length}.
\end{enumerate}
\end{rem}
\noindent
This finishes the study of Case (a). \qed\\

Now we study the other case.\\

\noindent
{\bf Case (b).}
Consider an embedded tropical curve $(\Xi, h_{\Xi})$
 with two vertices $\alpha'$ and $\beta'$.
Let $w_E$ be the weight of the edge $E$ connecting these vertices 
 and we assume that the length of the image $h_{\Xi}(E)$
 is also $w_E$ as above.
By Lemma \ref{lem:twovertex},
 it can be obtained from a standard tropical curve with two vertices
 $\Gamma_s$ whose unique bounded edge has length one
 by an integral affine linear map $\mathcal L$.
 

%

Let $\varphi_0\colon C_{0, \alpha'}\cup C_{0, \beta'}\to X_0$
 be a pre-log curve of type $(\Xi, h_{\Xi})$ and 
 $\varphi_k\colon (C_{k, \alpha'}\cup C_{k, \beta'})^*\to \mathfrak X$
 be a $k$-th order lift of it (on a suitable open subset), 
 see Subsection \ref{subsec:general_edge}.
 
The situation we need to consider is the following.  
Assume we have 
 the analytic continuation of functions relevant to the calculation of 
 obstructions (namely, the part of the function
 denoted by $Z$ in Case (a))
 on one of the components (say $C_{k, \alpha'}^*$).
Thus, we have an expression of it
 using a parameter on $C_{k, \alpha'}^*$ as $(\ddagger)$ above.
Then the problem is to obtain the similar expression on the other 
 component $C_{k, \beta'}^*$.

As we argued in Subsection \ref{subsec:general_edge}, 
 the map $\varphi_k\colon (C_{k, \alpha'}\cup C_{k, \beta'})^*\to \mathfrak X$
 factors through a map 
 $\varphi_{s, k}\colon (C_{k, \alpha'}\cup C_{k, \beta'})^*\to \mathfrak X_s$
 so that $\varphi_k = \widetilde{\Phi}_{\mathcal L}\circ\varphi_{s, k}$, 
 where $\widetilde{\Phi}_{\mathcal L}$ is 
 the map $\mathfrak X_s\to \mathfrak X$ induced from $\mathcal L$.

Using this, all the functions on $\mathfrak X$ (like $Z$)
 can be pulled back to functions on $\mathfrak X_s$, 
 where (a deformation of) the curve $\varphi_{s, 0}$
 corresponding to the standard tropical curve $\Gamma_s$ lies.
Since we can express the image of $\varphi_{s, 0}$
 explicitly on $\mathfrak X_s$, it is possible to extend the analytic continuation of 
 the function of the form $(\ddagger)$ by the same calculation as in Case (a).

\begin{rem}\label{rem:Z}
	Again as in Case (a), it will be difficult to calculate the analytic continuation 
	precisely, but the calculation of the leading term is quite possible and it depends 
	only on the data of $\varphi_0$.
\end{rem}

The rest of the calculation is the same as Case (a) and we omit the details.
Note that even if the direction vectors of the edges emanating from $h(\alpha')$
 and $h(\beta')$ span one or two dimensional subspace, 
 the pull back of the components of $\varphi_0(C_0)$ corresponding to these vertices 
 to $\mathfrak X_s$ is written in the same form as in Case (a).
Therefore, the calculation does not depend on 
 the dimension of the spaces spanned by these edges. \qed

\subsubsection{The calculation of the leading term}\label{subsubsec:leadterm}
Now we summarise the calculation so far.
Given a superabundant tropical curve $(\Gamma, h)$ of genus one, 
 we are calculating the pull back of functions on $\mathfrak X$
 to suitable components of $C_k$ through analytic continuation 
 (which will eventually be equivalent to the obstruction associated to 
 a pre-log curve of type $(\Gamma, h)$).
Here $(\Gamma, h)$ is an immersive tropical curve
 (although the argument so far assumed $(\Gamma, h)$
 is an embedding, there is essentially no change required in the case of
 an imbedding).

First, in view of (2) of Remark \ref{rem:generalweight}, 
 we re-define the path length on a tropical curve.
In the following definition, $(\Gamma, h)$
 is an arbitrary tropical curve 
 (in particular, it need not be immersive, 
 and some loops may be contracted.
However, we assume that the abstract graph $\Gamma$
 is a 3-valent graph in the sense of Definition \ref{def:tropical curve}).
\begin{defn}\label{def:path length}
Let $(\Gamma, h)$ be an arbitrary tropical curve defined over integers.
Assume that if $E\subset \Gamma$ is a bounded edge with weight $w_E$, 
 then the integral length $r_E$ of the image $h(E)$ is an integer multiple of $w_E$
 ($r_E$ can be 0).
Let 
\[
\mathcal P\colon v_1\mapsto v_2\mapsto\cdots\mapsto v_m
\]
 be a path in $\Gamma$, 
 here $\{v_i\}$ 
 is a set of consecutive vertices.
Let $E_i$ be the edge connecting $v_i$ and $v_{i+1}$.
Then, we define the length $\ell_{(\Gamma, h)}(\mathcal P)$ of $\mathcal P$ by
\[
\ell_{(\Gamma, h)}(\mathcal P)=\sum_{i=1}^{m-1}\frac{r_{E_i}}{w_{E_i}}.
\]
We extend the definition of $\ell_{(\Gamma, h)}(\mathcal P)$
 to arbitrary
 tropical curves defined over $\Bbb R$ in the obvious way.
\end{defn}
Note that when $(\Gamma, h)$ is defined over integers
 as above,
 then $\ell_{(\Gamma, h)}(\mathcal P)$ 
 is an integer by the assumption
 about $r_E$ and $w_E$.
We assume this condition for $(\Gamma, h)$ if not specified.
In fact, when $(\Gamma, h)$ is immersive, 
 we usually assume that by adding 2-valent vertices, 
 the length $r_E$ equals to $w_E$ for each bounded edge $E$ of $\Gamma$.
\begin{rem}\label{rem:length}
Note that the length is defined on a path in $\Gamma$.
In particular, there can be different paths $\mathcal P$ and $\mathcal P'$
 such that their images $h(\mathcal P)$ and $h(\mathcal P')$ are the
 same, but there lengths are different.
\end{rem}

Now
 by the calculation so far, one observes the following.
Let $\mathcal P$ be the unique shortest path from
 a 1-valent vertex $\alpha_0 = \alpha$ of $\overline\Gamma'$ 
 to the loop
 of $\Gamma$, see Subsection \ref{subsec:convention}
 and Figure \ref{fig:Gamma'}.
Let $\alpha_N$ be a vertex on $\mathcal P$ and 
 let $\alpha, \alpha_1, \dots, \alpha_N$
  be the vertices between $\alpha$ and $\alpha_N$. 
Let $q_{i, i+1}$ be the node between
 $C_{0,\alpha_i}$ and $C_{0, \alpha_{i+1}}$
 (in particular, $q_{0, 1} = p$ in the notation of Lemma \ref{lem:coordinate}).
For each $\varphi_{0}(q_{i, i+1})$,
 we can find a pre-log curve 
\[
\psi_{0, i}\colon C_{0, \alpha_i}\cup C_{0, \alpha_{i+1}}\to 
 \mathfrak X_{s, i}
\]
 of type $\Gamma_{s, i}$ in a suitable toric
 degeneration $\mathfrak X_{s, i}$ such that its image is mapped to 
 the union of curves $\varphi_0(C_{0, \alpha_i}\cup C_{0, \alpha_{i+1}})$
 by a toric morphism from a subset of $\mathfrak X_{s, i}$ to $\mathfrak X$
 as in Cases (a) and (b) above (see Remark \ref{rem:split} below for more details).
Here $\Gamma_{s, i}$ is a standard tropical curve with two vertices
 and $\mathfrak X_{s, i}$ is a toric degeneration defined respecting it.
\begin{rem}\label{rem:split}
Usually there will  not be a toric morphism from the whole 
 $\mathfrak X_{s, i}$ or a blow up of it to $\mathfrak X$.
Recall that $\mathfrak X_{s, i}$ is defined by 
 a fan obtained as the closure of the cones over 
 a polyhedral decomposition $\mathcal P$ of $\Bbb R^n$.
The unbounded regions of $\mathcal P$ correspond 
 to families of toric strata over the base $\Bbb C$.
Removing these strata from $\mathfrak X_{s, i}$, there will be a toric morphism
 from it to $\mathfrak X$.
Since the image of this region in $\mathfrak X$
 contains  
 $\varphi_k((C_{k, \alpha_i}\cup C_{k, \alpha_{i+1}})^*)$,
 restricting the map to it
 suffices for our purpose.	
\end{rem}

Moreover, if $\varphi_{k}$ is a lift of $\varphi_0$ over $\Bbb C[t]/t^{k+1}$, 
 then the map $\psi_{0, i}$ also lifts on the open subset 
 $(C_{k, \alpha_i}\cup C_{k, \alpha_{i+1}})^*$
 in the way that the composition of it with the above toric morphism is
 $\varphi_k|_{(C_{k, \alpha_i}\cup C_{k, \alpha_{i+1}})^*}$.
We write these lifts by $\psi_{k, i}$.

For each $i$, there is a set of functions
 $\{x_1^{(i)}, y_1^{(i)}, x_{2}^{(i)}, y_2^{(i)}, z_{1}^{(i)}, \dots, z_{n-3}^{(i)}\}$  
  on 
 $\mathfrak X_{s, i}$ around $\psi_{0, i}(q_{i, i+1})$ 
 and coordinates on the components 
 $C_{0, \alpha_i}$ and $C_{0, \alpha_{i+1}}$ 
 with the same properties as in Lemma \ref{lem:coordinate}.
Especially, the set of functions
 should satisfy the following properties.
\begin{itemize}
\item For each $i = 0, 1, \dots, N-1$,
 the image $\psi_{0, i}(C_{0, \alpha_i})$
 is defined by the equations of the form
\[
k_iy_1^{(i)}+l_ix_1^{(i)}+m_i = 0,\;\; y_2^{(i)} = 0,\;\;
 x_2^{(i)} = -\frac{m_{i+1}}{l_{i+1}},\;\;
 \;\; z_{1}^{(i)} = c_{i, 1},\dots, z_{n-3}^{(i)} = c_{i, n-3},
\]
 and the image $\psi_{0, i}(C_{0, \alpha_{i+1}})$
 is defined by the equations of the form
\[
k_{i+1}y_2^{(i)}+l_{i+1}x_2^{(i)}+m_{i+1} = 0, \;\;
 x_1^{(i)} = 0,\;\;
 y_1^{(i)} = -\frac{m_{i}}{k_{i}},
 \;\; z_{1}^{(i)} = c_{i, 1},\dots, z_{n-3}^{(i)} = c_{i, n-3},
\]
 around the node $\psi_{0, i}(q_{i, i+1})$.
Here $k_j, l_j, m_j$ and $c_{i, j}$ are non-zero constants.
\item In this notation, the node $\psi_{0, i}(q_{i, i+1})$ corresponds to 
 $x_1^{(i)} = 0$ on $\psi_{0, i}(C_{0, \alpha_i})$
 and $y_2^{(i)} = 0$ on $\psi_{0, i}(C_{0, \alpha_{i+1}})$.
Moreover, the node $\psi_{0, i-1}(q_{i-1, i})$ corresponds to $y_1^{(i)} = 0$
 and the node $\psi_{0, i+1}(q_{i+1, i+2})$ corresponds to $x_2^{(i)} = 0$.
\item The functions $x_1^{(i)}$ and $y_{2}^{(i)}$ satisfy 
\[
x_1^{(i)}y_{2}^{(i)}  =t, \;\;
 i= 0, \dots, N-1.
\]
\end{itemize}

With these notations, we can state the following proposition.
Namely, repeating
 the calculation of the analytic continuation of the part of the pull back of 
 the function $Z$ (as in Subsection \ref{subsubsec:nonstd})
 until we reach the loop, 
 we obtain the terms of the form
\[
1+O(t)+ \frac{\chi(t) t^M}{V_N}+ \frac{\chi_1(t) t^{M_1}}{V_N^2}
 +\frac{\chi_2(t) t^{M_2}}{V_N^3}+
 \cdots
\]
 on the open subset $C_{k, \alpha_N}^*$ of $C_{k}$
 around the puncture corresponding to the edge connecting
 $\alpha_{N-1}$ and $\alpha_N$.
Here $V_N$ is the parameter on $C_{k, \alpha_N}^*$
 given by $\psi_{k, N-1}^*(y_2^{(N)})$
 and $\chi(t), \chi_1(t), \chi_2(t), \cdots$
  are series of $t$ with non-negative
 exponents.
Moreover, $M, M_i, i = 1, 2, \dots$ are positive integers satisfying
\[
M<M_1<M_2<\cdots.
\]
Then the constant term of $\chi(t)$, which will give the leading term
 of the analytic continuation,
 and the exponent $M$ of $t$
 are given by the following.
\begin{prop}\label{prop:obstweight}
The integer $M$ is equal to the length $\ell_{(\Gamma, h)}(\mathcal P)$
 of the path $\mathcal P$ connecting 
 $\alpha$ and $\alpha_N$ in the sense of Definition \ref{def:path length}.
The constant term of the 
 coefficient $\chi(t)$ is given by
\[
\frac{l_0}{m_0}
 \cdot\frac{k_1}{m_1}\cdot \frac{l_1}{m_1}
 \cdot\cdots\cdot
 \frac{k_{N-1}}{m_{N-1}}
 \cdot \frac{l_{N-1}}{m_{N-1}}. 
\]
 up to a constant which is bounded by a constant
 $C_{(\Gamma, h)}$ independent of the path $\mathcal P$.
Here $k_i, l_i, m_i$ are coefficients of the defining equations of the pull back of the 
 component $\varphi_0(C_{0, \alpha_i})$ introduced above
 (see Remarks
 \ref{rem:coeff} and \ref{rem:dependence}).\qed
 \end{prop}
The bounded constant in the statement is due to the 
 factor $c$ appeared in the calculation in Case (a) above.

\begin{rem}\label{rem:total}
As we noted at the beginning of Subsection \ref{subsubsec:nonstd},
 our goal is to describe the pull back 
 of the function $Z$ by $\varphi_k$ to the components
 $C_{k, \alpha_i}^*$.
The term described in Proposition \ref{prop:obstweight}
 gives a part of the leading term of it.
The actual leading term is given by the sum of similar terms
 contributed from the 1-valent vertices $p_j$ of $\overline\Gamma'$
 satisfying $\alpha_N>p_j$ in the notation 
 introduced before Lemma \ref{lem:prinpart}.
\end{rem}

\begin{rem}\label{rem:coeff}
As we noted in Remark \ref{rem:coord}, individual coefficients
 $k_i, l_i, m_i$ depend on the choice of 
 functions in Lemma \ref{lem:coordinate}.
Such a choice can change coordinates $x_j^{(i)}$ and $y_j^{(i)}$ 
 by multiplying non-zero constants, or  coordinates $z_{i, j}$
 which are constant on 
 $\psi_{0, i}((C_{0,\alpha_i}\cup C_{0, \alpha_{i+1}})^*)$.
However, to preserve the relation $x_1^{(i)}y_2^{(i)} = t$, 
 we cannot change these coordinates independently.
It follows that the product of the pair of factors 
 $\frac{l_i}{m_i}\cdot \frac{k_{i+1}}{m_{i+1}}$ does not change
 under such a coordinate change.
The last factor $\frac{l_{N-1}}{m_{N-1}}$ depends on the choice of 
 the coordinates.
Also, the individual coefficients depend on the choice of a lift on 
 each open cover 
 (see the argument before Remark \ref{rem:generalweight}).
However, such dependence eventually vanishes, too.
See Remark \ref{rem:dependence}. 
\end{rem}

\noindent
\subsection{Calculation of obstructions in embedded cases. Step 5: Obstruction at the loop.}\label{subsec:step5}
\begin{figure}[h]
\includegraphics{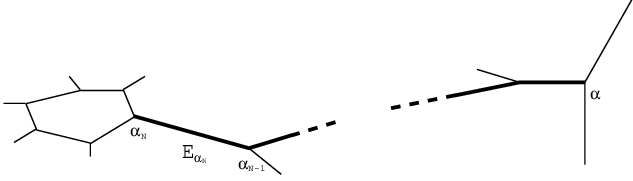}
\caption{The path $\mathcal P$ is drawn by bold lines.}\label{fig:residue}
\end{figure}

We use the notations in Subsection \ref{subsec:step4}.
Let $\alpha_N$
 be the vertex on the loop $L$ nearest to $\alpha$ (it is determined
 uniquely).
Let $C_{0, \alpha_N}$ be the component of 
 $C_0$ corresponding to $\alpha_N$.
Let $E_{\alpha_N}$ be the edge attached to $\alpha_N$ which is 
 the last edge in the path $\mathcal P$ from $\alpha$ to $\alpha_N$
 (see Figure \ref{fig:residue}).
We assume that we have a global $(k-1)$-th order lift 
 $\varphi_{k-1}\colon C_{k-1}\to \mathfrak X$ of $\varphi_0$, 
 and also have its $k$-th order lift $\varphi_{k}$
 on the union of the closure of the connected
 component of $C_{0, \Gamma}\setminus C_{0, L}$
 containing $C_{0, \alpha}$ and the open subset $C_{0, \alpha_N}^*$. 
Note that since this union is an open subset of a tree of rational curves,
 such a lift exists as in Lemma \ref{lem:liftT}.
Let $V_{N}=\psi_{k, N-1}^*(y_2^{(N-1)})$
 be the affine coordinate on $C_{k, \alpha_{N}}^*$.

Now fix a basis $\{\mathfrak v_1, \dots, \mathfrak v_{n-r}, 
 \mathfrak w_1, \dots, \mathfrak w_{r}\}$ of $N$
 as in Lemma \ref{lem:trivsub}.
Let $\{\mathfrak v_1^{\vee}, \dots, \mathfrak v_{n-r}^{\vee}, 
\mathfrak w_1^{\vee}, \dots, \mathfrak w_{r}^{\vee}\}$
 be its dual basis.
Then, by Proposition \ref{prop:obstweight} and Remark \ref{rem:total},
 the pull back 
 of the function on $\mathfrak X$ corresponding to $\mathfrak w_1^{\vee}$
 (see Remark \ref{rem:function} below for the precise meaning)
 by $\varphi_k$ to  
 $C_{k, \alpha_{N}}^*$ has an expansion 
 at the puncture corresponding to $V_N = 0$ (mod $t$)
 whose principal part is the finite series of
 the form
\[
(1+O(t))t^{\ell_{(\Gamma, h)}(\mathcal P)}
 \left(C\cdot \frac{k_0}{m_0}\cdot\frac{l_1}{m_1}
 \cdot \frac{k_1}{m_1}\cdot\cdots\cdot
 \frac{l_{N-1}}{m_{N-1}}\cdot \frac{k_{N-1}}{m_{N-1}}\cdot 
 \frac{1}{V_N}+ \xi\left(\frac{t}{V_N}\right) \right)+ 
 \sum_{i=2}^at^{\ell_{(\Gamma, h)}(\mathcal P_i)}\frac{1}{V_N}
 \xi_i\left(\frac{t}{V_N}\right),
\] 
 where 
 $\mathcal P_i$, $i\geq 2$ is the path from
 the 1-valent vertex $p_i$ of $\overline\Gamma'$ to $\alpha_{N}$
 such that $\alpha_{N}>p_i$ and $p_i\neq \alpha$
 (see  the paragraph before Lemma \ref{lem:prinpart} for the notation).
Also, $\xi$ and $\xi_i$ are polynomials, where $\xi$ does not have the constant term.
Moreover, $C$ is a bounded constant which is determined from the
 data of $(\Gamma, h)$.
See the comment
 after Proposition \ref{prop:obstweight}.
The vector $\mathfrak w_1$ is chosen so that the subspace of $N_{\Bbb R}$
 spanned by $\bar A$ and the directions of the edges emanating 
 from the vertex $p_1 = \alpha$ coincides with the 
 subspace spanned by $\bar A$ and $\mathfrak w_1$
 ($p_1 = \alpha$ is taken so that the length $\ell_{(\Gamma, h)}(\mathcal P)$
 is (one of) the shortest among the 1-valent vertices in the same connected component
 of $\overline\Gamma'\setminus L$).

\begin{rem}\label{rem:function}
Recall that the degeneration $\mathfrak X$ is defined
	by a fan in $(N\oplus\Bbb Z)\otimes\Bbb R$
	having nonnegative part with respect to $\Bbb Z\otimes\Bbb R$.
	Let $\mathfrak u$ be the positive generator of the $\Bbb Z$-summand and
	consider the dual basis $\mathfrak v_1^{\vee}, \dots, \mathfrak v_{n-r}^{\vee}, 
	\mathfrak w_1^{\vee}, \dots, \mathfrak w_{r}^{\vee}, \mathfrak u^{\vee}$ of the basis
	$\mathfrak v_1, \dots, 
	\mathfrak v_{n-r}, \mathfrak w_1, \dots, \mathfrak w_{r}, \mathfrak u$ 
	of $N\oplus\Bbb Z$ (see Lemma \ref{lem:trivsub} and the paragraph before it for the notation).
	Then we take a function on $\mathfrak X$ corresponding to 
	the vector $\mathfrak w^{\vee}+j\mathfrak u^{\vee}$ 
	for a unique integer $j$ so that
	it is not constantly 0 or $\infty$ on the component $X_{0, \alpha}$
	of the central fiber $X_0$ of $\mathfrak X$ corresponding to the 
	vertex $\alpha$.
	We call this function the function corresponding to 
	$\mathfrak w^{\vee}$, with the vertex $\alpha$ understood.
\end{rem}

Recall (see Subsection \ref{subsec:abundance})
 that
 the de Rham representatives of the dual of the
 cohomology group 
 $H^1(C_0, \mathcal N_{C_0/X_0})^{\vee}
  (\cong H^0(C_0, \mathcal N_{C_0/X_0}^{\vee}
   \otimes\omega_{C_0}))$
 in which the dual classes of the obstruction lie
 are given by vector valued logarithmic differential forms
 whose differential form parts are of the form
 $\frac{dU}{U}$ on each component of the loop of rational curves.
Here $U$ is a coordinate
 which takes the value $0$ or $\infty$ at the intersection with the other components
 of the loop.
Therefore, we have to use such a coordinate instead of 
 $V_N$
 when we calculate the obstruction.
 
Using the notations in Subsection \ref{subsec:step4},
 the component $\psi_{0, N-1}(C_{0, \alpha_N})$ 
 corresponding to the vertex
 $\alpha_N$ is defined by the equations
\[
k_Ny_2^{(N-1)}+l_Nx_2^{(N-1)}+m_N = 0,\;\; x_1^{(N-1)} = 0,\;\; 
 y_1^{(N-1)} = -\frac{m_{N-1}}{k_{N-1}},\;\;
 z_{1}^{(N-1)} 
  = c_{N-1, 1},\dots, z_{n-3}^{(N-1)} = c_{N-1, n-3}.
\]
The function $y_2^{(N-1)}$
 is pulled back to $V_N$ by $\psi_{k, N-1}$ as above.
We use the same letter $V_N$ to represent $\psi_{0, N-1}^*(y_2^{(N-1)})$
 (which is the reduction of $V_N$ to a function over $\Bbb C[t]/t$)
 for notational simplicity.
Then $x_2^{(N-1)}$ is pulled back to another affine coordinate
\[
U_N = -\frac{1}{l_N}(k_NV_N+m_N)
\]
 by $\psi_{0, N-1}$.
Note that this coordinate satisfies the conditions imposed on the
 particular coordinate $U$ mentioned above which is suitable
 for the calculation 
 of the obstruction.
Then when we rewrite the term $\frac{1}{V_N}$ using $U_N$, we have
\[
\frac{1}{V_N} = -\frac{k_N}{l_NU_N+m_N} =
 - \frac{k_N}{m_N}\cdot\frac{1}{1+\frac{l_N}{m_N}U_N}.
\]
Thus, it has the first order pole at $U_N = -\frac{m_N}{l_N}$
 with residue $\frac{k_N}{m_N}$
 (more precisely, the 1-form
 $- \frac{k_N}{m_N}\cdot\frac{1}{1+\frac{l_N}{m_N}U_N}
 \frac{dU_N}{U_N}$ has the 
 residue $\frac{k_N}{m_N}$ at $U_N = -\frac{m_N}{l_N}$).
When we use the function $\psi_{k, N-1}^*(x_2^{(N-1)})$ instead of 
 $\psi_{0, N-1}^*(x_2^{(N-1)})$, the function
 $- \frac{k_N}{m_N}\cdot\frac{1}{1+\frac{l_N}{m_N}U_N}$ 
 is perturbed by terms of positive order with respect to $t$.

We summarize the calculation so far.
Fix an integral linear combination $\mathfrak w^{\vee}$ 
of the vectors	$\mathfrak w_1^{\vee}, \dotsm \mathfrak w_r^{\vee}$.
Let $Z$ be the function on $\mathfrak X$ corresponding to $\mathfrak w^{\vee}$
 in the sense of Remark \ref{rem:function}.
We normalize it so that it is the constant function 
 with the value one when it is pulled back to
 the trivial lifts of $\varphi_0|_{C_{0, \alpha}^*}$ for any vertex $\alpha$ in $\Gamma'$.
\begin{rem}\label{rem:normalization}
In terms of the explicit correspondence between
 the set of local sections of the log normal sheaves and the set of local lifts
 given in Lemmas \ref{lem:lift3h} and \ref{lem:lifthigh}, 
 the pull back of a function corresponding to $\mathfrak w_j^{\vee}$
 to the trivial lifts of $C_{0, \alpha}^*$ has the constant value 
 which is the product of some of the constants $a_{i}$.
This value does not depend on the choice of a vertex in $\Gamma'$
 nor the order of the lifts. 
Therefore, the above normalization is well-defined and
 can be simultaneously fixed
 for functions corresponding to any $\mathfrak w^{\vee}$. 
\end{rem}
Let $p_1 = \alpha, p_2, \dots, p_a$ 
 be the 1-valent vertices of $\overline{\Gamma}'$
 which satisfy $\alpha_N>p_i$
 (see the paragraph before Lemma \ref{lem:prinpart} and Figure \ref{fig:0z}).
 
Recall that we are assuming we have a $(k-1)$-th order lift $\varphi_{k-1}$
 of $\varphi_0$.
By Lemma \ref{lem:canlift}, the restriction $\varphi_{k-1}|_{(C_{k-1, \alpha_N}\cup
 C_{k-1, \alpha_{N-1}})^*}$ has a $k$-th order lift (we may take the trivial lift
 as in Lemma \ref{lem:canlift}, but any lift will do).
Let us write such a lift as $\varphi_{k}|_{(C_{k, \alpha_N}\cup
	C_{k, \alpha_{N-1}})^*}$.
Pulling back the function $Z$ by $\varphi_{k}|_{(C_{k, \alpha_N}\cup
	C_{k, \alpha_{N-1}})^*}$,
 we obtain a function defined over $\Bbb C[t]/t^{k+1}$.

\begin{rem}
For the calculation of the obstruction, it suffices to know
 only the lift on ${(C_{k-1, \alpha_N}\cup
 	C_{k-1, \alpha_{N-1}})^*}$, 
 since the dual obstruction classes are supported on the loop
 (see Remark \ref{rem:suppob}).
On the other hand, in the explicit calculation above, we assumed a lift of $\varphi_{k-1}$
 in a larger subset of $C_{k-1}$.
This is because directly calculating the pull back $\varphi_k^*Z$
 of a suitable function $Z$ on $\mathfrak X$ to $(C_{k, \alpha_N}\cup
 	C_{k, \alpha_{N-1}})^*$ will be difficult.
We can effectively calculate it only through an appropriate analytic continuation.

The lift $\varphi_{k}|_{(C_{k, \alpha_N}\cup
 C_{k, \alpha_{N-1}})^*}$ in the paragraph before this remark
 can be either the restriction of 
 this given lift,
  or a trivial lift of $\varphi_{k-1}|_{(C_{k-1, \alpha_N}\cup
 	C_{k-1, \alpha_{N-1}})^*}$ as in that paragraph, or even any other lift.
The pull back $\varphi_k^*Z$ will change 
 according to these choices, but the obstruction does not depend on 
 them.
\end{rem}

\begin{defn}\label{def:preob}
The \emph{pre-obstruction at $\alpha_N$
 of the direction $\mathfrak w^{\vee}$}
 is the principal part of the function $Z$ 
 expanded using the coordinate $U_N$ on $C_{k, \alpha_N}^*$
 at the puncture
 $U_N = -\frac{m_N}{l_N}$.
We call the part of the
 residue of $Z\frac{dU_N}{U_N}$ 
 at that point which is of order $k$ with respect to $t$
 the \emph{obstruction at $\alpha_N$
	of the direction $\mathfrak w^{\vee}$}.
We write it as $o(\mathfrak w^{\vee}, \alpha_N; \varphi_{k-1})$.
Furthermore, we define the \emph{obstruction of the direction $\mathfrak w^{\vee}$}
 as the sum
\[
\sum_{\gamma}o(\mathfrak w^{\vee}, \gamma; \varphi_{k-1}),
\]
 where $\gamma$ runs over the set of vertices on the loop $L$ of $\Gamma$.
We write it by $o(\mathfrak w^{\vee};\varphi_{k-1})$.
\end{defn}

Now we have two definitions of the obstruction, Definitions  \ref{def:cechobst} 
 and \ref{def:preob}.
Definition \ref{def:cechobst} is the one from the general theory, but
 difficult to calculate in general.
The one in Definition \ref{def:preob} is specialized to the degenerate
 situation, and more calculable.
We will see these are equivalent in the following subsections.

Recall that the pre-obstruction at $\alpha_N$ (in any direction)
 is the sum of contributions from $p_1, \dots, p_a$ (see Lemma \ref{lem:prinpart}).
Namely, let $F_i$ be the edge of $\overline{\Gamma}'$ one of whose
 ends is $p_i$, and take the non-principal part of $Z$
 at the puncture of $C_{k, p_i}^*$ corresponding to the edge $F_i$.
Then its contribution to the pre-obstruction at $\alpha_N$ 
 of direction $\mathfrak w^{\vee}$ is given by the principal part of the
 analytic continuation of this non-principal part expanded using 
 the coordinate $U_N$ as in Definition \ref{def:preob}.

\begin{defn}\label{def:preobcontri}
We call this contribution to the pre-obstruction the
 \emph{pre-obstruction at $\alpha_N$ 
 	of direction $\mathfrak w^{\vee}$
 	contributed from $p_i$}.
Also, we call the associated residue the
 \emph{obstruction at $\alpha_N$ 
 	of direction $\mathfrak w^{\vee}$
 	contributed from $p_i$}.
We write it as $o(\mathfrak w^{\vee}, \alpha_N, p_i; \varphi_{k-1})$.
\end{defn}

By definition, $o(\mathfrak w^{\vee}, \alpha_N; \varphi_{k-1})
 = \sum_{i=1}^a o(\mathfrak w^{\vee}, \alpha_N, p_i; \varphi_{k-1})$.
Let us take $\mathfrak w_1^{\vee}$ and $\alpha$ as above.
We have the following from the above calculation.

\begin{prop}\label{prop:obstruction1}
The pre-obstruction at $\alpha_N$
 of direction $\mathfrak w_1^{\vee}$ contributed from
 $\alpha$ is given by
\[
\begin{array}{ll}
(\spadesuit) \hspace{.4in}
-(1+O(t))t^{\ell_{(\Gamma, h)}(\mathcal P)}
  & C\cdot \frac{l_0}{m_0}\cdot\frac{k_1}{m_1}
 \cdot \frac{l_1}{m_1}\cdot\cdots\cdot
 \frac{k_{N-1}}{m_{N-1}}\cdot \frac{l_{N-1}}{m_{N-1}}\cdot 
 \frac{k_N}{m_N}\cdot\frac{1}{1+\frac{l_N}{m_N}U_N}\\
 & +t^{\ell_{(\Gamma, h)}(\mathcal P)-1}
  \left(\frac{t}{1+\frac{k_N}{m_N}U_N}\right)^2
 O(t)\xi_2\left(\frac{t}{1+\frac{k_N}{m_N}U_N}\right),
\end{array}
\]
 here $\xi_2$ is a polynomial and $C$ is a bounded constant.
In particular, when $k = \ell_{(\Gamma, h)}(\mathcal P)$, the obstruction at $\alpha_N$
 of direction $\mathfrak w_1$ contributed from
 $\alpha$ is of the form
\[
t^{\ell_{(\Gamma, h)}(\mathcal P)}
 C\cdot \frac{l_0}{m_0}\cdot\frac{k_1}{m_1}
\cdot \frac{l_1}{m_1}\cdot\cdots\cdot
\frac{k_{N-1}}{m_{N-1}}\cdot \frac{l_{N-1}}{m_{N-1}}\cdot 
\frac{k_N}{m_N}
\]
 and when $k<\ell_{(\Gamma, h)}(\mathcal P)$, it is zero.
\qed
\end{prop}
\begin{rem}\label{rem:dependence}
As we noted in Remark \ref{rem:coord}, individual coefficients
 $k_i, l_i, m_i$ depend on the choices of 
 functions in Lemma \ref{lem:coordinate}.
However, as we claimed in Remark \ref{rem:coeff}, the overall coefficient
 $\frac{l_0}{m_0}\cdot\frac{k_1}{m_1}
 \cdot \frac{l_1}{m_1}\cdot\cdots\cdot
 \frac{k_{N-1}}{m_{N-1}}\cdot \frac{l_{N-1}}{m_{N-1}}\cdot 
 \frac{k_N}{m_N}$
 does not depend on the choice.
\end{rem}
\begin{defn}\label{def:leadingterm}
In the expression $(\spadesuit)$
 of Proposition \ref{prop:obstruction1},
 we call the term which has the 
 lowest order with respect to $t$ the \emph{leading term}.
\end{defn}
Thus, the leading term of $(\spadesuit)$ is 
 $-t^{\ell_{(\Gamma, h)}(\mathcal P)}C\cdot\frac{l_0}{m_0}\cdot\frac{k_1}{m_1}
 \cdot \frac{l_1}{m_1}\cdot\cdots\cdot
 \frac{k_{N-1}}{m_{N-1}}\cdot \frac{l_{N-1}}{m_{N-1}}\cdot 
 \frac{k_N}{m_N}\cdot\frac{1}{1+\frac{l_N}{m_N}U_N}$.
\begin{cor}\label{obstorder}
The leading term of the pre-obstruction 
 is determined by the configuration of the image $\varphi_0(C_0)$.
In particular, when $\varphi_k$ and $\varphi_k'$ 
 are two lifts of $\varphi_0$,
 the pre-obstructions at $\alpha_N$ of direction $\mathfrak w_1$
 contributed from 
 $\alpha$
 associated to them differ only 
 in higher order terms with respect to $t$.\qed
\end{cor}

\subsubsection{Residue and obstruction}\label{subsec:resobs}
In this subsection, we clarify the relation between the obstruction to lift
 the map $\varphi_k$ and the above calculation of analytic continuations of functions.
Recall that the obstruction is defined as a \v{C}ech 1-cocycle
 constructed using a suitable open cover of $C_k$ 
 (see Definition \ref{def:cechobst}).
In general, it is not easy to see directly that a given \v{C}ech cocycle is
 cohomologically trivial or not.
In our case, we know very well about the dual of the obstruction cohomology classes
 (see Corollary \ref{cor:obstbasis}).
Thus, we can see whether the given cocycle is zero or not by computing its
 pairing with a basis of the dual classes.

The pairing is defined as follows.
Let $\{C_{0,T_i}^*, C_{0,\beta\gamma}\}$ be an open cover of $C_0$.
Here $T_i$ is a part of the graph $\Gamma$ introduced in 
 Subsection \ref{subsec:step3} (see Figure \ref{fig:T}), and
 $C_{0, T_i}^*$ is the union of corresponding irreducible components of $C_0$ 
 with the node corresponding to the edge $E_i'$ in 
 Figure \ref{fig:T} removed.
Also, $\beta$ and $\gamma$ are adjacent vertices of $\Gamma$ at least
 one of which is not contained in any $T_i$, and 
 $C_{0,\beta\gamma}$ is the open subset $(C_{0, \beta}\cup C_{0, \gamma})^*$
 in the notation of the previous sections.

On the other hand, let us take locally closed subsets of $C_0$ as follows.
First, take $C_{0, T_i}^*$ as above. 
Second, for any $C_{0, \beta\gamma}$, let $p_{\beta\gamma}$ be the unique node.
Then take $C_{0, \beta}^*\cup p_{\beta\gamma}$ and 
 $C_{0, \gamma}^*\cup p_{\beta\gamma}$.
We write these locally closed subsets by $C_{0, \beta\gamma, 1}$ and
 $C_{0, \beta\gamma, 2}$, respectively.
Also, we write by $\{Y_l\}$ the set of these locally closed subsets.
By construction, to each $Y_l$, we can naturally associate the unique
 open subset from the covering above
 (in particular, to $C_{0, \beta\gamma, i}$, we associate $C_{0, \beta\gamma}$).
We write this open subset by $\tilde Y_l$.

Now take a section $\mu_l$ of the sheaf
 $\mathcal N_{C_0/X_0}|_{Y_l}(\infty p_{\beta\gamma}) = 
 (\varphi_0^*\Theta_{X_0}/\Theta_{C_0})|_{Y_l}(\infty p_{\beta\gamma})$
 restricted to each locally closed subset $Y_l$ of the form $C_{0, \beta\gamma}$.
That is, $\mu_l$ is a section of $\mathcal N_{C_0/X_0}$ on $Y_l$
which can have a pole of any order at the point $p_{\beta\gamma}$.
When $Y_l$ is of the form $C_{0, T_i}^*$, then take $\mu_l$ to be a section
 in $\Gamma(C_{0, T_i}^*, \mathcal N_{C_0/X_0})$.
Let $Y_1$ and $Y_2$ be two of these locally closed subsets whose
 intersection is an open subset (note that in this case $\tilde Y_1\neq \tilde Y_2$).
Then associate the section $\mu_1-\mu_2$ of $\mathcal N_{C_0/X_0}$ 
 to the open subset 
 $Y_1\cap Y_2$ (the order also matters. That is, we associate $\mu_2-\mu_1$
 to $Y_2\cap Y_1$).
Note that we have
\[
Y_1\cap Y_2 = \tilde Y_1\cap \tilde Y_2.
\]
Therefore, associating $\mu_1-\mu_2$ to $\tilde Y_1\cap \tilde Y_2$, 
 the set of sections $\{\mu_l\}$ determines a \v{C}ech 1-cocycle 
 with values in $\mathcal N_{C_0/X_0}$
 for the covering $\{Y_l\} = \{C_{0, T_i}^*\}\cup \{C_{0, \beta\gamma}\}$ above.

Conversely, any class in $H^1(C_0, \mathcal N_{C_0/X_0})$ can be represented 
 in this way (this is a consequence of Lemma \ref{lem:nondegeneracy} below).

Now we define the pairing between $H^1(C_0, \mathcal N_{C_0/X_0})$
 and $H = H^0(C_0, \mathcal N_{C_0/X_0}^{\vee}\otimes\omega_{C_0})$.
First, take an element $f\otimes \xi$ of the space $H$, where
 $f\in (\bar A_{\Bbb C})^{\perp}$ and $\xi$ is the generator of 
 the space $H^0(C_0, \omega_{C_0})$, see Corollary \ref{cor:obstbasis} and 
 the paragraph before it.
Take a family of sections $\{\mu_l\}$ as above.
Let $Y_l$ be a locally closed subset which is contained in $C_{0, \Gamma'}$.
Then the log normal sheaf splits as 
 $\mathcal N_{C_0/X_0}|_{Y_l}\cong
  \mathcal N'|_{Y_l}\oplus \mathcal O_{Y_l}^{\oplus r}$
 as in Subsection \ref{subsec:step3}.

Since each $\mu_l$ has the value in the sheaf $\mathcal N_{C_0/X_0}$,
 the $(\bar A_{\Bbb C})^{\perp}$-factor $f$ of $f\otimes \xi$ couples with it,
 resulting in a meromorphic 1-form on $Y_l$
 (recall $\Theta_{X_0}$ is naturally identified with $N\otimes \mathcal O_{X_0}$.
 Since $\Theta_{C_0}$ considered as a subsheaf of $\Theta_{X_0}$
 is annihilated by $(\bar A_{\Bbb C})^{\perp}\subset (N_{\Bbb C})^{\vee}$, this coupling is well
 defined).
Then we define the pairing by taking the sum of the residues of the poles of 
 these sections:
\[
(\heartsuit) \hspace{.5in}(f\otimes \xi, \{\mu_l\}) = \sum_lres_q(f\otimes\xi, \mu_l),
\]
 here $q$ runs through all the poles of $(f\otimes\xi, \mu_l)$ on $Y_l$.
We note the following.
\begin{rem}\label{rem:suppob}
The pairing $(f\otimes\xi, \mu_l)$ is identically zero unless 
 $Y_l$ is contained in the loop part $C_{0, L}$ of $C_0$, 
 since outside $C_{0, L}$, the 1-form $\xi$ is identically zero.
\end{rem}

One can check the following (see \cite{NY} for the full detail).
\begin{lem}
The definition $(\heartsuit)$ gives a well-defined pairing between 
 $H^1(C_0, \mathcal N_{C_0/X_0})$
 and $H$. \qed
\end{lem}
The well-definedness means the following two properties:
\begin{itemize}
\item First, different families of sections $\{\mu_l\}$ and $\{\mu'_l\}$ may give
 the same \v{C}ech 1-cocycle.
In this case, the definition $(\heartsuit)$ gives the same value.
That is, $(\heartsuit)$ is well-defined at the level of cocycles.
\item Second, $(\heartsuit)$ only depends on the cohomology class of 
 the \v{C}ech 1-cocycle.
\end{itemize}

Moreover, the following is also easy to see.
\begin{lem}\label{lem:nondegeneracy}
The pairing $(\heartsuit)$ is nondegenerate.
\end{lem}
\proof
Fix a basis $\{\mathfrak v_1, \dots, \mathfrak v_{n-r}, 
 \mathfrak w_1, \dots, \mathfrak w_{r}\}$ of $N$
 and its dual basis
 $\{\mathfrak v_1^{\vee}, \dots, \mathfrak v_{n-r}^{\vee}, 
 \mathfrak w_1^{\vee}, \dots, \mathfrak w_{r}^{\vee}\}$
 as in the beginning of Subsection \ref{subsec:step5}.
Then $\{\mathfrak w_1^{\vee}\otimes\xi, \dots, \mathfrak w_{r}^{\vee}\otimes\xi\}$
 gives a basis of the space $H$.

Now take the family $\{\mu(i)_l\}$, $i = 1, \dots, r$ as follows.
Namely, let $\beta, \gamma$ be a neighboring vertices on the loop of $\Gamma$.
We have associated locally closed subsets
 $C_{0, \beta\gamma, 1}$ and $C_{0, \beta\gamma, 2}$ to them.
We write these as $Y_1$ and $Y_2$.
Then let $\mu(i)_1$ on $Y_1$ be the image of the section of 
 $\varphi_0^*\Theta_{X_0}$ corresponding to the vector $\mathfrak w_i$
 in $\mathcal N_{C_0/X_0}$.
Take all the other $\mu(i)_l$ to be zero.

Then clearly 
\[
(\mathfrak w_j^{\vee}\otimes \xi, \{\mu(i)_l\}) = \delta_{ij}
\]
 holds.
This proves that the pairing is non-degenerate.\qed\\

Thus, a \v{C}ech 1-cocycle $\{\mu_l\}$ is zero if and only if its pairing with 
 any of elements of $H$ is zero.

Now let us apply this pairing to our calculation of the obstruction.
Assume that we have a lift $\varphi_k$ of $\varphi_0$.
By Lemmas \ref{lem:canlift} and \ref{lem:liftT}, 
 the map $\varphi_k$ restricted to the subsets 
 $\{\tilde Y_l\} = \{C_{0, T_i}^*\}\cup \{C_{0, \beta\gamma}\}$
 (or their suitable lifts)
 has lifts.
Choosing any one of lifts for each of the open covers $\{\tilde Y_l\}$, 
 the differences of them on the intersections of these covers 
 define the obstruction class.

To calculate the 
 pairing between this class and a class in $H$,
 we need to present the obstruction class by sections $\{\mu_l\}$ on 
 locally closed subsets $\{Y_l\}$. 
Such sections can be obtained from local lifts as follows.
Recall that on an open subset of the form $C_{0, \beta\gamma}$,
 there is a canonical 
 one to one correspondence between the set of $k$-th order lifts of 
 $\varphi_0$ restricted to it and the suitable set of sections of the normal sheaf
 $\mathcal N_{C_0/\mathfrak X}$ with its coefficients extended
 (see Lemma \ref{lem:lifthigh} and the paragraphs before it).

In particular, on a set of the form $C_{0, \beta\gamma}$, given a $k$-th order lift of
 $\varphi_0|_{C_{0, \beta\gamma}}$, the set of $(k+1)$-th order lifts which reduce
 to the given $k$-th order lift canonically corresponds to the set of pairs of sections
 $(\mathfrak m_{k, \beta}, \mathfrak m_{k, \gamma})$, where
 $\mathfrak m_{k, \beta}\in t^{k+1}\Bbb C[t]/t^{k+2}\otimes_{\Bbb C}
  \Gamma(C_{0, \beta\gamma, 1}, \mathcal N_{C_0/X_0})$ and
  $\mathfrak m_{k, \gamma}\in t^{k+1}\Bbb C[t]/t^{k+2}\otimes_{\Bbb C}
    \Gamma(C_{0, \beta\gamma, 2}, \mathcal N_{C_0/X_0})$
  such that they coincide at the node of $C_{0, \beta\gamma}$.

For simplicity, we assume that a neighborhood of the edge connecting $\beta$ and $\gamma$
 is standard in the sense of Definition \ref{def:twovertex} 
 and the length of it is one (otherwise, taking a suitable 
 covering as in Subsection \ref{subsec:general_edge}, we can reduce the situation  
 to the standard one.
See also Remark \ref{rem:coord} (1) and Remark \ref{rem:secorder2} (3)).
Assume that the given $k$-th order lift of $\varphi_0|_{C_{0, \beta\gamma}}$
 corresponds to a section $(\mathfrak n_{k-1, \beta}, \mathfrak n_{k-1, \gamma})$
 in the expression of Lemma \ref{lem:lifthigh}.
We write its components explicitly as follows:
\[\begin{array}{ll}
\mathfrak n_{k-1, \beta}
  & = (y\partial_y-\frac{\kappa y}{\lambda w_1}\cdot w_1\partial_{w_1})
      + c(s)(x\partial_x-y\partial_y) + c_0(s)z\partial_z + 
        c_1(s)w_1\partial_{w_1}+\cdots +c_{n-2}(s)w_{n-2}\partial_{w_{n-2}}
      \\
     & \in 
 (y\partial_y-\frac{\kappa y}{\lambda w_1}\cdot w_1\partial_{w_1})
    + t\Bbb C[t]/t^{k+1}\otimes_{\Bbb C}
  \Gamma(C_{0, \beta\gamma, 1}, \mathcal N_{C_0/X_0}),
 \end{array}
 \]
 \[\begin{array}{ll}
 \mathfrak n_{k-1, \gamma}
   & = (x\partial_x-\frac{k x}{l z}\cdot z\partial_{z})+
   c'(u)(y\partial_y-x\partial_x) + c'_0(u)z\partial_z + 
     c'_1(u)w_1\partial_{w_1}+\cdots +c'_{n-2}(u)w_{n-2}\partial_{w_{n-2}}\\
    &\in 
   (x\partial_x-\frac{k x}{l z}\cdot z\partial_{z})
     + t\Bbb C[t]/t^{k+1}\otimes_{\Bbb C}
    \Gamma(C_{0, \beta\gamma, 2}, \mathcal N_{C_0/X_0}).
  \end{array}
  \] 

From these expressions, we can obtain sections $\mu_l$ for
 the locally closed subsets $C_{0, \beta\gamma, 1}$
 and $C_{0, \beta\gamma, 2}$.
Namely, first consider the section $\mathfrak n_{k-1, \beta\gamma}$ of 
 $\Bbb C[t]/t^{k+1}\otimes \mathcal N_{C_0/X_0}|_{C_{0, \beta\gamma}}$
 given by
\[\begin{array}{l}
-\frac{\kappa u}{\lambda w_1}\cdot w_1\partial_{w_1}
  - \frac{ks}{lz}\cdot z\partial_z + (c(s)-c'(u))(x\partial_x-y\partial_y)\\
 \hspace{.2in} + (c_0(s)+c'_0(u)-\gamma_0)z\partial_z + (c_1(s)+c'_1(u)-\gamma_1)w_1\partial_{w_1}
  +\cdots + (c_{n-2}(s)+c_{n-2}(u)-\gamma_{n-2})w_{n-2}\partial_{w_{n-2}},
\end{array}
\]
 here $\gamma_i$ is the constant term of $c_i(s)$
  (and also of $c_i'(u)$, see Lemma \ref{lem:lifthigh}).
Also, the functions $z$ and $w_1$ in the coefficients $-\frac{\kappa u}{\lambda w_1}$ and
 $- \frac{ks}{lz}$ of the first two terms in fact mean their
 pull back to $C_{0, \beta\gamma}$
 by the map $\varphi_0$, so that they are actually functions in $s$ and $u$
 (with nonzero constant terms), respectively.
\begin{rem}
Note that this does not depend on the choice of functions
 satisfying the conditions of Lemma \ref{lem:coordinate}.
Namely, other such set of functions has the form
\[
(X_1, X_2, X_3, \dots, X_{n+1}) = 
 (xs_1, ys_2, s_3, 
  \dots, s_{n+1}),
\]
 where $s_1, \dots, s_{n+1}$ are monomials of $z, w_1, \dots, w_{n-2}$.
In this case, we have $X_1\partial_{X_1} = x\partial_x$ and $X_2\partial_{X_2} = y\partial_y$, 
 and this guarantees the independence of the above section with respect to 
 the choice of functions.
\end{rem}

Now as in Remark \ref{rem:param}, the parameters $s$ and $u$ can be 
 regarded as parameters of the deformed curve $C_{k+1, \beta\gamma}$
 so that they are related by $su = t$.
In other words, we can regard the section
 $\mathfrak n_{k-1, \beta\gamma}$ above as a section of 
 the log normal sheaf $\mathcal N_{C_{k+1, \beta\gamma}/\mathfrak X}$
 which is evaluated to zero by the covector $\frac{dt}{t}$.
Here $C_{k+1, \beta\gamma}$ is mapped to $\mathfrak X$ by the 
 trivial lift of $\varphi_0|_{C_{0, \beta\gamma}}$ (see Lemma \ref{def:trivlift2}).

Now write $\mathfrak n_{k-1, \beta\gamma}$ in terms of $s$ using the relation
 $u = \frac{t}{s}$ and take the parts of order $k+1$ with respect to 
 the power of $t$.
Then restricting it to the component $C_{0, \beta}^*$, we obtain a section
 of the sheaf $\mathcal N_{C_0/X_0}(\infty p_{\beta\gamma})|_{C_{0, \beta\gamma, 1}}$.
Similarly, by writing $\mathfrak n_{k-1, \beta\gamma}$ in terms of $u$,
 taking the parts of order $k+1$ with respect to 
 the power of $t$ and then restricting it to $C_{0, \gamma}^*$, we obtain a section
 of the sheaf $\mathcal N_{C_0/X_0}(\infty p_{\beta\gamma})|_{C_{0, \beta\gamma, 2}}$.
These give the sections $\mu_l$ for the locally closed subsets 
 $C_{0, \beta\gamma, i}$.
By similar calculation, we can attach a section $\mu_l$ to 
 each locally closed subset $Y_l$.
For subsets of the form $C_{0, T_i}^*$, we can attach any section of 
 $\mathcal N_{C_0/X_0}|_{C_{0, T_i}^*}$
 (since their support do not intersect the support of representatives of the space $H$,
 any choice gives the same value of the pairing).

Also, one can see that the residues of the resulting sections do not depend
 on the choice of a lift of $\varphi_{k}$ on $C_{k, \beta\gamma}$
 (in the above paragraph, we take the trivial lift introduced in Lemma \ref{lem:canlift}),
 since other lifts differ from it only by sections of 
 $t^k\Bbb C[t]/t^{k+1}\otimes
 \Gamma(C_{0, \beta\gamma}, \mathcal N_{C_0/X_0})$, and they do not contribute to
 the residues at $p_{\beta\gamma}$.

We can show the following by a straightforward inspection of the definition.
\begin{prop}\label{prop:obstdef}
The \v{C}ech cohomology class defined by $\{\mu_l\}$ coincides with 
 the obstruction class given in Definition \ref{def:cechobst}. \qed
\end{prop}

\subsubsection{Calculation of the pairing}\label{subsec:pairing}

By Lemma \ref{lem:nondegeneracy} and Proposition \ref{prop:obstdef}, 
 the map $\varphi_k$ can be lifted to the next order if and only if 
 the pairings between $\{\mu_l\}$ and all the classes in $H$ vanish.
Therefore, let us calculate the pairing.

Now let $\alpha_N$ be a vertex on the loop and $\alpha_{N-1}$ be its adjacent
 vertex which is not on the loop, see Figure \ref{fig:residue}.
For simplicity, we assume that $\alpha_{N-1}$ is not a 1-valent vertex of $\overline\Gamma'$
 (we can always assume this by applying a base change and adding 2-valent vertices).
Consider the image of the open subset $C_{k, \alpha_N\alpha_{N-1}}$
 by $\varphi_k$.
 
We fix a basis $\{\mathfrak v_1, \dots, \mathfrak v_{n-r}, 
 \mathfrak w_1, \dots, \mathfrak w_{r}\}$ of $N$
 and its dual basis
 $\{\mathfrak v_1^{\vee}, \dots, \mathfrak v_{n-r}^{\vee}, 
 \mathfrak w_1^{\vee}, \dots, \mathfrak w_{r}^{\vee}\}$
 as before. 
We take a set of functions
 $\{x^{(N-1)}_1, y_2^{(N-1)}, x_2^{(N-1)}, y_1^{(N-1)}, z_{1}^{(N-1)}, \dots, z_{n-3}^{(N-1)}\}$
  on $\mathfrak X$ 
 (or on its suitable covering $\mathfrak X_s$ as in Subsection \ref{subsec:general_edge},
 when the neighborhood of the edge connecting $\alpha_N$ and $\alpha_{N-1}$
 is not standard in the sense of Definition \ref{def:twovertex})
 as before
 around the image of the node of $C_{0, \alpha_N\alpha_{N-1}}$.
We can assume that the functions $z_{n-r-2}^{(N-1)}, \dots, z_{n-3}^{(N-1)}$
 correspond to the vectors $\mathfrak w_1^{\vee}, \dots, \mathfrak w_r^{\vee}$.
That is, the covectors $\frac{dz_{n-r-3+i}^{(N-1)}}{z_{n-r-3+i}^{(N-1)}}$, seen as an element of 
 $N^{\vee}$ by the isomorphism $\Omega_{\mathfrak X}\cong \mathcal O_{\mathfrak X}
 \otimes (N\oplus \Bbb Z)^{\vee}\cong
  \mathcal O_{\mathfrak X}\oplus(N^{\vee}\oplus\Bbb Z^{\vee})$,
 equal to $\mathfrak w_i^{\vee}$, $i = 1, \dots, r$.

We calculate the pairing between the obstruction class associated to $\varphi_k$ 
 and the class $\eta_i = \mathfrak w_i^{\vee}\otimes\xi$
 in the dual obstruction space $H$,
  at the node $p_{\alpha_N\alpha_{N-1}}$ of $C_{0, \alpha_N\alpha_{N-1}}$.
In this case, we only need to take account of the 
 $z_{n-r-3+i}^{(N-1)}\partial_{z_{n-r-3+i}^{(N-1)}}$-component $\mu_{l, i}$ of the 
 obstruction cohomology class $\{\mu_l\}$.

The restriction of $\varphi_k$ to $C_{k, \alpha_N\alpha_{N-1}}$
 corresponds to a section $(\mathfrak n_{k-1, \alpha_N}, \mathfrak n_{k-1, \alpha_{N-1}})$
 as argued above.
First note that $\eta_i$ is identically zero on the component $C_{0, \alpha_N\alpha_{N-1}, 2}$
 (that is, the union of $C_{0, \alpha_{N-1}}^*$ and the node $p_{\alpha_N\alpha_{N-1}}$).
Therefore, there is a contribution to the pairing only from the component 
 $C_{0, \alpha_N\alpha_{N-1}, 1}$.
Note also that the class $\eta_i$ does not have a zero or a pole at the node
 $p_{\alpha_N\alpha_{N-1}}$.
Thus, from the definition of the pairing $(\heartsuit)$ above, we only need to 
 calculate the first order pole of the section $\mu_{l, i}$ on 
 $C_{0, \alpha_N\alpha_{N-1}, 1}$ at $p_{\alpha_N\alpha_{N-1}}$.

According to the explicit correspondence between the set of sections of the normal sheaf
 and
 the set of local lifts given in Lemma \ref{lem:lifthigh} and Proposition \ref{prop:obstdef}, 
 we can check the following.

\begin{prop}\label{prop:anal-pair}
Let $\{\mu_l\}$ be the obstruction cohomology class associated to $\varphi_k$ as above.
Its pairing with $\eta_i$ at the point $p_{\alpha_N\alpha_{N-1}}$ of $C_{0, \alpha_N\alpha_{N-1}, 1}$
 coincides with the $t^{k+1}$-part of the
 residue of $\varphi_{k+1}^*z_{n-r-3+i}^{(N-1)}
 \cdot\frac{dU_N}{U_N}$ at $p_{\alpha_N\alpha_{N-1}}$ (which is the obstruction at $\alpha_N$
 of the direction $\mathfrak w_i^{\vee}$ in the sense of Definition \ref{def:preob}),
 here $\varphi_{k+1}$ is any local lift of $\varphi_k$ on $C_{k, \alpha_N\alpha_{N-1}}$.
Also, the function $z_{n-r-3+i}^{(N-1)}$ is normalized according to 
 Remark \ref{rem:normalization}.\qed
\end{prop}

We have calculated the value of $\varphi_{k+1}^*z_{n-r-3+i}^{(N-1)}$
 on $C_{k, \alpha_N\alpha_{N-1}}$ by the analytic continuation. 
In particular, the contribution to the pairing from the pole $p_{\alpha_N\alpha_{N-1}}$
 is the same as the sum of the residues of the pre-obstructions $(\spadesuit)$
 in Proposition \ref{prop:obstruction1} contributed from 1-valent vertices $\{p_i\}$
 of $\overline\Gamma'$
 such that $\alpha_N>p_i$ in the notation introduced before Lemma \ref{lem:prinpart}.

These are the only contribution to the pairing, as the following claim shows.

\begin{lem}
Let $\beta$ and $\gamma$ be neighboring vertices on the loop part $L$ of
 $\Gamma$.
The contribution to the
 pairing between the obstruction class $\{\mu_l\}$ and the dual class $\eta_i$
 at the node $p_{\beta\gamma}$ in $C_{0, \beta\gamma}$
 is zero.
\end{lem}
\proof
First note that at the node $p_{\beta\gamma}$, the differential form $\eta_i$
 has a pole of order one.
Therefore, the constant part of $\{\mu_l\}$ contributes to the pairing there.

Now by definition of the sections $\mu_l$ in Subsection \ref{subsec:resobs},
 the constant terms of $\mu_l$ on the locally closed subsets
 $C_{0, \beta\gamma,1}$ and
 $C_{0, \beta\gamma,2}$
 coincide.
On the other hand, when $s$ and $u$ are local coordinates of 
 the branches $C_{0, \beta\gamma,1}$ 
 and $C_{0, \beta\gamma,2}$ at the node
 $p_{\beta\gamma}$, the relation 
 $\frac{ds}{s}+\frac{du}{u} = 0$ holds.
Therefore, the contributions to the pairing from these branches
 cancel, proving the lemma.\qed.\\

Therefore, we have the following.
\begin{cor}\label{cor:anal-ob}
The pairing $(\eta_i, \{\mu_l\})$ equals the obstruction $o(\mathfrak w_i^{\vee}, \varphi_k)$ 
 of the direction 
 $\mathfrak w_i^{\vee}$ in the sense of Definition \ref{def:preob}.\qed
\end{cor}

This explains the reason why we call $o(\mathfrak w_i^{\vee}, \varphi_k)$ the obstruction.
Namely, By Proposition \ref{prop:obstdef},
 the obstruction to lift $\varphi_k$ to one step further vanishes if and only if
 $o(\mathfrak w_i^{\vee}, \varphi_k)$ vanishes for all directions.

\begin{rem}
From this corollary, one can see that the vanishing of the 
 obstruction $o(\mathfrak w_i^{\vee}, \varphi_k)$ is equivalent to the 
 condition that there is a rational 1-form on $C_{0, L}$ with 
 prescribed poles at the points as $p_{\alpha_N\alpha_{N-1}}$
 whose residue is given by the value of the pairing between
 $\{\mu_l\}$ and $\eta_i$ there.
Since there is a natural canonical 1-form on $C_{0, L}$, 
 this is also the same as the existence of a rational function of the
 prescribed pole.

In the cases of curves of higher genus too, the vanishing of the obstruction
 is reduced to the same kind of Cousin type problems.
\end{rem}
\subsection{Calculation of obstructions in embedded cases.  Step 6: Cancelation of the obstruction and existence of smoothings.}\label{subsec:step6}
In the previous subsection, we identified the obstruction $o(\mathfrak w^{\vee}, \varphi_k)$
 to lift $\varphi_k$ to one step further with the sum of residues of suitable functions.
In this subsection, we look closely at this sum and deduce the necessary and sufficient 
 condition under which the obstruction vanishes.

\begin{lem}\label{lem:linearity}
The function $o(\mathfrak w^{\vee};\varphi_{k})$ 
 is linear with respect to the variable $\mathfrak w^{\vee}$.
\end{lem}
\proof
This follows from Corollary \ref{cor:anal-ob}
 since the pairing which calculates the obstruction 
 is linear with respect to elements of the dual obstruction space $H$.\qed

\begin{rem}\label{rem:kuranishi}
\begin{enumerate}
\item We can extend the domain $(\bar A)^{\perp}\cap N^{\vee}$ of the map
 $o(\cdot;\varphi_k)$ to $(\bar A_{\Bbb C})^{\perp}$ by linearity.
Then recalling that the dual of the obstruction class
 $H = (H^1(C_0, \mathcal N_{C_0/X_0}))^{\vee}$ is isomorphic to
 $(\bar A_{\Bbb C})^{\perp}$, 
 the map $o(\cdot;\cdot)\colon (\bar A_{\Bbb C})^{\perp}\times 
 H^0(C_0, \mathcal N_{C_0/X_0})\to
 t^{k+1}\Bbb C[t]/t^{k+2}$ is the Kuranishi map of order $k$ in
 Subsection \ref{subsec:Kuranishi}.
Here the space $H^0(C_0, \mathcal N_{C_0/X_0})$
 is identified with the space of lifts $\varphi_{k}$ with a given
 $\varphi_{k-1}$.
\item The map $o(\cdot;\cdot)\colon (\bar A_{\Bbb C})^{\perp}\times 
 H^0(C_0, \mathcal N_{C_0/X_0})\to
 t^{k+1}\Bbb C[t]/t^{k+2}$ is affine linear with respect to  
 the part $H^0(C_0, \mathcal N_{C_0/X_0})$.
\end{enumerate}
\end{rem}
Then we have the following from the definition of the obstruction.
\begin{prop}\label{prop:liftability}
The map $\varphi_{k}$ is liftable to a map of order $k+1$
 if and only if 
 $o(\mathfrak w^{\vee};\varphi_{k})$  vanishes for any 
 $\mathfrak w^{\vee}$.\qed
\end{prop}
\begin{cor}\label{cor:basis}
The map $\varphi_{k}$ is liftable to a map of order $k+1$
 if and only if $o(\mathfrak w_i^{\vee};\varphi_{k})$ vanish for a basis
 $\mathfrak w_1^{\vee}, \cdots \mathfrak w_{r}^{\vee}$ of 
 $N^{\vee}\cap (\bar A)^{\perp}$.
\end{cor}
\proof
This follows from the linearity of 
 $o(\mathfrak w^{\vee};\varphi_{k})$
 with respect to the variable 
  $\mathfrak w^{\vee}$ (Lemma \ref{lem:linearity}).\qed

\begin{rem}\label{rem:obstbasis}
In fact, the set of vectors $\mathfrak w_1^{\vee}, \cdots \mathfrak w_{r}^{\vee}$
 need not be a basis of $N^{\vee}\cap (\bar A)^{\perp}$.
It should only be a basis of $(\bar A)^{\perp}$ over $\Bbb R$ by linearity.
However, we will at least take $\mathfrak w_i^{\vee}$ in $N^{\vee}$
 since otherwise the correspondence between the obstruction and
 the residue of functions becomes unclear.
\end{rem}

\begin{rem}\label{rem:basis-ob}
It is clear from the corollary that the vanishing or non-vanishing of the obstruction 
 does not depend on the choice of a basis
 $\{\mathfrak v_1, \dots, \mathfrak v_{n-r}, \mathfrak w_1, \dots, \mathfrak w_r\}$.
See Remark \ref{rem:direction}.
\end{rem}

The precise calculation of 
 $o(\mathfrak w^{\vee}; \varphi_k)$ will be hard.
Fortunately, if we know its leading term,
 we will have a good control of the entire obstruction (see the proof of 
 Theorem \ref{thm:immersive} below),  
 and the calculation of the leading term has essentially already been  
 done.

Let $\alpha$ be a 1-valent vertex 
 of $\overline\Gamma'$.
We choose $\alpha$ so that
 the edge length of the path $\mathcal P$
\[
\alpha\to \alpha_1\to \cdots \to \alpha_N\in L
\]
 from $\alpha$ to the loop is one of the shortest among the 1-valent
 vertices of $\overline\Gamma'$.
Let $A_1$ be the minimal affine subspace of $N_{\Bbb R}$ containing
 $\overline\Gamma'$ and the edges emanating from $\alpha$.
Obviously, 
\[
\dim A_1 = \dim A+1
\]
 holds.
Let $\bar A_1$ be the linear subspace of $N_{\Bbb R}$ parallel to $A_1$.
The following is clear from the calculation in Subsection \ref{subsec:step5}.
\begin{lem}\label{lem:transversality}
Let $\{\mathfrak w_1^{\vee}, \dots, \mathfrak w_r^{\vee}\}$ be a basis of 
 $(\bar A)^{\perp}\cap N^{\vee}$ such that 
 $\{\mathfrak w_2^{\vee}, \dots, \mathfrak w_r^{\vee}\}$ is a basis of 
 $(\bar A_1)^{\perp}\cap N^{\vee}$.
Then the pre-obstruction at $\alpha_N$ of the direction $\mathfrak w_1^{\vee}$
 contributed from $\alpha$ 
 (see Definition \ref{def:preobcontri})
 is 
 dominated by the data of the vertices $\alpha, \alpha_1, \dots, \alpha_N$
 in the sense that the leading term of 
 it with respect to 
 the exponent of $t$ 
 only depends on the 
 parameters of the defining equations of $\varphi_0(C_{0, \alpha})$
 $\varphi_0(C_{0, \alpha_1}), \dots, \varphi_0(C_{0, \alpha_N})$.
See also Remark \ref{rem:dependence}.
\qed
\end{lem}

Now let $\overline\Gamma^{(2)}$ be the subgraph of $\Gamma$
 which is the unique connected component of 
 $\Gamma\cap A_1$
 containing the loop $L$.
The graph $\overline\Gamma^{(2)}$ also
 has several 1-valent vertices.
Let $\alpha^{(2)}$ be one of these 1-valent vertices
 with the property that 
 the path length from $\alpha^{(2)}$ to the loop is the shortest 
 among them.
Let $\alpha^{(2)}_M$ be the vertex on the loop closest to 
 $\alpha^{(2)}$.

Then do the same calculation as in Subsection \ref{subsec:step5}
 to this vertex $\alpha^{(2)}$.
This calculates an obstruction whose direction is
 contained in $(\bar A_1)^{\perp}\cap N^{\vee}$.
Namely, take $\mathfrak w_3^{\vee}, \dots, \mathfrak w_r^{\vee}$
 in the statement of Lemma \ref{lem:transversality}
 so that they consist a basis of $(\bar A_2)^{\perp}\cap N^{\vee}$.
Then the leading term of the pre-obstruction at $\alpha_M^{(2)}$ of the direction
 $\mathfrak w_2^{\vee}$ contributed from $\alpha^{(2)}$
 is dominated by the data of the vertices on the unique shortest
 path from $\alpha^{(2)}$
 to $\alpha^{(2)}_M$ as in Lemma \ref{lem:transversality}.

Next, take the minimal affine subspace $A_2$
 containing $A_1$ and the edges 
 emanating from $\alpha^{(2)}$, and continue the same process.

Through this process, the flag of affine subspaces
\[
A = A_0\subset A_1\subset A_2\subset\cdots
 \subset A_{r-1}\subset A_r = N_{\Bbb R}
\]
 is defined (if there are several 1-valent
 vertices with the property above,
 the flag of subspaces may not be uniquely determined.
 But it does not matter to the 
 result below).
Note that the codimension of the subspace $A_i$ is $r-i$.
Let $\bar A_i$ be the linear subspace of $N_{\Bbb R}$
 parallel to $A_i$.

Accordingly, we can choose a sequence of vectors 
\[
\mathfrak w_1^{\vee}, \mathfrak w_2^{\vee}, \dots, \mathfrak w_r^{\vee}
\]
 in $(\bar A)^{\perp}\cap N^{\vee}$
 with the following properties:
\begin{itemize}
\item The vector 
 $\mathfrak w_i^{\vee}$ belongs to $(\bar A_{i-1})^{\perp}\setminus
 (\bar A_i)^{\perp}$.
\item The set $\{\mathfrak w_1^{\vee}, 
 \mathfrak w_2^{\vee}, \dots, \mathfrak w_r^{\vee}\}$
 composes a basis of
 $(\bar A)^{\perp}\cap N^{\vee}$
\end{itemize}
Also,
 we have a series of total obstructions
 $o(\mathfrak w_1^{\vee};\varphi_k), o(\mathfrak w_2^{\vee};\varphi_k),
  \cdots, o(\mathfrak w_r^{\vee};\varphi_k)$
  associated to these vectors.

As we stated in Proposition \ref{prop:liftability},
 the necessary and sufficient condition for the existence of a lift of
 $\varphi_{k}$ to the $(k+1)$-th order
 is the vanishing of these obstructions.
Note that each element of $(\bar A)^{\perp}$ gives a hyperplane
 in $N_{\Bbb R}$, that is, the set of vectors annihilated by that element.
 
When we assume that $(\Gamma, h)$ is an immersion
 (we have been assuming that $(\Gamma, h)$ is an embedding
 in this section, but relaxing the assumption
 to the immersive case requires essentially no change)
 and all the edge weights are one,
 the desired condition precisely 
 coincides with the \emph{well-spacedness condition} considered by
 Speyer \cite{S}.
When the edge weights are general, it is modified using the path length of 
 Definition \ref{def:path length}.
\begin{defn}\label{well-spaced}
An immersive superabundant tropical curve
 $(\Gamma, h)$ of genus one is said to be \emph{well-spaced}
 with respect to an affine hyperplane 
 $\mathcal H$
 of $N_{\Bbb R}$ containing $h(\overline\Gamma')$
  if the following condition is satisfied. 
Let 
\[
\Gamma_{\mathcal H} \subset h(\Gamma)\cap \mathcal H
\]
 be the connected component containing $h(\overline\Gamma')$ and
 let 
\[
p_1^{\mathcal H}, \dots, p_j^{\mathcal H}
\]
 be the 1-valent vertices of it.
Denote by 
\[
\mathcal P^{\mathcal H}_i,\;\; i = 1, \dots, j
\]
 the unique path connecting
 $p_i^{\mathcal H}$ to the loop $h(L)$.
Then the set
\[
\{\ell_{(\Gamma, h)}(\mathcal P^{\mathcal H}_1),
  \dots, \ell_{(\Gamma, h)}(\mathcal P^{\mathcal H}_j)\}
\]
  of positive integers contains at least two minimum.

We call $(\Gamma, h)$ well-spaced if it is well-spaced with respect to all 
 hyperplanes containing $h(\overline{\Gamma}')$.
\end{defn}
\begin{thm}\label{thm:immersive}
Assume that $(\Gamma, h)$ is an immersive
 superabundant tropical curve of genus one.
Then it is smoothable (see Definition \ref{def:smoothable})
 if and only if 
 $(\Gamma, h)$ is well-spaced.
\end{thm}
\proof
Let $\varphi_0\colon C_0\to \mathfrak X$ be a pre-log curve of 
 type $(\Gamma, h)$.
(it always exists: If $A$ is the minimal dimensional affine subspace
containing the loop part of $h(\Gamma)$, then $h(\Gamma')$ can be seen as
a tropical curve in $\bar A$, which is regular in the sense of 
Definition \ref{def:nonsuperabundant2}.
It is not difficult to see that for regular tropical curves, there always exist
corresponding pre-log curves (in fact, it is a part of the claim of 
Theorem \ref{thm:regsm}).
Since the complement of $h(\Gamma')$ in $h(\Gamma)$ is a tree, 
it is easy to extend such a pre-log curve corresponding to 
$h(\Gamma')$ to a pre-log curve corresponding to $h(\Gamma)$.).
  
Assume that we have its $k$-th order lift $\varphi_k$.
We will show that we can perturb the map $\varphi_k$ to $\varphi_k'$
 so that the obstruction $o(\mathfrak w^{\vee};\varphi_k')$ vanishes for any 
 $\mathfrak w^{\vee}\in N^{\vee}\cap (\bar A)^{\perp}$
 if and only if $(\Gamma, h)$ is well-spaced.
\begin{rem}\label{rem:lowob}
To be more precise, for small $k$, the obstruction to lift $\varphi_k$
 does not exist if the distance from the 1-valent vertices of 
 $\overline{\Gamma}'$ to the loop is large.
Therefore, the 'only if' part of the claim does not hold for such $k$.
These cases are not essential to the calculation of the obstruction, 
 and we omit cumbersome case statement associated to them.
\end{rem}

As usual,
 we assume that the direction vectors of the edges of $h(\Gamma)$ span
 $N_{\Bbb R}$.
The pre-obstruction at a vertex on the loop 
 defined in Definition \ref{def:preob}
 can be written in the form
\[
\zeta_0(\mathfrak w^{\vee};\varphi_0) + 
 \zeta_1(\mathfrak w^{\vee}),
\]
 here 
 $\zeta_0(\mathfrak w^{\vee};\varphi_0)$ is the leading term,
 which depends only on the data of 
 $\varphi_0(C_0)$ according to
 Proposition \ref{prop:obstruction1}. 
The part $\zeta_1(\mathfrak w^{\vee})$ is the sum of the higher order 
 terms.

First let us assume that the dual space 
 of obstructions $H=(H^1(C_0, \mathcal N_{C_0/\mathfrak X}))^{\vee}$
 has dimension one
 for simplicity.
In this case, 
 we can think of the pre-obstruction as a polynomial in the variable $t$.  
 
Let $a$ be the order of the leading term
 $\zeta_0(\mathfrak w^{\vee};\varphi_0)$ with respect to $t$. 
Suppose we perturb $\varphi_0$ by adding terms of order $t^b$
 to the coefficients of the defining equations, here $b$ is a positive integer.

Explicitly, 
 let $\alpha$ be a 1-valent vertex of $\overline\Gamma'$
 which is one of the closest to the loop among such 1-valent vertices.
Then 
 using the notation in Subsection \ref{subsec:step4}, 
 we have an equation 
\[
k_0y_1^{(0)}+l_0x_1^{(0)}+m_0 = 0,
\]
 which is a part of the defining equations of 
 the image of the map $\psi_{0, 0}$,  
 which is mapped to the image $\varphi_0(C_{0, \alpha}\cup C_{0, \alpha_1})$
 by a toric morphism (see the paragraph after Remark \ref{rem:length}).

Let $\mathcal U_1, \mathcal U_2$ be the components of 
 $\Gamma\setminus\overline\Gamma'$
 which contain the vertex $\alpha$ in their closures, and let
 $C_{0, \mathcal U_1}, C_{0, \mathcal U_2}$ be the union of components
 of $C_0$ corresponding to them.
When $\mathcal U_1$ or $\mathcal U_2$ is an unbounded edge, 
 then take $C_{0, \mathcal U_1}$ or $C_{0, \mathcal U_2}$
 to be the empty set.

We modify the above equation to
\[
k_0y_1^{(0)}+(l_0+l_0')x_1^{(0)}+m_0 = 0
\]
 by a non-zero complex number $l_0'$.
 
Then we can construct a pre-log curve of type $(\Gamma, h)$
 from $\varphi_0(C_0)$
 without modification except the part 
 $\varphi_0(C_{0, \alpha})\cup \varphi_0(C_{0, \mathcal U_1})
 \cup\varphi_0(C_{0, \mathcal U_2})$.
The defining equations of the parts
 $\varphi_0(C_{0, \mathcal U_1})$ and $\varphi_0(C_{0, \mathcal U_2})$
 are also
 modified, while those of the other components of $C_0$ are not modified. 
Then the coefficient of the leading term of the pre-obstruction
\[
t^a\frac{l_0}{m_0}\cdot\frac{k_1}{m_1}
 \cdot \frac{l_1}{m_1}\cdot\cdots\cdot
 \frac{k_{N-1}}{m_{N-1}}\cdot \frac{l_{N-1}}{m_{N-1}}\cdot 
 \frac{k_N}{m_N}
\]
 contributed from the vertex $\alpha$ is modified to 
\[
t^a\frac{l_0+l_0'}{m_0}\cdot\frac{k_1}{m_1}
 \cdot \frac{l_1}{m_1}\cdot\cdots\cdot
 \frac{k_{N-1}}{m_{N-1}}\cdot \frac{l_{N-1}}{m_{N-1}}\cdot 
 \frac{k_N}{m_N}.
\]
While we took $l_0'$ to be a complex number, we can also take it to be 
 $\varepsilon_0t^b$ for some complex number $\varepsilon_0$ 
 and a positive integer $b$.
In this case, the coefficient of the leading term is modified to 
\[
t^a\frac{l_0}{m_0}\cdot\frac{k_1}{m_1}
 \cdot \frac{l_1}{m_1}\cdot\cdots\cdot
 \frac{k_{N-1}}{m_{N-1}}\cdot \frac{l_{N-1}}{m_{N-1}}\cdot 
 \frac{k_N}{m_N}.
  + t^{a+b}\frac{\varepsilon_0}{m_0}\cdot\frac{k_1}{m_1}
 \cdot \frac{l_1}{m_1}\cdot\cdots\cdot
 \frac{k_{N-1}}{m_{N-1}}\cdot \frac{l_{N-1}}{m_{N-1}}\cdot 
 \frac{k_N}{m_N}.
\]
Note that the leading term of 
 the pre-obstruction does not change through this perturbation.
Moreover, the part 
 $\zeta_1(\mathfrak w^{\vee})$ is modified at most
 in order $t^{a+b+1}$.

The obstruction to lift $\varphi_k$
 is calculated from the pre-obstruction by taking the sum of 
 appropriate residues of order $k+1$ contributed from 1-valent vertices
 of $\overline\Gamma'$, see Definitions \ref{def:preob} and 
 \ref{def:preobcontri}.
Recall that we write the order of the leading term of the pre-obstruction
 by $a$.
By the calculation above, if $k\geq a$, then 
 the obstruction $o(\varphi_k)$, 
 which is of order $k+1>a$, can be cancelled  
 by modifying 
 the coefficients of the defining equations of $\varphi_0$ by terms of
 order $k+1-a$ with respect to $t$.
 
Therefore, for the vanishing of the obstruction at any order,
 it is enough to show that we can choose $\varphi_0$
 so that the obstruction $o(\varphi_{a-1})$ vanishes
 (note that for $c<a-1$, any deformation $\varphi_c$ of
 $\varphi_0$ does not have an obstruction to lift, see Remark \ref{rem:lowob}).
By definition, $o(\varphi_{a-1})$ is given by the sum of the 
 residues of the leading terms of the pre-obstructions contributed from 
 1-valent vertices of $\overline\Gamma'$.
 
By Proposition \ref{prop:obstruction1}, it is easy to see
 that the necessary and sufficient 
 condition for this is the well-spacedness of the
 tropical curve $(\Gamma, h)$.
Namely, in this case we can change the coefficient $l_0$ of
 the defining equation to another suitable constant $\tilde l_0$
 to cancel the residues contributed from other 1-valent vertices of $\overline\Gamma'$.

If $\dim H$ is larger than one, 
 we apply the above argument inductively to the 
 vertices of $\Gamma\setminus \Gamma'$, 
 starting from the ones closer to the loop, 
 as in the argument after Lemma \ref{lem:transversality}.
Using the same notation as there, we have
 a basis $\{\mathfrak w_1^{\vee}, \cdots,\mathfrak w_r^{\vee}\}$
 of $(\bar A)^{\perp}\cap N^{\vee}$.
Then by Corollary \ref{cor:basis}
 the necessary and sufficient condition 
 for the map $\varphi_k$ to have a lift $\varphi_{k+1}$ of 
 order $k+1$ is the vanishing of the obstructions 
 $o(\mathfrak w_i^{\vee};\varphi_k)$ for all $i = 1, \dots, r$.
It suffices to prove that
 the existence of a map $\varphi_k$ satisfying 
 this condition is equivalent to the well-spacedness
 condition for $(\Gamma, h)$.
 
If $(\Gamma, h)$ is well-spaced, then by
 definition the well-spacedness conditions hold for  
 all the hyperplanes corresponding to the vectors $\mathfrak w_i^{\vee}$. 
Then the above argument for the case $\dim H = 1$ extends straightforwardly to
 show that we can cancel these obstructions by perturbing $\varphi_k$.
 
Conversely, it is clear that for the vanishing of these obstructions, 
 the well-spacedness conditions must be satisfied by 
 the hyperplanes corresponding to the vectors 
 $\mathfrak w^{\vee}_i$.
Moreover, if the well-spacedness condition is satisfied by
 any affine hyperplane corresponding to the vectors 
 $\mathfrak w^{\vee}_i$, 
 it is satisfied by any affine hyperplane $\mathcal H$
 containing $h(\overline\Gamma')$.
Explicitly, let 
 $\mathfrak w^{\vee}\in (\bar A)^{\perp}\cap N^{\vee}$ be a vector written as
\[
\mathfrak w^{\vee} = a_i\mathfrak w^{\vee}_i+\cdots + 
 a_r\mathfrak w^{\vee}_r
\]
 with $a_i\neq 0$.
Let $\mathcal H_{\mathfrak w^{\vee}}$ be the affine hyperplane
 containing $h(\Gamma')$ parallel to the hyperplane of $N_{\Bbb R}$
 annihilating $\mathfrak w^{\vee}$.
Then
 the set of 1-valent vertices of the connected component of the graph 
 $\Gamma_{\mathcal H_{\mathfrak w^{\vee}}} = 
  h(\Gamma)\cap \mathcal H_{\mathfrak w^{\vee}}$ containing $h(\Gamma')$ 
 is the same as that of the graph $h(\Gamma)\cap \mathcal H_i$, 
 where $\mathcal H_i$ is the affine
 hyperplane containing $h(\Gamma')$ corresponding to the vector 
 $\mathfrak w^{\vee}_i$.
 \qed

 \section{Correspondence theorem for superabundant curves  of genus one II:  Existence of smoothings for general cases}\label{sec:II}
In this section, we extend the results of the previous section
 to not necessarily immersive tropical curves.
\subsection{Calculation of the obstruction for general genus one cases}\label{subsec:general_genus_one}
Here we remove the assumption that $(\Gamma, h)$ is immersive.
\emph{However, we still assume Assumption A} 
 (essentially this does not give restrictions to tropical curves
 which correspond to some classical curves, 
 see Subsection
 \ref{subsec:resol} below).

In this case, some of the edges of $\Gamma$ can be contracted by 
 the map $h$.
This produces higher valent vertices in $h(\Gamma)$.
 
In \cite{N1}, we studied tropical curves with higher valent vertices and
 their associated pre-log curves. 
In particular, we calculated the dual obstruction spaces of such curves,
 extending Theorem \ref{thm:obstruction}.
The existence of higher valent vertices produces several new phenomena.
When there are no higher valent vertices on the loop, 
 there are always pre-log curves of the given type, 
 and the description of the dual obstruction is the same as the immersive cases.
However, higher valent vertices outside the loop affect the calculation of the
 obstruction.

On the other hand, when there are higher valent vertices on the loop, 
 the situation changes rather drastically.
First, there are cases where no pre-log curve of the given type exists, 
 (see \cite[Example 81]{N1}).
Second, these vertices also affect the dual obstruction.

Studying these points carefully, we can deduce the necessary and sufficient condition 
 for the smoothability of tropical curves of genus one even when 
 higher valent vertices exist.
See Theorem \ref{thm:general}.

In this subsection, we study the cases when there are no higher
 valent vertices on the loop part of $h(\Gamma)$ and 
 extend the calculation in the previous section to these cases.
The other cases are studied in Subsection \ref{subsec:4-vloop}.

Before studying a general tropical curve $(\Gamma, h)$, 
 we deal with simple cases where the image has only one
 vertex. 
This is the basic piece of general tropical curves with 
 higher valent vertices which we will study.
 
\subsubsection{Higher valent tropical curves whose images have only one vertex.}\label{subsec:1val}
Let $(\Gamma, h)$, $h\colon \Gamma\to \Bbb R^m$
 be a tropical curve of genus zero 
 (that is, $\Gamma$ is a tree) whose image
 $h(\Gamma)$ has only one vertex.
We assume that this vertex is placed at the origin of $\Bbb R^m$. 
Here $m$ is a positive integer.
Recall that we always assume 
 $\Gamma$ is a 3-valent graph.
Thus, if the valence of the unique vertex $v$ of the image 
 $h(\Gamma)$ is $k$ in the sense of Definition \ref{def:trop},
 the graph $\Gamma$ has $k-2$ vertices and 
 $k-3$ bounded edges. 

Let $\Gamma_0$ be a graph which is combinatorially the same as 
 $\Gamma$ but whose edge weights are all one.
Then the map $h_0\colon \Gamma_0\to \Bbb R^{k-1}$ 
 which maps all the bounded edges and vertices to the origin and
 the unbounded edges to the rays generated by
\[
-e_1 = (-1, 0, \dots, 0),\;\; -e_2 = (0, -1, 0, \dots, 0),
 \;\; -e_{k-1} = (0, \dots, 0, -1),\;\; (1, \dots, 1)
\]
 gives $\Gamma_0$ a structure of a tropical curve.
There is a unique linear map $G\colon \Bbb R^{k-1}\to \Bbb R^m$
 which makes the following diagram commutes:
\[
\xymatrix{
\Gamma_0 \ar[d]_{I} \ar[r]^{h_0} & \Bbb R^{k-1}\ar[d]_G\\
\Gamma \ar[r]_h & \Bbb R^{m} 
}
\]
 here $I$ is the identification of the graphs $\Gamma_0$ and
 $\Gamma$.

Correspondingly, algebraic curves of type $(\Gamma, h)$ 
 (see Definition \ref{def:typev})
 can be obtained from algebraic
 curves of type $(\Gamma_0, h_0)$ by
 the map between toric varieties associated to the map $G$.
Thus, to understand algebraic curves of type $(\Gamma, h)$,
 it suffices to study 
 algebraic
 curves of type $(\Gamma_0, h_0)$.
On the other hand, 
 algebraic curves of type $(\Gamma_0, h_0)$ are linear curves
 in $\Bbb P^{k-1}$, and we can represent these curves very explicitly.
This enables us to study functions on curves corresponding to
 $(\Gamma, h)$ by pulling them back to lines in $\Bbb P^{k-1}$, 
 as we did in the previous section.

\subsubsection{General remarks on the study of tropical curves
 with higher valent vertices}\label{subsec:setup}

Let $(\Gamma, h)$ be a tropical curve of genus one satisfying Assumption A, 
 which does not have higher valent vertices on the loop.
Let $\varphi_0\colon C_0\to \mathfrak X$ be a pre-log curve of type 
 $(\Gamma, h)$.
Note that such a pre-log curve always exists
 (see the beginning part of the proof of Theorem \ref{thm:immersive}).
Our purpose is to calculate the obstruction to lift $\varphi_0$
 to a general fiber of the toric degeneration $\mathfrak X$.

The general setup for this calculation is the same as before.
Namely:
\begin{enumerate}
\item Construct an open cover of $C_0$ consisting of 
 suitable open subsets of two neighboring components of $C_0$.
\item Take a lift of $\varphi_0$ (or more generally of $\varphi_k$
 when it is constructed)
 on each of these open coverings.
Their differences on the overwrappings compose the obstruction class.
\item Calculate the obstruction 
 through analytic continuation of functions
 on a suitable lift of $C_0$.
Using the simple nature of degenerate curves, this can be done explicitly.
\end{enumerate}
Here we make some remarks concerning the differences
 which appear 
 when we deal with higher valent vertices.

The part (1) is the same as before.
As we noted in Subsection \ref{subsec:1val}, a higher valent vertex
 corresponds to a curve in a suitable toric variety and
 it can be obtained as the image of a linear curve in a projective space
 by a map between toric varieties.
Thus, 
 a pair of neighboring components of $C_0$
 corresponds to a pair of neighboring vertices in $h(\Gamma)$, 
 and their suitable open subsets compose an open covering of $C_0$.
 
As for (2), even when the vertices are higher valent, 
 we can prove the existence of trivial lifts as in Lemma \ref{lem:canlift}
 and the canonical correspondence between the set of local sections of the 
 normal sheaf and the set of local lifts.
This enables us to identify the \v{C}ech theoretic obstruction with the
 analytically defined obstruction as in Subsection \ref{subsec:pairing}.

As for (3), 
 such a calculation of the analytic continuation is done by the analogue of the argument in 
 Subsections \ref{subsec:step4} and \ref{subsec:step5}, extended to the cases with higher valent
 vertices. 
Next we discuss this issue.

\subsubsection{The case when $h(\Gamma)$ has 4-valent vertices}\label{subsec:highvalent1}
We first study the case when the tropical curve 
 $(\Gamma, h)$ has 4-valent vertices in its image
 and later we study the case with general higher valent vertices.
We will extend the calculation of 
 Section \ref{sec:I} to these cases.

Consider a 3-valent graph $\Gamma_0$ 
 with two vertices $\alpha_0, \beta_0$, 
 and one bounded edge $E_0$ connecting
 the vertices.
The other edges are non-compact, and let $E_1, E_2$ be the 
 unbounded edges emanating from $\alpha_0$
 and $E_3, E_4$ be the 
 unbounded edges emanating from $\beta_0$.
We take all the edge weights to be one.

\begin{figure}[h]
\includegraphics{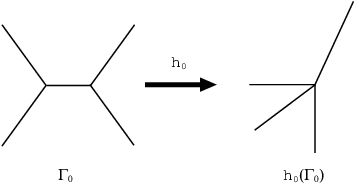}
\caption{}\label{fig:four_valent}
\end{figure}

Then take a map $h_0\colon\Gamma_0\to \Bbb R^3$
 which contracts the bounded edge and gives $\Gamma_0$
 a structure of a tropical curve for which
 the associated toric variety (see Definition \ref{toric})
 defined by the complete fan whose rays are the rays of $h_0(\Gamma_0)$
 is smooth.
This means that we can take the directions of the edges of $h_0(\Gamma_0)$
 to be
\[
(-1, 0, 0),\;\; (0, -1, 0),\;\; (0, 0, -1),\;\; (1, 1,1),
\]
 and the associated toric variety is 
 isomorphic to $\Bbb P^3$. 
The image of $h_0$ is a cone consisting of the apex and four rays
 (see Figure \ref{fig:four_valent}).

Now consider a superabundant
 tropical curve $(\Gamma, h)$ of genus one
 which may have 4-valent vertices on the part away from the loop.
In Subsection \ref{subsec:step6}, 
 we defined a flag of affine subspaces 
\[
A = A_0\subset A_1\subset A_2\subset\cdots
 \subset A_{r-1}\subset A_r = N_{\Bbb R},
\]
 where $A_0$ is the minimal affine subspace containing the loop of
 $h(\Gamma)$ (as before, there may be several ways to choose
 $A_i$, and we fix any one of them).
The construction can be trivially generalized to cases 
 with higher valent vertices and we use the same notation. 
Then let
  $v\in h(\Gamma)$ be a 4-valent vertex
 with the following property.

\begin{itemize}
\item There is some affine subspace
 $A_i$ of $N_{\Bbb R}$ 
 and $v$ is one of the 
 vertices of the component of the subgraph $A_i\cap h(\Gamma)$
 containing the loop
 such that
 some of the edges emanating from $v$ are 
 not contained in $A_i$.
\end{itemize}
\begin{rem}
The vertex $v$ is the analogue of the 1-valent vertex $\alpha$ of
 $\overline\Gamma'$ in 
 the argument of Subsection \ref{subsec:weight one}.
In the current case, such a vertex can be 2-valent as a vertex 
 of the graph $A_i\cap h(\Gamma)$.
See also the argument below.
\end{rem}
Let $\mathfrak{E}_1,\mathfrak E_2,\mathfrak E_3,\mathfrak
  E_4$ be the edges of $h(\Gamma)$
 emanating from $v$ and let $w_1, w_2, w_3, w_4$ be the weights
 of these edges.
Note that two of these edges may coincide.
In this case, the weight of the edge is a pair of integers
 (see Definition \ref{def:trop}).

Let $\Gamma_1$ be the graph topologically isomorphic to $\Gamma_0$, 
 but the edge weights of the unbounded edges 
 are $w_1, w_2, w_3, w_4$
 (the weight of the bounded edge is determined from these weights
 and the combinatorial type of $\Gamma_1$ (see Remark \ref{rem:type}), 
 but we do not need to specify it).
Let $B$ be the minimal
 affine subspace of $N_{\Bbb R}$ containing 
 the edges $\mathfrak{E}_1,\mathfrak E_2,\mathfrak E_3,\mathfrak
  E_4$
  and $\bar B$ be the linear subspace parallel to $B$. 
This subspace $\bar B$ is
 two or three dimensional
 and has a natural integral structure $\bar B_{\Bbb Z}$ coming from 
 the affine lattice $B\cap N$.
Then there is an obvious map $h_1\colon \Gamma_1\to \bar B$
 which contracts the unique bounded edge and
 gives $\Gamma_1$ a structure of a tropical curve
 isomorphic to a neighborhood of the vertex $v$ of $h(\Gamma)$.
We identify (a suitable open subset of) 
 $h_1(\Gamma_1)$ as a subgraph of $h(\Gamma)$
 in the obvious way.
 
By the property of the vertex $v$ above, 
 there are one or two edges of $h_1(\Gamma_1)$ 
 contained in $A_i$
 (if there are three, then all the edges must be contained in $A_i$
 by the balancing condition).
If there is only one edge contained in $A_i$, then take it
 as $h_1(E_1) = \mathfrak E_1$.
If there are two edges contained in $A_i$, 
 let $h_1(E_1) = \mathfrak E_1$ be the edge closer to the loop $h(L)$.
More accurately, there are following three cases:
\begin{enumerate}[(a)]
\item There are two edges of $h_1(\Gamma_1)$ contained in $A_i$.
In this case, there are two subcases:
\begin{enumerate}
\item[(a-1)] The space $B$ is three dimensional.
\item[(a-2)] The space $B$ is two dimensional.
\end{enumerate}
\item Only one edge of $h_1(\Gamma_1)$ is contained in $A_i$.
Let $\mathfrak E_2, \mathfrak E_3, \mathfrak E_4$ be the other edges of 
 $h_1(\Gamma_1)$ (two of them can coincide),
 and let $v_2, v_3, v_4$ be the direction vectors of them.
Then the dimension of the subspace of $N_{\Bbb R}$ 
 spanned by $\bar A_i$ and $\{v_2, v_3, v_4\}$ is $\dim A_i+1$.
Here $\bar A_i$ is the linear subspace of $N_{\Bbb R}$ parallel to $A_i$.
In this case, there are again two subcases:
\begin{enumerate}
\item[(b-1)]  The space $B$ is three dimensional.
\item[(b-2)] The space $B$ is two dimensional.
\end{enumerate}
\item Only one edge of $h_1(\Gamma_1)$ is contained in $A_i$, 
 and the dimension of the subspace 
 spanned by $\bar A_i$ and $\{v_2, v_3, v_4\}$ is $\dim A_i+2$.
In this case, the space $B$ is three dimensional.
\end{enumerate}

We choose the edge $E_1$ of $\Gamma_1$ as above.
Then its image
 $\mathfrak E_1$ must be
 a bounded edge in $h(\Gamma)$.

Let $\varphi_0\colon C_0\to \mathfrak X$ be a pre-log curve
 of type $(\Gamma, h)$, and let 
 $\varphi_{0, v}\colon C_{0, v}\to X_{0, v}\subset\mathfrak X$ be the component
 corresponding to the 4-valent vertex $v$ of $h(\Gamma)$.
Here $X_{0, v}$ is the component of the central 
 fiber of $\mathfrak X$ corresponding to the vertex $v$.
The following is shown by the same argument as in 
 Subsection \ref{subsubsec:leadterm}.
\begin{lem}\label{lem:pullback}
The map $\varphi_{0, v}$ factors through a torically transverse map 
 to the open subset of $\Bbb P^3$ obtained by deleting all
 toric strata of codimension at least two. \qed
\end{lem}
Let $\mathfrak E_1, \mathfrak E_2, \mathfrak E_3$
 be the
 edges 
 emanating from the unique vertex $v$ of $h_1(\Gamma_1)\subset B$
 which are all different.
We can take these since otherwise all the edges emanating from $v$ will be contained
 in $A$. 
Let 
 $u_1, u_2, u_3$ be the primitive integral vectors 
 of the directions of them.
Then the map in Lemma \ref{lem:pullback} can be defined so that
 $\varphi_{0, v}$ is written as $\Phi_{\mathcal L}\circ \varphi_{0, v, s}$, 
 where $\varphi_{0, v, s}\colon
 C_{0, v}\to \Bbb P^3$ is a torically transverse
 map, and $\Phi_{\mathcal L}$ is 
 the map between toric varieties defined on the open subset of $\Bbb P^3$
 and is obtained
 from 
 the linear map
\[
\mathcal L \colon \Bbb Z^3\to \bar B_{\Bbb Z}
\]
 determined by extending the map
\[
(-1, 0, 0)\mapsto w_1u_1,\;\; (0, -1, 0)\mapsto w_2u_2,\;\;
 (0, 0, -1)\mapsto w_3u_3
\]
 linearly (here we take the open subset of $\Bbb P^3$ since in general
 there may not be a map to $X_0$  from 
 the whole $\Bbb P^3$).

Based on this observation, we study each case (a), (b) and (c).\\

\noindent
{\bf Case (a-1).} 
First we assume that two edges $\mathfrak E_1$ and $\mathfrak E_2$
 of $h_1(\Gamma_1)$ are contained in $A_i$.
Also we assume that we have a $k$-th order lift $\varphi_k$ of $\varphi_0$
 and we are going to calculate the obstruction to lift $\varphi_k$ 
 one step further.
As we argued in the proof of Theorem \ref{thm:immersive}, we cannot
 usually expect that the obstruction to lift $\varphi_k$ vanishes.
Instead, if the residues of 
 the leading terms of the pre-obstruction contributed from suitable set of 
 vertices of
 $\Gamma\setminus L$ sum up to zero, then we can perturb $\varphi_k$
 so that the obstruction vanishes.

On the other hand, 
 as in the calculation in Subsection \ref{subsec:step4},
 since the leading term of the pre-obstruction is determined by 
 the data of $\varphi_0$, 
 we do not need to specify what the lift $C_k$ or $\varphi_k$ is, 
 and we concentrate on the calculation of the leading term.

Take a set of vectors $\{\mathfrak v_1, \dots, \mathfrak v_{n-r},
\mathfrak w_1, \dots, \mathfrak w_r\}$ so that 
$\{\mathfrak v_1, \dots, \mathfrak v_{n-r},
\mathfrak w_1, \dots, \mathfrak w_i\}$
is a basis of $\bar A_i\cap N$, $i = 0, 1, \dots, r$.
Let $\{\mathfrak v_1^{\vee}, \dots, \mathfrak v_{n-r}^{\vee}, 
\mathfrak w_1^{\vee}, \dots, \mathfrak w_r^{\vee}\}$
be the dual basis.
As the argument in Subsection \ref{subsec:step4}, 
for the calculation of the leading term of the pre-obstruction
contributed from the vertex $v$, we need to study the
analytic continuation of the functions corresponding to the
vectors $\{\mathfrak w_{i+1}^{\vee}, \dots, \mathfrak w_r^{\vee}\}$
pulled back to the domain curve (see Remark \ref{rem:function} for the precise meaning).

By Lemma \ref{lem:pullback},
 we can pull back the functions
 on $\mathfrak X$ corresponding to the vectors
 $\mathfrak w_j^{\vee}$
  to a linear curve on $\Bbb P^3$. 
Let $x, y, z$ be the affine coordinates on $\Bbb P^3$ corresponding
 to the dual basis of the basis 
\[
\{(-1, 0, 0),\;\; (0, -1, 0),\;\; (0, 0, -1)\}
\]  
 of $\Bbb Z^3$.
Then a torically transverse 
 linear curve in $\Bbb P^3$ is written in the form
\[
ax + y+ b = 0,\;\;
cx + z+ d = 0,
\]
 where $a, b, c, d$ are non-zero constants satisfying $\frac{a}{b}\neq \frac{c}{d}$.
Let $S$ be a parameter on $C_{0, v}$
 given by the pull back of the function $x$.
Then the node of $C_{0}$ corresponding to the edge $E_1$
 is given by the parameter $S = 0$.

On the other hand, choose a basis $\{f_1, f_2, f_3\}$
 of $\bar B_{\Bbb Z}$ in the following way.
Namely, we take $f_1 = u_1$ and
 $f_2$ to be a vector on the plane spanned by 
 $u_1$ and $u_2$, which is contained in $\bar A_i$ by definition.
Then choose $f_3$ to be any vector so that 
 the set $\{f_1, f_2, f_3\}$ becomes 
 a basis of $\bar B_{\Bbb  Z}$.

Using these bases 
 for $\Bbb Z^3$ and 
 $\bar B_{\Bbb Z}$, 
 the map $\mathcal L$ is represented by the matrix
\[
\begin{pmatrix}
w_1 & \alpha & \beta\\
0 & \gamma & \delta\\
0 & 0 & \varepsilon
\end{pmatrix},
\]
 where $\alpha, \beta, \gamma, \delta$ and $\varepsilon$
 are constants and
 $\gamma\varepsilon\neq 0$.
Recall that $w_1$ is the weight of the edge $E_1$.
 
Note that we can take a basis of $N$ in the form 
\[
\{\mathfrak v_1, \dots, \mathfrak v_{n-r},
 \mathfrak w_1, \dots, \mathfrak w_i, 
 \mathfrak w_{i+1} = pf_3+q\mathfrak w', 
 \mathfrak w_{i+2}, \dots, \mathfrak w_{r}\},
\]
 where $\mathfrak w'\in \bar A_i\cap N$ and $p, q$ are rational numbers.
 
Take dual vectors $\mathfrak w_j^{\vee}$ ($i+1\leq j\leq r$) as above.
Since $\mathfrak w_j^{\vee}$ annihilates the subspace $\bar A_i$
 (therefore it annihilates the vectors $f_1$ and $f_2$ under the obvious
 embedding of $\bar B_{\Bbb Z}$ into $N$), 
 it follows that the pull back 
 of the function on $\mathfrak X$ corresponding to 
 $\mathfrak w_{i+1}^{\vee}$ to the linear curve in $\Bbb P^3$ above
 has the form
\[
\epsilon z^{K}
\]
 for some non-zero constant $\epsilon$ and a non-zero integer
 $K$, 
 and the pull back of the functions corresponding to 
 $\mathfrak w_{j}^{\vee}$, $j\geq i+2$ are constants.
Using the parameter $S$ on $C_{0, v}$ as above, 
 $z^K$ becomes
\[
(-1)^{K}(cS+d)^{K}
 = (-1)^{K}(d^{K}+Kd^{K-1}cS+O(S^2))
\]
 up to a constant multiple.
The leading term of the pre-obstruction of the direction $\mathfrak w_j^{\vee}$
 contributed from the vertex $v$ is calculated by the analytic continuation 
 of this function as in Subsection \ref{subsec:weight one}.
Then, since $c$ and $d$ are non-zero, this gives a nontrivial contribution.\\

\noindent
{\bf Case (a-2).}
In this case, the vector $u_2$ satisfies $u_2 = -u_1$.
We take a basis $\{f_1, f_2\}$ of $\bar B_{\Bbb Z}$ as follows.
Namely, take $f_1 = u_1$ as above and take $f_2$ arbitrarily
 so that the pair $\{f_1, f_2\}$ becomes a basis of $\bar B_{\Bbb Z}$.
We use the same basis for $\Bbb Z^3$ (the domain of the map $\mathcal L$)
 as above.

Using these bases, the map $\mathcal L$ is represented by the matrix
\[
\begin{pmatrix}
w_1 & -w_2 & \beta\\
0 & 0 & \gamma
\end{pmatrix}.
\]
Then, taking a basis of $N$ as in Case (a-1),  
 the function on $\mathfrak X$ corresponding to the vector
 $\mathfrak w_{i+1}^{\vee}$ 
 is pulled back to a function of the form
\[
\epsilon z^{K},
\]
 and again the leading term of the pre-obstruction does not vanish.
The functions 
 corresponding to the vectors 
 $\mathfrak w_{j}^{\vee}$, $j\geq i+2$ are pulled back to constants.\\
 
\noindent
{\bf Case (b-1).}
In this case, we  take a basis of $\bar B_{\Bbb Z}$ as follows.
Namely, we take $f_1 = u_1$ as before.
Also, we take $f_2$ from vectors contained in $\bar A_i$.
But in this case the vector $u_2$ is not contained in the 
 linear span of $f_1$ and $f_2$ or 
 $u_1$ and $u_3$.
Then we choose $f_3$ so that the set $\{f_1, f_2, f_3\}$
 becomes a basis of $\bar B_{\Bbb Z}$.
Thus, in this case the matrix of $\mathcal L$ is written in the form
\[
\begin{pmatrix}
w_1 & \alpha & \beta\\
0 & \gamma & \delta\\
0 & \varepsilon & \varphi
\end{pmatrix}.
\]
Note that we have $\varepsilon\neq 0$, $\varphi\neq 0$
 and $\varepsilon+\varphi\neq 0$,
 since otherwise two of the edges emanating from $v$ will be contained in 
 $\bar A_i$.

Take a basis of $N$ as in Case (a).
Then the function on $\mathfrak X$ corresponding to
 $\mathfrak w_{i+1}^{\vee}$
 is pulled back to a function of the form
\[
\epsilon\cdot y^Kz^L,
\]
 for some non-zero constant $\epsilon$ and non-zero integers
 $K, L$ such that $K+L\neq 0$..
Using the parameter $S$ on $C_{0, v}$, this becomes
\[\begin{array}{ll}
(aS+b)^K(cS+d)^L
  &= b^Kd^L+b^{K-1}d^{L-1}(Kad+Lbc)S + O(S^2)\\
  \end{array}
\]
 up to a constant multiple.
Then, 
 as in Proposition \ref{prop:obstruction1},
 the leading term of the pre-obstruction contributed from the vertex $v$
 has the factor $Kad+Lbc$. 
This can be cancelled by
 suitably choosing the configuration of the image of $C_{0, v}$
 (in other words, suitably choosing the non-zero coefficients
 $a, b, c, d$), since $K\neq -L$ and
 the only condition for the constants $a, b, c, d$ is $ad-bc\neq 0$. 

The functions corresponding to the vectors 
 $\mathfrak w_j^{\vee}$, $j\geq i+2$ are pulled back to constants.\\
 
\noindent
{\bf Case (b-2)}.
In this case, we take a basis of $\bar B_{\Bbb Z}$ as follows.
Namely, we take $f_1 = u_1$ as above, and take $f_2$ so that
 the pair $\{f_1, f_2\}$ becomes a basis of $\bar B_{\Bbb Z}$.
Using this basis, the map $\mathcal L$ has the form
\[
\begin{pmatrix}
w_1 & \alpha & \beta\\
0 & \gamma & \delta
\end{pmatrix}.
\]
 with $\gamma\neq 0, \delta\neq 0$ and $\gamma + \delta\neq 0$.
Then as in Case (b-1), we see that the 
 leading term of the pull back of the function corresponding to 
 the vector $\mathfrak w_{i+1}^{\vee}$ can be cancelled when we take 
 suitable configurations of $C_{0, v}$,
 and the functions corresponding to the vectors 
 $\mathfrak w_j^{\vee}$, $j\geq i+2$ are pulled back to constants.\\

\noindent
{\bf Case (c)}.
In this case, we take a basis 
 $\{f_1, f_2, f_3\}$ of $\bar B_{\Bbb Z}$ as follows.
Namely, set $f_1 = u_1$ as before,
 and choose $f_2$ so that it is contained in the plane
 spanned by 
 the vectors $u_1$ and $u_2$.
Then choose $f_3$ so that $\{f_1, f_2, f_3\}$ 
becomes a basis of $\bar B_{\Bbb Z}$.

Using these bases of $\Bbb Z^3$ and $\bar B_{\Bbb Z}$,
 the matrix of $\mathcal L$ has the form
\[
\begin{pmatrix}
w_1 & \alpha & \beta\\
0 & \gamma & \delta\\
0 & 0 & \varphi
\end{pmatrix}, 
\]
 where the constants $\gamma$ and $\varphi$ are non-zero.

Note that we can take a basis of $N$ of the form
\[
\{\mathfrak v_1, \dots, \mathfrak v_{n-r},
 \mathfrak w_1, \dots, \mathfrak w_i, 
 \mathfrak w_{i+1}=af_2+b\mathfrak w', 
 \mathfrak w_{i+2} = cf_2+df_3+e\mathfrak w'', 
 \mathfrak w_{i+3}, \dots, \mathfrak w_{n}\},
\] 
 where $\mathfrak w', \mathfrak w''\in \bar A_i\cap N$
 and $a, b, c, d, e$ are rational numbers.

In this case the function corresponding to the vector 
 $\mathfrak w_{i+1}^{\vee}$ is pulled back to a function of the form
\[
y^Kz^L
\]
 up to a constant multiple, where $K, L$ are integers and $K$ is not
 zero. 
The function corresponding to the vector 
 $\mathfrak w_{i+2}^{\vee}$ is pulled back to
\[
z^{M}
\]
 up to a constant multiple, where $M$ is a nonzero integer.
The functions corresponding to $\mathfrak w_{j}^{\vee}$, $j\geq i+3$
 are pulled back to constants.
As in Case (a), the leading term of $z^M$ cannot vanish
 no matter how we take the configuration of $C_{0, v}$.\\

Summarizing, we see the following.
\begin{prop}
The leading term of the
 pre-obstruction contributed from the vertex $v$ can be cancelled 
 if it is of type (b) above, by suitably choosing $C_0$ and $\varphi_0$.
When it is of type (a) or (c), then it cannot be cancelled unless there 
 are other vertices of $A_i\cap h(\Gamma)$ from which edges not contained in $\bar A_i$
 emanate
 whose distances to the loop $h(L)$ are less than or
 equal to that from the vertex $v$ to the loop.\qed 
\end{prop} 
In the latter case (Cases (a) and (c)),
 the existence of the 4-valent vertex $v$ implies that 
 such a tropical curve is contained in a higher codimensional locus
 in the space of smoothable tropical curves, because 
 the contraction of the bounded edge of $\Gamma_1$ does not
 contribute to the well-spacedness and it just gives a 
 redundant condition to the map $h$. 
Thus, these cases do not contribute to the enumeration problems.\\

Now we study tropical curves with general higher valent vertices.

\subsubsection{The case when $h(\Gamma)$ has 4- or higher valent vertices}
\begin{figure}[h]
\includegraphics{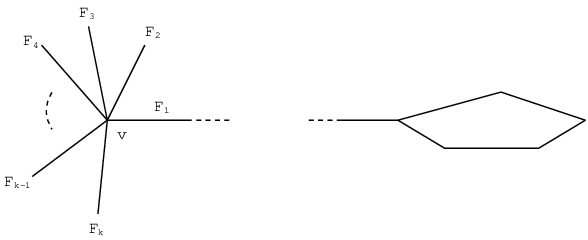}
\caption{$k$-valent vertex $v$ and the loop of
 $h(\Gamma)$.}\label{fig:mult}
\end{figure}

This case can be understood by generalizing the calculation 
 in the case with 4-valent vertices above.
Let $(\Gamma, h)$ be a superabundant tropical curve of genus one 
 in $N_{\Bbb R}$ which has a 
 $k$-valent vertex $v\in h(\Gamma)$, where $k\geq 4$.
As in the previous subsection, we assume that
 $v$ is away from the loop of $h(\Gamma)$. 
The inverse image of $v$ consists of closed subgraphs of 
 $\Gamma$ which are of genus zero.
For simplicity, we assume the inverse image is connected.
Then it is a tree 
 with $k-2$ vertices 
\[
\{V_1, \dots, V_{k-2}\}
\]
 and $k-3$ bounded edges 
\[
\{E_{1}, \dots, E_{k-3}\}
\]
 both of 
 whose ends are contained in the $k-2$ vertices above.
There are $k$ edges 
\[
\{F_1, \dots, F_k\}
\]
 in $\Gamma$
 exactly one of whose ends is contained in 
 $\{V_1, \dots, V_{k-2}\}$, and these edges are mapped to the
 edges emanating from $v$.
Let $\Gamma_2$ be the graph obtained from the union 
\[
\{V_1, \dots, V_{k-2}\}\cup \{E_{1}, \dots, E_{k-3}\}\cup \{F_1, \dots, F_k\}
\]
 by removing one of the ends from each of the edges $F_i$
 which is not contained in $\{V_1, \dots, V_{k-2}\}$.
Let $h_2\colon \Gamma_2\to \Bbb R^n$ be the map giving 
 $\Gamma_2$ the structure of a tropical curve which is 
 obtained by the restriction of $h$ to a neighborhood of 
 $h^{-1}(v)$ in the obvious way.
We assume that the edge $F_1$ is mapped to the 
 edge of $h(\Gamma)$ which is the nearest to the loop
 among the images of $F_i$ (see Figure \ref{fig:mult}).
We regard $h_2(\Gamma_2)$ as a subgraph of $h(\Gamma)$.
\noindent

As before, we make use of a model of a neighborhood of 
 $v$ in $h(\Gamma)$.
Namely, let $\Gamma_0$ be a graph combinatorially identical to
 $\Gamma_2$, 
 but 
 each edge of $\Gamma_0$ has weight one.
We write the edges and vertices of $\Gamma_0$
 by $V_{0, i}, E_{0, i}, F_{0, i}$, parallel to the notations
 for the graph $\Gamma_2$.
Let 
\[
h_0\colon\Gamma_0\to \Bbb R^{k-1}
\]
 be the map which has the following properties.
\begin{itemize}
\item The map $h_0$ contracts the edges
 $\{E_{0, 1}, \dots, E_{0, k-3}\}$ to a point.
In particular, the union 
 $\{V_{0,1}, \dots, V_{0,k-2}\}\cup \{E_{0,1}, \dots, E_{0,k-3}\}$
 is mapped to a point (which we take to be the origin of $\Bbb R^{k-1}$).
\item The map $h_0$ sends the edge 
 $F_{0,i}$ (with one of the ends removed as above)
 to a ray of the direction
\[
e_i = (0, \cdots, 0, \stackrel{i}{\breve{-1}}, 0, \cdots, 0),\;\;
 i = 1, \dots, k-1,
\]
 and
\[
(1, \dots, 1),\;\; i = k.
\]
\end{itemize}
Then the map $h_0$ gives $\Gamma_0$ a structure of a tropical curve.
These conditions fix the map $h_0$ up to isomorphisms of 
 tropical curves.

The image $h_2(\Gamma_2)$ of the
 tropical curve $(\Gamma_2, h_2)$ is obtained from the image of 
 $(\Gamma_0, h_0)$ by the linear map 
\[
G\colon\Bbb R^{k-1}\to \Bbb R^n
\]
 defined by
\[
e_i\mapsto w_iu_i, \;\; i =1, \dots, k-1, 
\]
 where $u_i$ is the primitive integral generator of the direction of 
 the image of the edge $F_i$ by $h$
 and $w_i$ is the weight of $F_i$.

Correspondingly, pre-log curves of type 
 $(\Gamma_2, h_2)$ can be obtained from pre-log curves of 
 type $(\Gamma_0, h_0)$, which are linear curves
 in $\Bbb P^{k-1}$,
 as the image of the map between toric varieties
\[
\Phi_G\colon \Bbb P^{k-1}\to \Bbb P
\]
 associated to the map $G$.
Here $\Bbb P$ is the toric variety which is defined by a fan
 whose rays coincide with the image $h_2(\Gamma_2)$.

Let $\mathcal P$ be the unique shortest path from the vertex $v$ to the loop 
 of $h(\Gamma)$.
Let $A_v$ be the minimal affine subspace of $N_{\Bbb R}$ 
 which contains the union of $\mathcal P$ and the loop of $h(\Gamma)$.
Let $Q_v$ be the minimal affine subspace of $N_{\Bbb R}$ 
 which contains $A_v$
 and all the edges emanating from $v$.
In the argument below, 
 when $P$ is an affine subspace of $N_{\Bbb R}$, 
 then we write by $\bar P$ the linear subspace of $N_{\Bbb R}$
 parallel to $P$.
\begin{prop}\label{prop:higher}
Let $\bar R$ be a subspace of 
 $(\bar A_v)^{\perp}$ in $N_{\Bbb R}^{\vee}$
 containing $(\bar Q_v)^{\perp}$.
Then 
 there is a configuration of $C_{0, v}$ for which the leading term
 of the pre-obstruction contributed from $v$ 
 in any direction contained in 
 $\bar R$ vanishes if and only if $\bar R$
 satisfies the following condition.
\begin{itemize}
\item Let $\mathfrak w^{\vee}$ be any vector in 
 $\bar R\setminus(\bar Q_v)^{\perp}$.
Let $\bar H_{\mathfrak w^{\vee}}$
 be the hyperplane of $N_{\Bbb R}$ annihilated by
 $\mathfrak w^{\vee}$. 
Let $H_{\mathfrak w^{\vee}}$ be the affine hyperplane parallel to 
 $\bar H_{\mathfrak w^{\vee}}$ which contains $A_v$.
Then there are at least three edges emanating from $v$ which are not
 contained in $H_{\mathfrak w^{\vee}}$.
\end{itemize}
\end{prop}
\proof
The proof proceeds along the same line as in the calculation 
 of Subsection \ref{subsec:highvalent1}.
We use the same notation 
 introduced in the paragraphs before the proposition.

First, let $\mathfrak w^{\vee}$ be a vector in $\bar R\cap N^{\vee}$
 such that there are only two edges emanating from 
 $v$ which are not contained in $H_{\mathfrak w^{\vee}}$.
In this case, take a basis $\{e_1, \dots, e_{k-1}\}$ of $\Bbb Z^{k-1}$ 
 so that the vectors $e_i$, $i = 1, \dots, k-2$, 
 are mapped by $G$ to vectors contained in 
 $\bar H_{\mathfrak w^{\vee}}$, and 
 take a basis $\{f_1, \dots, f_n\}$ of $N$ so that 
 $\{f_1, \dots, f_{n-1}\}$ is a basis of $\bar H_{\mathfrak w^{\vee}}$.

Then the $n\times (k-1)$ matrix representing $G$ has
 the form
\[
\begin{pmatrix}
*& \cdots & \cdots & * & *\\
*& \cdots & \cdots & * & *\\
\vdots & \vdots & \vdots & \vdots &\vdots \\
*& \cdots & \cdots & * & *\\
0& \cdots & \cdots & 0 & a
\end{pmatrix}
\]
 where $a$ is not zero.

Note that a pre-log curve of type $(\Gamma_0, h_0)$
 can be represented by defining equations of the form
\[
x_1+a_2x_2-b_2 = 0,\;\; x_1+a_3x_3-b_3 = 0, \;\;
 \cdots, \;\; x_1+a_{k-1}x_{k-1}-b_{k-1} = 0,
\]
 where $\{x_1, \dots, x_{k-1}\}$ is the inhomogeneous coordinate
 system of $\Bbb P^{k-1}$ corresponding to the dual basis of 
 $\{e_1, \dots, e_{k-1}\}$,
 and $a_i, b_i$ are nonzero constants satisfying ${b_i}\neq {b_j}$,
 $\forall i, j$.
Let us write by $\psi_0\colon C_{0, v}\cong\Bbb P^1\to \Bbb P^{k-1}$
 a parameterization of this linear curve.   
We can assume that $x_1$ is zero at the node corresponding to the edge 
 $F_1$.
We use the pull back of $x_1$ to $C_{0, v}$ as an affine coordinate
 on $C_{0, v}$ and write it by $S$.

Let $\{f_1^{\vee}, \dots, f_n^{\vee}\}$ be the dual basis of 
 $\{f_1, \dots, f_n\}$.
Note that $\mathfrak w^{\vee}$ is proportional to $f_n^{\vee}$.
Then the function on $\mathfrak X$ corresponding to $f^{\vee}_n$
 in the sense of Remark \ref{rem:function}
 is pulled back by $\Phi_G\circ\psi_0$ to a function of the form
\[
\varepsilon \left(-\frac{S-b_{k-1}}{a_{k-1}}\right)^L,
\]
 where $\varepsilon$ is a nonzero constant and $L$ is a nonzero integer. 

As in Case (a) of Subsection \ref{subsec:highvalent1},
 the analytic continuation of this function
 to the loop part will give a nonzero contribution
 to the leading term of the pre-obstruction
 of direction $\mathfrak w^{\vee}$.
This proves that the condition in the statement of the proposition
 is necessary for the vanishing of the leading term of the pre-obstruction
 contributed from $v$ in any direction
 which belongs to
 the subspace $\bar R$.

Conversely, assume that the condition in the statement of the proposition
 is satisfied for a subspace $\bar R$.
Take a set of vectors in $N$ in the following way.
Namely, let $\{f_1, \dots, f_s\}$ be a basis of 
 the subspace $\bar A_v$.
Choose $f_{s+1}$ from one of the directions of the edges
 emanating from $v$ which is not contained in $\bar A_v$.
Next, choose $f_{s+2}$ so that it is another one of the 
 directions of the edges
 emanating from $v$ which is not contained in the span of 
 $\bar A_v$ and $f_{s+1}$.
Then inductively choose $f_{s+3}, f_{s+4}, \dots, f_{s+l}$
 so that $f_{s+j}$ is one of the 
 directions of the edges
 emanating from $v$ which is not contained in the span of 
 $\bar A_v$ and $f_{s+1}, f_{s+2}, \dots, f_{s+j-1}$.
Here $s+l = \dim Q_v$.
Then add vectors $f_{s+l+1}, \dots, f_n$ so that 
 $\{f_1, \dots, f_n\}$ becomes a basis of $N_{\Bbb Q}$.
Let $\{f_1^{\vee}, \dots, f_n^{\vee}\}$ be the dual basis of 
 $\{f_1, \dots, f_n\}$ in $N_{\Bbb Q}^{\vee}$.

We take a basis 
\[
\{g_1, \dots, g_b\}
\]
 of $\bar R$ in the following way, here $b = \dim \bar R$.
First, note that $\{f_{s+l+1}^{\vee}, \dots, f_n^{\vee}\}$
 is a basis of $(\bar Q_v)^{\perp}$.
We set $g_1 = f_{s+l+1}^{\vee}, \dots, g_{n-s-l} = f_n^{\vee}$.
Note that each of the hyperplanes $\bar H_{g_i}$, $i = 1, \dots, n-s-l$
 contains a parallel transport of the union of the loop of $h(\Gamma)$, 
 the path $\mathcal P$ and the edges emanating from $v$.
We write one of such parallel transports of the union by $\bar U$. 
We use the same notation ($v, F_i, \mathcal P$, etc.)
 for the corresponding objects of $\bar U$. 

Now consider the flag
\[
\begin{array}{ll}
(\bar Q_v)^{\perp} = \langle g_1, \dots, g_{n-s-l}\rangle &
 \subset \langle g_1, \dots, g_{n-s-l}, f_{s+l}^{\vee}\rangle\\
 & \subset \langle g_1, \dots, g_{n-s-l},
  f_{s+l}^{\vee}, f_{s+l-1}^{\vee}\rangle\\
 & \subset \cdots \\
 & \subset \langle g_1, \dots, g_{n-s-l}, f_{s+l}^{\vee}, f_{s+l-1}^{\vee}, 
 \dots, f_{s+1}^{\vee}\rangle
  = (\bar A_v)^{\perp}.
\end{array}
\]
 in $N_{\Bbb R}^{\vee}$,
 here for vectors $h_1, \dots, h_j$,
 the notation
 $\langle h_1, \dots, h_j\rangle$ means the subspace spanned by
 them.
We write this flag by
\[
C_0 = (\bar Q_v)^{\perp} \subset C_1\subset\cdots\subset
 C_l = (\bar A_v)^{\perp}.
\]
The intersection of this flag with $\bar R$ gives a flag 
 in $\bar R$ which we write by
\[
C_0\cap \bar R = \bar R_0 \subset
 \bar R_1\subset\cdots \subset \bar R_{b-(n-s-l)} = \bar R.
\]
Note that $\bar R_0 = C_0$ since $\bar R$ contains
 $C_0 = (\bar Q_v)^{\perp}$ by assumption.
Then take a basis $\{g_1, \dots, g_b\}$ 
 of $\bar R$ defined over $\Bbb Q$ so that 
 $\{g_1, \dots, g_{n-s-l+i}\}$ is a basis of $\bar R_i$.

By definition of $f_i$ and $g_i$, 
 if $f_{s+l-j}^{\vee}\in C_{j+1}\cap \bar R = \bar R = \bar R_{b-(n-s-l)}$ and 
 $f_{s+l-j}^{\vee}\notin C_{j}\cap \bar R = \bar R_{b-(n-s-l)-1}$ for some $0\leq j\leq l-1$, 
 then the edge emanating from $v$ whose direction $f_{s+l-j}$
 is contained in $(\bar R_{b-(n-s-l)-1})^{\perp}$, but not contained in 
 $(\bar R_{b-(n-s-l)})^{\perp}$.

Inductively, 
 one sees that 
 the basis $\{g_1, \dots, g_b\}$ of $\bar R$ has the following 
 property.
\begin{lem}\label{lem:flag}
For any $j$ with $1\leq j\leq b-(n-s-l)-1$, 
 there is at least one edge emanating from $v$
 which is contained in $(\bar R_j)^{\perp}$ but not contained in 
 $(\bar R_{j'})^{\perp}$ with $j<j'$.\qed
\end{lem}


Multiply the vectors $g_i$, $1\leq i\leq b$, by suitable integers so that 
 they belong to $N^{\vee}$ and write the result
 by $\{g_1', \dots, g_b'\}$.
We consider the pull back of the functions on $\mathfrak X$
 corresponding to the vectors 
 $g_i'$  (in the sense of Remark \ref{rem:function})
 by the map $\Phi_G\circ\psi_0$ to functions on $C_{0, v}$.
We use the same parameter $S$ of $C_{0, v}$ as above.

First, the functions corresponding to $g_i'$, 
 $1\leq i\leq n-s-l$ are pulled back to constant functions.
Therefore, they do not contribute to the leading term of the pre-obstruction.

Consider the hyperplane $H_{g_{b}}$ in $N_{\Bbb R}$.
Since $g_{b}\in \bar R$, there are at least three edges emanating
from $v$ which are not contained in 
$H_{g_{b}}$ by assumption.
Then it follows 
 from the same argument as in Subsection \ref{subsec:highvalent1}
 that
 the function corresponding to $g_{b}$
 is pulled back to a function of the form
\[
\varepsilon_b \left(-\frac{S-b_{b, 1}}{a_{b, 1}}\right)^{L_{b, 1}}
\left(-\frac{S-b_{b, 2}}{a_{b, 2}}\right)^{L_{b, 2}}\cdots
\left(-\frac{S-b_{b, i_b}}{a_{b, i_b}}\right)^{L_{b, i_b}}
\]
 for some $i_b\geq 2$.
Here $\varepsilon_b$ is a nonzero constant,
 $L_{b, j}$ are nonzero integers
 and $\{a_{b, j}\}$ and $\{b_{b, j}\}$ are 
 subsets of $\{a_2, \dots, a_{k-1}\}$ and $\{b_2, \dots, b_{k-1}\}$.
 
Expanding this as a series of $S$ near $S = 0$, 
 the term of order one with respect to $S$ contributes to the leading term
 of the pre-obstruction contributed from $v$ of the direction 
 $g_{b}$.
It is easy to see that there are nonzero numbers
 $\{a_{b, j}\}$ and $\{b_{b, j}\}$ such that the coefficient
 of $S$ in the expansion vanishes. 
Note that the situation where $i_b = 2$ and $L_{b, 1} +L_{b, 2} = 0$
 does not happen, since in this case the number of edges not contained
 in $H_{g_b}$ will be only two.
This assures the existence of $\{a_{b, j}\}$ and $\{b_{b, j}\}$ satisfying
${b_{b, j}}\neq {b_{b, j'}}$, when $j\neq j'$.
This means that if we take the configuration of the image of $C_{0, v}$
 so that one of the components of its inverse
 image by $\Phi_G$ is defined by the defining equations
 $x_1+a_2x_2-b_2 = 0,\;\; x_1+a_3x_3-b_3 = 0, \;\;
 \cdots, \;\; x_1+a_{k-1}x_{k-1}-b_{k-1} = 0$ with the property that the
 part $\{\{a_{b, j}\}, \{b_{b, j}\}\}$ of the coefficients 
 satisfies the above condition, then the leading term of the
 pre-obstruction in the direction $g_{b}$ vanishes.

Now consider the pull back of the function on $\mathfrak X$
 corresponding to the vector $g_{b-j}$, $1\leq j\leq b-(n-r-l)-1$.
It has the similar form as above:
\[
\varepsilon_{b-j} \left(-\frac{S-b_{b-j, 1}}{a_{b-j, 1}}\right)^{L_{b-j, 1}}
\left(-\frac{S-b_{b-j, 2}}{a_{b-j, 2}}\right)^{L_{b-j, 2}}\cdots
\left(-\frac{S-b_{b-j, i_{b-j}}}{a_{b-j, i_{b-j}}}\right)^{L_{b-j, i_{b-j}}}.
\]
Since the hyperplane $H_{g_{b-j}}$ also satisfies the property that
 there are at least three edges emanating
 from $v$ which are not contained in 
 $H_{g_{b-j}}$, we have $i_{b-j}\geq 2$ as before.

The important point is that by Lemma \ref{lem:flag}, 
 there is at least one pair $(a_{b-j, p}, b_{b-j, p})$ which does not
 appear as the pair $(a_{b-j', p'}, b_{b-j', p'})$
 for $j'<j$, $1\leq p'\leq i_{b-j'}$.
Then, as in the case for $H_{g_b}$ above,
 it follows that one can cancel the coefficient of $S$ in the 
 expansion of the function above by suitably choosing nonzero
 coefficients
 $a_{b-j, p}$ and $b_{b-j, p}$, without changing the other coefficients
 $a_{b-j', p'}$ and $b_{b-j', p'}$.

By induction, this implies that one can choose the configuration of 
 $C_{0,v}$ (in other words, the coefficients $\{a_{i}\}$ and $\{b_{i}\}$
 of the above equations) so that
 the contributions to the leading terms of 
 the pre-obstructions from the vertex $v$ of the directions
 $g_1, \dots, g_{n-r-l}, g_{n-r-l+1}, \dots, g_{b}$ simultaneously vanish.
Note that the leading term of the pre-obstruction 
 contributed from $v$ is linear with respect to the vectors 
 $g_{n-r-l+1}, \dots, g_{b}$.
This follows from the same argument as in the proof of Lemma \ref{lem:linearity},
 since for any direction $\mathfrak w^{\vee}$, the leading term has the 
 same order with respect to $t$ (which is the path length of $\mathcal P$
 in the sense of Definition \ref{def:path length}).

Although the set 
 $\{g_1, \dots, g_{n-r-l}, g_{n-r-l+1}, \dots, g_{b}\}$
 is not an integral basis of $\bar R$ 
 in general, by this linearity, 
 it follows that the leading term of the pre-obstruction contributed from 
 the vertex $v$ in any direction contained in $\bar R$ vanishes
 for this configuration of $C_{0, v}$.
This proves the proposition. \qed\\

Even when the leading terms of 
 the pre-obstructions contributed from individual higher valent vertices
 cannot be cancelled, 
 by the argument used in the proof of
 this proposition, one can evaluate the
 leading terms contributed
 from these vertices which are not on the loop of $h(\Gamma)$
 as in Subsections \ref{subsec:step4} and \ref{subsec:step5},
 where we did the calculation for the cases corresponding to
 immersive tropical curves.
However, these cases are not general in the space of smoothable tropical curves
 as we mentioned at the last of Subsection \ref{subsec:highvalent1}.
In the next subsection, we recall the case with higher valent vertices
 on the loop from  \cite{N1}.

\subsection{4- or higher valent vertices on the loop}\label{subsec:4-vloop}
Here we consider the case when the superabundant
 tropical curve
 $(\Gamma, h)$ of genus one has 4- or higher valent vertices on the loop.
We assume Assumption A as usual.
This case was studied in detail in the paper \cite{N1}, 
 and degenerate curves become smoothable 
 by a reason very different from the calculation we have done so far 
 in this paper.
Now we recall results of \cite{N1} relevant to the current situation.

\subsubsection{Smoothability of pre-log curves corresponding to tropical curves of genus one with higher valent vertices on the loop}\label{subsec:loopvertobst}
First we recall the result on the smoothability of pre-log curves corresponding to 
 tropical curves of genus one with higher valent vertices on the loop, 
 assuming the existence of such pre-log curves.

\begin{thm}\cite[Theorem 81, Remark 82]{N1}\label{thm:highvalloop}
Let $h\colon \Gamma\to\Bbb R^n$ be a tropical curve of genus one
 which satisfies Assumption A.
Let $\{F_j\}$ be the set of flags of $(\Gamma, h)$
 whose vertices are on the loop of $h(\Gamma)$.
We write by $v_j\in N$ the direction vector of $F_j$. 
Let $\bar B$ be the subspace of 
 $N_{\Bbb R}$ spanned by these directions
 $v_j$.
 
Assume that there exists a pre-log curve 
 $\varphi_0\colon C_0\to X_0$ of type $(\Gamma, h)$ 
 in the central fiber of a suitable toric degeneration $\mathfrak X\to \Bbb C$.
Then if the obstruction of the direction 
 contained in $(\bar B)^{\perp}\cap N^{\vee}$ (see Definition \ref{def:preob})
 to deform $\varphi_0$ 
 to a general fiber of $\mathfrak X$
 vanishes, the whole obstruction vanishes.
In other words, we only need to calculate the obstruction 
 whose direction is contained in 
 $(\bar B)^{\perp}\cap N^{\vee}$.
 \qed
\end{thm}
Even when there are higher valent vertices on the loop, 
 the calculation of the obstruction in the direction contained in $(\bar B)^{\perp}\cap N^{\vee}$
 can be done in the same way as in the previous subsection.

\begin{rem}\label{rem:immersive}
\begin{enumerate}
\item In \cite[Proposition 68]{N1}, we proved a general theorem 
 for the existence of pre-log curves associated to regular
 tropical curves.
When tropical curves are superabundant, 
 there are such curves for which pre-log curves of the given type do
 not exist, even when the tropical curves are embedded, as
 shown in \cite[Example 82]{N1}.
\item However, for an embedded (or immersed)
 tropical curve of genus one, it is not difficult to see that pre-log curves
 of the given type exist even if the tropical curve is superabundant.
This is sketched at the beginning of the proof of Theorem \ref{thm:immersive}.
\item For non-immersive tropical curves of genus one, there are 
 both cases where pre-log curves of the given types
 exist or do not exist, as shown in \cite[Example 81, Example 83]{N1}, 
 and we recall the general result in the next subsection.
The above theorem shows that the existence of 
 such a pre-log curve automatically implies the vanishing of 
 a part of the obstruction associated to the loop of the tropical curve. 
\end{enumerate}
\end{rem}
To see the importance of the above theorem to the correspondence
 type results, we consider an example shown to the author by
 Carolin Torchiani.
\begin{example}\label{ex:loop}
Consider tropical curves of genus one in $\Bbb R^3$ 
 with four unbounded edges of directions 
\[
v_1 = (-1, 1, -1), \;\; v_2 = (2, 1, 0),\;\; v_3 = (1, -2, -1),\;\; v_4 = (-1, 0, 1).
\]
Here the edge with the direction $(-1, 0, 1)$ has weight two and
 the other edges have weight one.
Consider the set of such curves which pass through the points
\[
P_1 = (0, 0, 0),\;\; Q_1 = (1, 4, 2).
\]
By direct calculation, one sees that there is only one
 such tropical curve given in Figure \ref{fig:type2}.
This tropical curve is regular and has multiplicity six
 (we can define the multiplicity of higher genus curves 
 with appropriate incidence conditions in the same way as in 
 \cite{NS}, so long as the tropical curve is regular).

\begin{figure}[h]
\includegraphics[width=5cm]{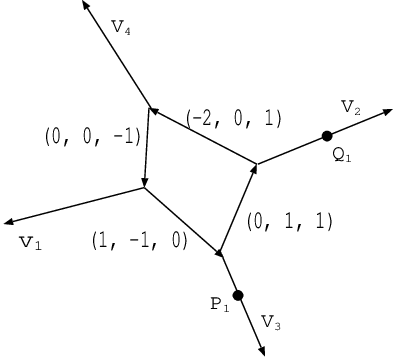}
\caption{}\label{fig:type2}
\end{figure}

On the other hand, consider the set of those tropical curves 
 which pass through the points
\[
P_2 = (0, 0, 0),\;\; Q_2 = (-2, -4, -1).
\]
Then there is one regular tropical curve
 given in Figure \ref{fig:type1}.
 
\begin{figure}[h]
\includegraphics[width=6cm]{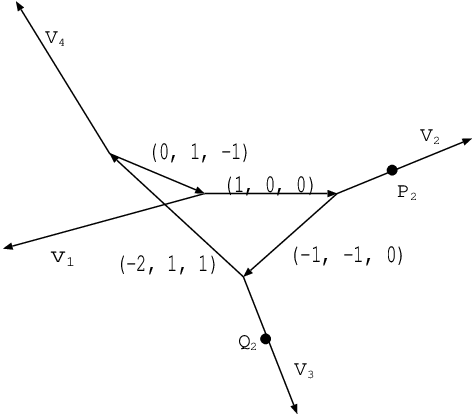}
\caption{}\label{fig:type1}
\end{figure}

However, this tropical curve has multiplicity two.
Therefore, there should be some other curves which compensate for the
 missing multiplicity four.
In fact, there is a family of superabundant tropical curves 
 given in Figure \ref{fig:type1'}. 
 
\begin{figure}[h]
\includegraphics[width=9cm]{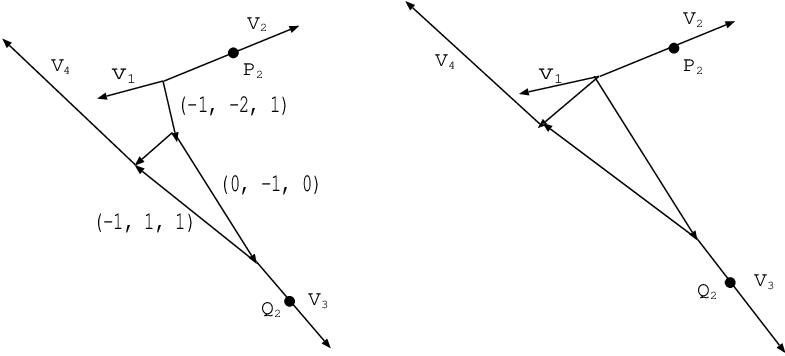}
\caption{The loop part can shrink or extend without breaking the 
 incidence conditions.}\label{fig:type1'}
\end{figure}
\end{example}

Using the calculations so far, 
 one sees that 
 all these superabundant tropical curves are not smoothable,
 except the one which has a 4-valent vertex on the loop
 (the picture on the right hand side of Figure \ref{fig:type1'}).
 
Let us write this tropical curve with a 4-valent vertex by $(\Gamma, h)$.
By Theorem \ref{thm:highvalloop}, a pre-log curve of type
  $(\Gamma, h)$ is, if it exists, smoothable.
Therefore, the problem is reduced to the existence of such pre-log curves.
According to Theorem \ref{thm:119}
 below, it follows that there is a unique family of pre-log curves
 of this type.
 
Moreover, we can define the multiplicity of the tropical curve
 $(\Gamma, h)$ as in the regular case
 (since now the length of the edge of direction $(-1, -2, 1)$
 is fixed to zero, the freedom to deform $(\Gamma, h)$
 is four, matching with the number of constraints given by the 
 points $P_2$ and $Q_2$), 
 and it turns out to be four.
 
On the side of pre-log curves, note that the 4-valent
 vertex corresponds to a rational curve with four special points.
Given a pre-log curve of type $(\Gamma, h)$ as above, 
 its component corresponding to the 4-valent vertex determines 
 a point in the moduli space of rational curves with four special points.
Fixing this point, we can transport each component by the torus action of the ambient space.
If these actions are compatible in the sense that the nodes of neighboring 
 components are kept, this gives a deformation of the pre-log curve.
In this way, we obtain a $(\Bbb C^*)^4$-torsor of pre-log curves
 of type $(\Gamma, h)$
 (this is the "family of pre-log curves" in the statement of Theorem \ref{thm:119}
 below. See Remark \ref{rem:numberofcurves}.).
See \cite[Subsection 5.2]{N1} for more details of this construction.

Then the same argument as in \cite[Proposition 5.7]{NS} shows there are
 four pre-log curves in this family satisfying constraints corresponding to $P_2$
 and $Q_2$.
As we noted above, each of them is smoothable to a general fiber keeping the constraint
 conditions.
Thus, we again obtain six smoothable curves matching the constraints.
\qed

\subsubsection{Existence of pre-log curves corresponding to tropical curves of genus one with higher valent vertices on the loop}\label{subsec:loopexist}
As we noted in Remark \ref{rem:immersive}, 
 in the case of regular tropical curves, 
 and also in the case of immersive superabundant 
 tropical curves of genus one, 
 the existence of pre-log curves of that type is unconditionally assured.
For general superabundant tropical curves, the study of the existence of 
 corresponding pre-log curves  occupies the half of 
 the problem of correspondence between tropical and algebraic curves
 (the other half is the deformation of the pre-log curves).
In \cite{N1}, we developed a formalism to study the existence of
 such pre-log curves. 
In particular, for the case of genus one,
 we obtained the general answer.
We recall it in this subsection.

To state the result, we fix some notations.
Let $(\Gamma, h)$, $h\colon \Gamma\to \Bbb R^n = N\otimes\Bbb R$
be a tropical curve of genus one satisfying Assumption A.
Let $\tilde L\subset h(\Gamma)$ be the union of all the open edges emanating from 
vertices on the loop 
 and the loop itself (in other words, $\tilde L$ is the
'open star' of the loop part of $h(\Gamma)$).

The part $\tilde L$ can be thought of as a tropical curve by extending the 
open edges to infinity.
We write it by $(\Gamma', h')$.
If there is a pre-log curve of type $(\Gamma', h')$, it is easy to 
see that there is also a pre-log curve of type $(\Gamma, h)$, 
since $\Gamma\setminus \Gamma'$ is a union of trees.
Thus, we assume there are no vertices outside the loop, so that
$h(\Gamma)$ equals $\tilde L$.
We also assume that the image $h(\Gamma)$ is not contained
in a proper affine subspace of $\Bbb R^n$.

Let $A$ be the minimal affine subspace of $N_{\Bbb R}\cong \Bbb R^n$ containing 
the loop $h(L)\subset h(\Gamma)$.
Let $\bar A$ be the linear subspace of $N_{\Bbb R}$ parallel to $A$.
Take a basis $\{e_1, \dots, e_n\}$ of $N$ so that 
$\{e_1, \dots, e_r\}$ is a basis of $\bar A\cap N$.
Let $\{e_1^{\vee}, \dots, e_n^{\vee}\}$ be the dual basis of $\{e_1, \dots, e_n\}$.

Let $\{v_1, \dots, v_k\}$ be the vertices on the loop of $h(\Gamma)$.
Among the edges emanating from $v_i$, let
 \[
 E_{v_i, 1}, \dots, E_{v_i, m_i},\;\; 1\leq i\leq k,
 \]
 be those not contained in $A$ (if not empty).
 We choose $m_i-1$ edges $E_{v_i, 1}, \dots, E_{v_i, m_i-1}$
 for each $i = 1, \dots, k$.
 Let $\mathcal E$ be the set of these edges.

\begin{rem}
There are choices of $m_i-1$ edges from 
 the $m_i$ edges $E_{v_i, 1}, \dots, E_{v_i, m_i}$
 for each vertex $v_i$, $1\leq i\leq k$.
However, these choices do not affect the result thanks to the
 balancing condition.
\end{rem}

Then we construct a flag of affine subspaces of $N_{\Bbb R}$
containing $A$ in the following way.
First choose an edge $E_{v_i, q}$ from $\mathcal E$.   
Let $A_1$ be the minimal affine subspace containing both $A$ and $E_{v_i, q}$.
Let $n_1$ be the number of the edges in $\mathcal E$ contained in $A_1$.
Next, choose an edge $E_{v_j, r}$ from $\mathcal E$ which is not contained
in $A_1$.
Let $A_2$ be the minimal affine subspace containing both $A_1$ and $E_{v_j, r}$.
Let $n_2$ be the number of the edges in $\mathcal E$ contained in $A_2$.

Continuing this process, we can construct a flag
\[
A = A_r\subset A_{r+1}\subset \cdots \subset A_{n-1}\subset A_{n} = N_{\Bbb R}.
\]
There is also an associated sequence of the numbers
\[
n_{r+1}<n_{r+2}<\cdots < n_{n-1}<n_{n} = \sum_{i=1}^k(m_i-1).
\]

Among these flags, we choose one for which the sequence of the positive integers
\[
n_{r+1}, n_{r+2}-n_{r+1}, \dots, n_{n}-n_{n-1}
\]
is the minimal with respect to the lexicographic order
(in other words, we choose a flag so that $n_{i+1}-n_i$
with smaller $i$ will be as small as possible). 
There may be several flags which give the same sequence of the integers, 
and any of them will do.
We write such a flag by   
$A_r\subset A_{r+1}\subset \cdots \subset A_{n-1}\subset A_{n}$. 

In this paper, we assume that all the numbers
 $n_{r+1}, n_{r+2}-n_{r+1}, \dots, n_{n}-n_{n-1}$ are one for simplicity.
We also assume that there are no edges emanating from $v_1, \dots, v_k$ 
 contained in the affine subspace $A$ except those contained in the loop.
This is the case relevant to enumerative problems.
See \cite[Subsection 4.2.2]{N1} for general cases.
In this case, $n_n = n-r$.

Now choose a basis $\{g_{r+1}, \dots, g_n\}$ of the sublattice 
 $N^{\vee}\cap \bar A^{\perp}$ 
 of the
 dual lattice $N^{\vee}$  
 in the following way.
First, let $g_{n}$ be one of the generators of the annihilating subspace
 of $\bar A_{n-1}$. 
Next, take $g_{n-1}$ so that $\{g_n, g_{n-1}\}$ will be a basis of the annihilating
subspace of $\bar A_{n-2}$.
Similarly, take $\{g_{r+1}, \dots, g_{n}\}$ so that 
$\{g_{n-i+1}, \dots, g_{n}\}$ is a basis of the annihilating subspace
of $\bar A_{n-i}$.

We order the edges $E_{1}, \dots, E_{n-r}$
so that the first $i_1 ( = m_1-1)$ edges of them emanates from the vertex $v_1$, 
the next $i_2-i_1$ edges emanates from $v_2$, and so on.
We write this reordered sequence of the edges by the same letters.
Note that some of $i_{j+1}-i_j ( = m_{j+1}-1)$ might be zero.

Let $v_m$ be the primitive integral generator of the 
direction of the edge $E_m$ and $w_m$ be the weight of it.  
The values of 
\[
w_mg_j(v_m), \;\; j = r+1, \dots, n
\]
give an element of $\Bbb Z^{n-r}$ for 
$m = 1, \dots, n-r$.   
Taking these vectors as column vectors, we make a
 square matrix $G$ of size $n-r$.
By construction, $G$ is invertible 
 (the inverse matrix may not be defined over the integers).

Let $I^{n-r} = (0, 1)^{n-r}$ be a higher dimensional 
 open unit cube.
Since the matrix $G$ is invertible, the image of $I^{n-r}$ by $G$ is a 
 higher dimensional open parallelotope $P\subset\Bbb R^{n-r}$
 spanned by the columns of 
 $G$.
For $a\neq b$, 
 $\sum_{c=1}^di_c+1\leq a, b\leq\sum_{c=1}^{d+1}i_c$ for some $0\leq d\leq k-1$, 
 let $H_{ab}$ be the hyperplane in $\Bbb R^{n-r}$ defined by
 \[
 H_{ab} = \left\{\begin{pmatrix} X_1\\ \vdots \\ X_{n-r}\end{pmatrix} \; \bigg| \; X_a = X_b\right\}
 \]
 and $H$ be their union
 \[
 H = \coprod H_{ab},
 \]
 where the pair $a, b$ runs through those which satisfy the above condition.
Then we have the following.
We use the same notation as in the above argument.
 
\begin{thm}\label{thm:119}$($\cite[Corollary 119]{N1}$)$
Given a tropical curve $(\Gamma, h)$ of genus one,
 there are pre-log curves of type $(\Gamma, h)$ 
 if and only if the set of integral points in $P\setminus G(H)$ is not empty.
 Moreover, the number of families of such pre-log curves is 
 given by the number of integral points in $P\setminus G(H)$.\qed
\end{thm}

 \begin{rem}\label{rem:numberofcurves}
 	\begin{enumerate}
 		\item The number of the isomorphism classes of the domain curve $C_0$
 		of pre-log curves $\varphi_0\colon C_0\to X_0$
 		of type $(\Gamma, h)$ is finite and given by the number of integral points in $P\setminus G(H)$.
 		\item The family of pre-log curves corresponding to an integral point in $P\setminus G(H)$
 		consists of curves obtained by gluing marked rational curves with fixed moduli associated 
 		to the vertices of $\Gamma$ (the marked points correspond to nodes or intersection
 		with toric divisors) as we explained at the last of Subsection \ref{subsec:loopvertobst}. 
 		\item
 		When $n_n$ is larger than $n-r$, and when there are pre-log curves of type $(\Gamma, h)$, 
 		then there is a continuum of isomorphism classes of the domain curves
 		of these pre-log curves.
 	\end{enumerate}
 \end{rem}

\subsection{Main result}\label{subsec:mainresult}
Let $(\Gamma, h)$ be an immersive 3-valent tropical curve of genus one  
 of fixed degree $\Delta$.
In the case of genus one, the degree
 $\Delta$ determines the expected dimension of
 the parameter
 space (without referring to
 the dimension of the ambient space). 
Namely, the expected dimension of the parameter space 
 of tropical curves of the fixed degree including 
 $(\Gamma, h)$
 is $|\Delta|$, the number of unbounded edges of 
 $\Gamma$. 

Let us take a toric degeneration $\mathfrak X$
 defined respecting $(\Gamma, h)$ (see Definition \ref{def:degeneration}).
Let $\varphi_0\colon C_0\to X_0$ be a pre-log curve of type 
 $(\Gamma, h)$, where $X_0$ is the central fiber of the
 degeneration $\mathfrak X$.
Let $r$ be the dimension of the space 
 $H^1(C_0, \mathcal N_{C_0/X_0})$.
The dimension of the parameter space of the tropical curve
 in the sense of Proposition \ref{prop:trop_moduli} is,
 if it is not empty, given by
 the dimension of the space $H^0(C_0, \mathcal N_{C_0/X_0})$
 (this claim is essentially proved in the proof of \cite[Proposition 7.3]{NS}.
 See also \cite[Section 5.2]{N1} for details.), 
 and the expected dimension above is given by
\[
\dim H^0(C_0, \mathcal N_{C_0/X_0}) - r. 
\]
Thus, the dimension of the parameter space equals
 $r+|\Delta|$.
\begin{rem}
In general, 
 the dimension of 
 the parameter space of tropical curves of a given type
 can be smaller than expected, see \cite[Example 19]{N1}.
However, if the parameter space contains at least one 
 immersive tropical curve, as we assume here, 
 then the parameter space has dimension equal to
 or greater than the expected dimension.
\end{rem}
 
To state our main result (Theorem \ref{thm:general} below), 
 we extend the notion of well-spacedness given in Definition
 \ref{well-spaced} (which was defined only for immersive tropical curves)
 to take the study in 
 Subsections \ref{subsec:general_genus_one} and  \ref{subsec:4-vloop} into account.

Let $(\Gamma, h)$ be a tropical curve of genus one
 such that the direction vectors of the edges span $N_{\Bbb R}$.
Also assume that it can be deformed into an immersive tropical curve.
Let $\bar B$ be the subspace of $N_{\Bbb R}$ spanned by
 the directions of the edges at least one of whose ends is
 on the loop, as in Theorem \ref{thm:highvalloop}.

For any vector $\mathfrak w^{\vee}\in (\bar B)^{\perp}\cap N^{\vee}$, 
 let $\bar H_{\mathfrak w^{\vee}}$ be the associated hyperplane in 
 $N_{\Bbb R}$ and let $H_{\mathfrak w^{\vee}}$ be the 
 affine hyperplane parallel to $\bar H_{\mathfrak w^{\vee}}$
 containing the loop part of $h(\Gamma)$.
We write by
\[
\Gamma_{\mathfrak w^{\vee}}\subset
 h(\Gamma)\cap H_{\mathfrak w^{\vee}}
\]
 the connected component containing the loop part of $h(\Gamma)$.
Let 
\[
p_1^{\mathfrak w^{\vee}}, \dots, p_j^{\mathfrak w^{\vee}}
\]
 be the vertices of $\Gamma_{\mathfrak w^{\vee}}$ with the 
 following properties.
\begin{itemize}
\item An edge emanates from 
 $p_i^{\mathfrak w^{\vee}}$ which is not contained
 in $H_{\mathfrak w^{\vee}}$.
\item There are no other vertices with this property lying on the 
 unique shortest path from $p_i^{\mathfrak w^{\vee}}$ to the loop.
\end{itemize}
By definition of $\bar B$, 
 the vertex $p_i^{\mathfrak w^{\vee}}$
 is not contained in the loop.

Let $\mathcal P_i^{\mathfrak w^{\vee}}$ be the unique shortest path 
 on $\Gamma_{\mathfrak w^{\vee}}$ from 
 $p_i^{\mathfrak w^{\vee}}$ to the loop.
Recall that we defined the length 
 $\ell_{(\Gamma, h)}(\mathcal P_i^{\mathfrak w^{\vee}})$ of such paths
 in Definition \ref{def:path length}.
Also recall that given a tropical curve $(\Gamma, h)$ of genus one, 
 we constructed a parallelotope $P$ and a union of hyperplanes
 $G(H)$ in Subsection \ref{subsec:loopexist}.
Note that when there are no higher valent vertices on the loop of 
 $h(\Gamma)$, the parallelotope $P$ is the empty set and
 the first half in the condition in the following definition 
 is also empty.

\begin{defn}\label{well-spaced2}
The curve $(\Gamma, h)$ is said to be \emph{well-spaced} 
 if the complement $P\setminus G(H)$ contains an integral point, 
 and if moreover 
 one of the following conditions is 
 satisfied for any vector 
 $\mathfrak w^{\vee}\in (\bar B)^{\perp}\cap N^{\vee}$.
\begin{enumerate}
\item The set $\{\ell_{(\Gamma, h)}(\mathcal P^{\mathfrak w^{\vee}}_1),
  \dots, \ell_{(\Gamma, h)}(\mathcal P^{\mathfrak w^{\vee}}_j)\}$ 
  contains at least two minimum.
\item The set
 $\{\ell_{(\Gamma, h)}(\mathcal P^{\mathfrak w^{\vee}}_1),
  \dots, \ell_{(\Gamma, h)}(\mathcal P^{\mathfrak w^{\vee}}_j)\}$
  contains only one minimum.
In this case, 
 let $p_i^{\mathfrak w^{\vee}}$ be the vertex of 
 $\Gamma_{\mathfrak w^{\vee}}$
 at which 
 $\ell_{(\Gamma, h)}(\mathcal P^{\mathfrak w^{\vee}}_i)$ takes the
 minimum among the numbers
 $\{\ell_{(\Gamma, h)}(\mathcal P^{\mathfrak w^{\vee}}_1),
  \dots, \ell_{(\Gamma, h)}(\mathcal P^{\mathfrak w^{\vee}}_j)\}$.
Then there are at least three edges emanating from 
 $p_i^{\mathfrak w^{\vee}}$
 which are not contained in $H_{\mathfrak w^{\vee}}$.
\end{enumerate}
In fact, this definition can be extended to curves of higher genus whose
 support of superabundancy has genus one (see Definition \ref{def:abundancysupport})
 in a straightforward manner.
In this case, the shortest path from a vertex to the loop which supports
 the representatives of the dual obstruction space $H$ may not be unique.
In such a case, each of those shortest paths contributes to the set 
 $\{\ell_{(\Gamma, h)}(\mathcal P^{\mathfrak w^{\vee}}_1),
 \dots, \ell_{(\Gamma, h)}(\mathcal P^{\mathfrak w^{\vee}}_j)\}$
 of lengths.
\end{defn}
\begin{rem}\label{rem:assumption}
\begin{enumerate}
\item 
In this definition, we do not assume that the map
 $h$ satisfies Assumption A.
However, we assume that $h$ can be deformed into an immersion.
\item The reason we assumed Assumption A 
 in previous sections is that 
 in some cases, without Assumption A, it is difficult to give 
 a correspondence between tropical curves and degenerate
 (usually, pre-log) curves.
In fact, as \cite[Example 19]{N1} shows, 
 in general the correspondence between tropical curves
 and degenerate algebraic curves breaks down when a loop is
 contracted by $h$.
However, when we assume the condition that the tropical curve
 can be deformed into an immersive curve 
 (which is (iii) of Assumption A), it seems that
 many of the results in this paper can be extended.
In fact, when the curve is genus one, any tropical curve which
 corresponds to some classical curve has this property
 (see Proposition \ref{prop:nonassumption} below).
\item 
As in Example \ref{ex:high} below, 
 when we consider higher genus curves in general, we need to discard (ii) of
 Assumption A, and it is important to find a reasonable assumption 
 to work with.
In these cases, assuming only (iii) of Assumption A might 
 serve as a nice set up.
\end{enumerate}
\end{rem}

Due to the calculation 
 in Subsections \ref{subsec:general_genus_one} and \ref{subsec:4-vloop},
 the next theorem follows from a straightforward extension of the 
 proof of Theorem \ref{thm:immersive}.
\begin{thm}\label{thm:general}
Let $(\Gamma, h)$ be a tropical curve of genus one
 satisfying Assumption A.
Then $(\Gamma, h)$ is smoothable if and only if it is well-spaced
 in the sense of Definition \ref{well-spaced2}.
\end{thm}
\proof
The smoothability of $(\Gamma, h)$ means two things.
Namely, the existence of a pre-log curve of type $(\Gamma, h)$ and the 
 vanishing of the obstruction to deform it.
Since the proof is basically the same as that of Theorem \ref{thm:immersive}, 
 we only give a sketch.
By Theorem \ref{thm:highvalloop}, 
 we only need to consider obstructions of directions contained in 
 $(\bar B)^{\perp}$.

First assume that $(\Gamma, h)$ is smoothable.
When there are higher valent vertices on the loop and the 
 set $P\setminus G(H)$ contains no integral points, then 
 there are no pre-log curve of type $(\Gamma, h)$, which contradicts to 
 the smoothability.
Therefore, assume the condition about the integral points of $P\setminus G(H)$
 is satisfied.
In this case, if $(\Gamma, h)$ is not well-spaced, then there is some 
 $\mathfrak w^{\vee}\in (\bar B)^{\perp}\cap N^{\vee}$ such that for
 the associated affine hyperplane $H_{\mathfrak w^{\vee}}$, 
 there is a vertex $p^{\mathfrak w^{\vee}}$ of $\Gamma_{\mathfrak w^{\vee}}$ 
 with the following properties.
\begin{itemize}
\item An edge emanates from 
 $p^{\mathfrak w^{\vee}}$ which is not contained
 in $H_{\mathfrak w^{\vee}}$.
\item The lengths of the paths 
 in the sense of Definition \ref{def:path length}
 from the other vertices of 
 $\Gamma_{\mathfrak w^{\vee}}$ with this property to the loop
 are strictly larger than that from the vertex
 $p^{\mathfrak w^{\vee}}$
 to the loop.
\item There are only two edges (counted with multiplicity) emanating
 from $p^{\mathfrak w^{\vee}}$ which are not contained in  
 $H_{\mathfrak w^{\vee}}$.
\end{itemize}

By Proposition \ref{prop:higher}, the leading term of the 
 pre-obstruction of the direction $\mathfrak w^{\vee}$ 
 contributed from the vertex $p^{\mathfrak w^{\vee}}$
 cannot be cancelled, and the contributions to the pre-obstruction 
 of this direction from the other vertices are of higher order
 with respect to the exponent of $t$, the parameter of the smoothing.
Thus, the obstruction of the direction $\mathfrak w^{\vee}$
 cannot vanish and the pre-log curve cannot be smoothed.
This proves that if $(\Gamma, h)$ is smoothable, 
 then it is well-spaced.
 
Conversely, assume $(\Gamma, h)$ is well-spaced.
Then we can choose a basis 
 $\{\mathfrak w_1^{\vee}, \dots, 
 \mathfrak w_r^{\vee}\}$ of $(\bar B)^{\perp}\cap N^{\vee}$ 
 as in the proof of Theorem \ref{thm:immersive}.
Let $\varphi_0\colon C_0\to X_0$ be a pre-log curve of type
 $(\Gamma, h)$ and assume that
 we have a $k$-th order lift 
 $\varphi_k\colon C_k\to \mathfrak X$ of $\varphi_0$.
Then by the proofs of Theorem \ref{thm:immersive}
 and Proposition \ref{prop:higher}, 
 perturbing the coefficients of the defining equations of 
 the pull-back of the components of $\varphi_k(C_k)$ corresponding to 
 appropriate vertices $p_j^{\mathfrak w_i^{\vee}}$, 
 we can cancel the obstructions $o(\mathfrak w_i^{\vee};\varphi_k)$
 simultaneously for all $i$.
Then all the obstructions vanish due to the linearity of the
 obstruction (Lemma \ref{lem:linearity}).
Thus, $\varphi_0$ can be smoothed up to any order.\qed

\begin{rem}
	In her thesis \cite{T}, Torchiani introduced a nice enumerative
	invariant for tropical curves of genus one with a desirable invariance 
	property, by constructing and studying
	combinatorial moduli spaces of such
	tropical curves. 
	This invariant is as nice as one can expect for the counting of tropical curves. 
	In particular, the invariant does not depend on the configurations of
	constraint conditions.
	However, it does not in general give the number of corresponding algebraic
	curves.
	Namely, in her definition, all tropical curves considered in Theorem
	\ref{thm:highvalloop} contribute to the counting number positively.
	However, corresponding algebraic curves need not always exist
	due to Theorem \ref{thm:119}.
	This implies that the counting of classical curves is subtler than the combinatorial 
	count in that one needs to implement the information of the integral points
	of $P\setminus G(H)$ into the study of the moduli spaces.
	We will discuss this problem elsewhere.  
\end{rem}
\begin{rem}\label{rem:mainext}
A part of Theorem \ref{thm:general} can be extended to the case where 
 a tropical curve has higher genus but its support of superabundancy is
 of genus one (see Definition \ref{def:abundancysupport}).
It is easy to see that the well-spacedness is necessary for the smoothability of
 these curves.
However, in general it is not sufficient.
\end{rem}

\subsection{Pure dimensionality of the parameter space of 
well-spaced tropical curves of a given degree}\label{subsec:resol}
Let $\varphi\colon C\to X$ be a 
 torically transverse stable map from an elliptic curve
 $C$ to a toric variety $X$.
The intersection of the image $\varphi(C)$ with the toric divisors
 of $X$ 
 determines the degree $\Delta$ (see Definition \ref{def:degree}).
It is easy to see that the expected dimension of the moduli space
 of stable maps containing $\varphi$ with the same tangency
 to the toric divisors as $\varphi$ is the same as the 
 (set theoretical) cardinality of the inverse image $\varphi^{-1}(\mathcal D)$, 
 where $\mathcal D$ is the union of the toric divisors of $X$.

As in \cite[Proposition 6.3]{NS}, 
 given a toric degeneration 
 $\mathfrak X'$ of $X$ over $\Bbb C$, 
 after a base change and toric blow-ups with centers in $X_0'$
 (we write by $\mathfrak X$ the modification of $\mathfrak X'$
 obtained 
 in this way), 
 the map $\varphi$ extends to 
 a family of stable maps $\Phi\colon \mathcal C\to\mathfrak X$ over
 $\Bbb C$ such that its restriction to the central fiber is a pre-log curve.
Then the central fiber $C_0$ of $\mathcal C$ determines an 
 abstract graph $\Gamma_{C_0}$ 
 as the dual intersection graph 
 (not necessarily trivalent and it might be topologically genus zero
 when some component of $C_0$ has arithmetic genus one)
 and the map $\varphi_0\colon C_0\to X_0$ 
 determines (up to parallel transforms)
 an imbedding $h\colon \Gamma_{C_0}\to N_{\Bbb R}$
 of degree $\Delta$
  which gives $\Gamma_{C_0}$ the structure of a tropical curve.

Now we introduce the following combinatorial operations
 on tropical curves
 to simplify the proofs of the results below.
In this subsection, 
 the tropical curves need not be defined over the integers, 
 and the lengths of the edges are defined in the extended sense
 as remarked at the last of Definition \ref{def:path length}.
 
Let $\Gamma_1$ be an embedded 
 weighted graph in $\Bbb R^n$ which has only
 one vertex $v$.
In particular, all the edges are half lines.
Moreover, assume that the slopes of all the edges are rational and 
 they satisfy the balancing condition at
 the vertex $v$.
Split the set of the edges of $\Gamma_1$ into two disjoint subsets 
 $I$ and $J$ so that the cardinalities of $I$ and $J$ are
 both larger than one.
\begin{defn}
We say that the embedded weighted graph $\Gamma_1'$ in $\Bbb R^n$ is
 obtained by \emph{splitting $J$ from $v$ by length $r$}, 
 if the following properties hold.
\begin{itemize}
\item The graph $\Gamma_1'$ has two vertices $v_1$ and $v_2$,
  and unbounded edges labelled by $I$ and $J$, which have the same directions as those 
  of $\Gamma_1$. 
The vertex $v_1$ is placed at the same place as $v$.
\item From the vertex $v_1$, the unbounded edges labelled by $I$ emanate.
Also, there is a unique bounded edge $E$ of length $r$
 connecting $v_1$ and $v_2$
 so that the union of $E$ and the edges labelled by $I$ satisfy the
 balancing condition at $v_1$.
This condition determines the weight, the length 
 and the direction of $E$ uniquely.
\item From the vertex $v_2$, the unbounded edges 
 labelled by $J$ and the bounded edge $E$ emanate.
\end{itemize}
\end{defn}
\noindent
Then the graph $\Gamma_1'$ satisfies the balancing condition at 
 $v_2$ automatically.

The following proposition shows that Assumption A in fact gives 
 no restriction if we are interested in tropical curves of genus one
 which correspond to classical curves.
\begin{prop}\label{prop:nonassumption}
The map $(\Gamma_{C_0}, h)$
 can be deformed into an immersive 3-valent
 tropical curve $(\Gamma_{C_0}', h')$.
\end{prop}
\proof
If the graph $\Gamma_{C_0}$ is topologically genus one, 
 then it is easy to see that we can deform higher valent vertices of 
 $\Gamma_{C_0}$ into 3-valent trees
 in a way that the map $h$ can also be deformed into an immersion.
Namely, if the vertex is not contained in the loop, then deform it 
 arbitrarily by repeating the above splitting
 until all the vertices become 3-valent.
If the vertex is on the loop, then there are edges
 $E_1, E_2$ emanating from it which are parts of the loop.
If $I$ is the set of the other edges emanating from that vertex, 
 split $I$ from $v$ by arbitrary positive length.
This reduces the problem to the former case.

Let us consider the case when the graph 
 $\Gamma_{C_0}$ is topologically a tree.
There is a vertex $v$ corresponding to the component of $C_0$
 of arithmetic genus one.
In view of the preceding paragraph, it suffices to prove the following.
\begin{claim}
Let $D$ be an irreducible curve of arithmetic genus one
 and $Y$ be a toric variety.
Then any torically transverse embedding $\psi\colon D\to Y$
 can be degenerated into the family 
 $\Phi\colon\mathcal D\to\mathfrak Y$ 
 as in the beginning of this subsection 
 in a way that the resulting 
 graph $\Gamma_{D_0}$ is topologically genus one.
\end{claim}
\proof
Note that to the map $\psi\colon D\to Y$ corresponds 
 a tropical curve $\Gamma_D\subset \Bbb R^k$ 
 which has only one vertex.
We prove the above claim
 by induction with respect to the number of 
 unbounded edges of $\Gamma_D$.
Since $\psi$ is an embedding, the image of $\psi$ is not
 contained in the closure of a one dimensional torus orbit of $Y$.

To begin the induction, assume that $\Gamma_D$ has three 
 unbounded edges.
In this case, the image $\psi(D)$ is necessarily contained in 
 the closure of a two dimensional torus orbit of $Y$, 
 and we can assume that $Y$ is a two dimensional toric variety
 from the first.

As we noted above, the expected dimension of the moduli space
 of stable maps containing $\psi$ is three. 
Therefore, $\psi$ has
 at least three dimensional freedom to deform.
Then as in \cite{NS}, we can impose a family of incidence conditions
 consisting  
 of three general sections of
 $Y\times(\Bbb C\setminus\{0\})\to \Bbb C\setminus\{0\}$
 and extend it together with $\psi$ to a family of stable maps over $\Bbb C$
 satisfying the incidence conditions (see \cite[Proposition 6.3]{NS}).
Then the dual intersection graph of the central fiber gives rise to 
 a tropical curve on $\Bbb R^2$
 of the given degree incident to general three points.
However, it is impossible to meet the incidence conditions unless the
 original vertex is resolved into a loop.
 
Now assume that we have proved the claim when the 
 number of the unbounded edges of $\Gamma_D$
 is at most $k(\geq 3)$.
Let us take a map $\psi\colon D\to Y$ such 
 that the corresponding graph
 $\Gamma_{D}$ has $k+1$ unbounded edges.
Let $A$ be the minimal affine subspace containing $\Gamma_D$
 and write $d = \dim A$.
Then $\psi(D)$ is contained in a $d$ dimensional torus orbit
 of $Y$, and we can assume $Y$ is $d$ dimensional from the first.

As in the above argument, the number of freedom to deform 
 $\psi$ is at least $k+1$, but clearly it is strictly bigger than $d$.
Then by taking appropriate incidence conditions of codimension $k+1$
 (now 
 consisting not only of points
 but also of suitable torus orbits, see \cite{NS})
 and degenerations
 of $Y$ and $\psi$ as above, the central fiber of the degeneration 
 gives rise to a tropical curve of the given degree which matches
 incidence conditions of codimension $k+1$ 
 given by affine linear subspaces.
However, since the tropical curve corresponding to $\psi$, 
 which has only one vertex, has only $d$ dimensional freedom to deform
 (that is, the freedom of parallel transport), 
 it is impossible to meet the incidence conditions unless the vertex
 resolves into a loop or a tree with more than one vertices.
In the former case the assertion is proved, and in the latter case
 it is clear that all the new vertices have valence smaller than 
 $k+1$.
Thus, by induction hypothesis the claim is proved as well. \qed\\

Now let us complete the proof of the proposition.
When the graph $\Gamma_{C_0}$ is a tree,
 then by the claim we can deform 
 the vertex of it corresponding to the component of $C_0$ 
 which has arithmetic genus one into a graph with a loop.
Since the complement of the vertex is a union of trees, 
 we can extend this to a global deformation
 of $(\Gamma_{C_0}, h)$ preserving the balancing condition. 
Then by the first paragraph of the proof, the resulting graph can be 
 deformed into an immersive trivalent tropical curve. \qed
\begin{cor}\label{cor:maxdeg}
Any embedded curve of arithmetic genus one in a toric variety 
can be maximally degenerated.
That is, if $\varphi\colon C\to X$ is an embedding of a curve of arithmetic genus one into 
 a toric variety, then there is a family of stable maps 
 $\Phi\colon \mathcal C\to \mathfrak X$ over $\Bbb C$ with the following properties:
\begin{itemize}
\item $\mathcal C$ is a family of pre-stable curves and
 $\mathfrak X$ is a toric degeneration of $X$.
\item The restriction of $\Phi$ to $1\in\Bbb C$ is isomorphic to $\varphi$.
\item Let $\varphi_0\colon C_0\to X_0$ be the restriction of $\Phi$
 to $0\in \Bbb C$.
Then $\varphi_0$ is a maximally degenerate curve (see Remark \ref{rem:type0}).
In particular, any irreducible component of $C_0$ is a nonsingular rational curve.  \qed
\end{itemize}
\end{cor}

\begin{rem}\label{rem:contract}
Let the map $h\colon \Gamma\to N_{\Bbb R}$ be a tropical curve of genus one
 which may not satisfy Assumption A, but assume that it can be deformed into 
 an immersive tropical curve (the latter condition is necessary for the smoothability
 due to Proposition \ref{prop:nonassumption}).
If the map $h$ contracts the loop part, we cannot apply Theorem \ref{thm:general}
 to it.
However, in view of the arguments in this paper,
 it is rather straightforward to see that well-spacedness
 is a necessary condition for the smoothability of a pre-log curve
 corresponding to such $(\Gamma, h)$ 
 (in this case, there is a component of the domain $C_0$ 
 of a pre-log curve which has 
 arithmetic genus one corresponding to the vertex which deforms
 into a loop).
In particular, since the curve $\Phi|_{C_0}\colon C_0\to X_0$ 
 at the beginning of this subsection is
 a pre-log curve of type $(\Gamma_{C_0}, h)$, the tropical curve
 $(\Gamma_{C_0}, h)$ is must be well-spaced.
Note that in this case, the loop part of the image $h(\Gamma_{C_0})$
 is the one point set consisting of th e vertex 
 which is the image of the loop of $\Gamma_{C_0}$.

In fact, it is quite plausible that the well-spacedness condition is also 
 the sufficient condition in this case too, but to show that we 
 need to calculate the dual obstruction
 (the main point is showing the analogue of
 Theorem \ref{thm:highvalloop} for the component of $C_0$
 with arithmetic
 genus one mentioned above).
It will not be very difficult to do it, but we omit it in this paper, 
 in view of Proposition \ref{prop:nonassumption}.
When we consider curves of higher genus, the importance of such a calculation 
 increases, as Example \ref{ex:high} below shows.
We will study this issue elsewhere.
\end{rem}

In \cite[Theorem 3.2.10]{T}, Torchiani proved the following.
\begin{thm}\label{prop:deformdeg}
Let $(\Gamma, h)$ be a well-spaced tropical curve of genus one
 which can be deformed into an immersive tropical curve. 
We assume $(\Gamma, h)$ is not contained in an affine hyperplane
 of $N_{\Bbb R}$.
Let $e$ be the number of unbounded edges.
Then $(\Gamma, h)$ has $e$ dimensional freedom to deform 
 preserving the well-spacedness condition.
\end{thm}
In fact, she proved that the space of these tropical curves has a natural 
 structure of a tropical variety.

Combining this with Proposition \ref{prop:nonassumption}, 
 we have the following assertion for classical elliptic curves.
\begin{cor}
Let $C$ be a 
 smooth elliptic curve and let $X$
 be a toric variety.
Let $\varphi\colon C\to X$ be a torically transverse 
 algebraic curve which is not contained in 
 the closure of an orbit of a proper subtorus of the torus
 acting on $X$.
Note that the map $\varphi$ determines the degree map
 $\Delta\colon N\setminus\{0\}\to \Bbb N$. 
Then the moduli space of flat deformations of $\varphi$
 preserving the torically transversality and degree $\Delta$
 is of pure dimension $|\Delta|$. 
\end{cor}
\proof
By Corollary \ref{cor:maxdeg}, 
 the map $\varphi$ can be extended to a degenerating family
 whose central fiber $\varphi_0$ is a maximally degenerate pre-log curve.
The dual intersection graph of this pre-log curve gives a tropical curve, 
 which is an immersive 3-valent tropical curve.
 
By the assumption, this tropical curve, 
 which we write by $(\Gamma_0, h_0)$,
 is not contained in a
 proper affine subspace of $N_{\Bbb R}$.
Also note that $(\Gamma_0, h_0)$ is well-spaced since it is smoothable by construction.
Then by Theorem \ref{prop:deformdeg}, the dimension of deformations
 of $(\Gamma_0, h_0)$
 preserving the well-spacedness is equal to
 $|\Delta|$.

Now if $\varphi$ has $k$ dimensional freedom of deformations, 
 the map $\varphi_0$ also has at least $k$
 dimensional freedom of deformations.
Since these deformations give deformations of $(\Gamma_0, h_0)$,
 $k$ is bounded from above by $|\Delta|$.
On the other hand, since $|\Delta|$ is the expected dimension 
 of the space of deformations of $\varphi$ preserving the degree
 $\Delta$, $\varphi$ has at least $|\Delta|$ dimensional freedom
 to deform. 
This proves the corollary.\qed

\subsection{Examples}
In this subsection we give several examples which will be of interest
 and related to the general theory developed in this paper.
\begin{example}\label{ex:high}
Here we consider an example of a tropical curve of higher genus,
 in which a situation where 
\begin{itemize}
\item Assumption A does not hold and
\item highly degenerate tropical curves inevitably appear 
 as a general member of the set of smoothable tropical curves
\end{itemize}
 happens.
Namely, consider the following Figure \ref{fig:highdeg}.

\begin{figure}[h]
\includegraphics{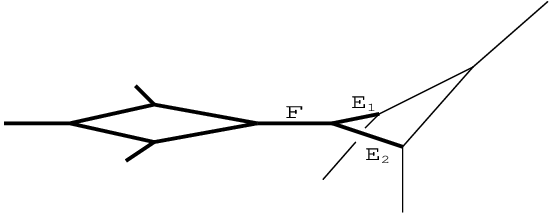}
\caption{Segments drawn by bold lines are contained in the affine subplane spanned
 by the directions of the edges in the loop.}\label{fig:highdeg}
\end{figure}

It is an immersed superabundant tropical curve of genus two in $\Bbb R^3$, 
 and $\dim H$ is equal to one, where $H$ is the dual space of obstructions.
The obstruction is localized at the loop on the left, drawn by bold lines.
The loop on the right is, 
 seen as a tropical curve of genus one, regular.
Assume that the integral length of the edge $E_2$ is twice as
 long as
 that of $E_1$ and the weights of these edges are the same.
Since the loop containing $E_i$ is regular 
 (and it has only four edges), 
 one sees that the ratio of the integral lengths of $E_2$ and $E_1$
 does not change when we deform the tropical curve.

Consider a condition for the smoothability of this tropical curve.
Although it is not of genus one, we can calculate the obstruction as in the calculation 
 in the previous subsections, and it is easy to see that,
 for the vanishing of the obstruction, 
 the integral lengths of the edges $E_1$
 and $E_2$ should be equal.
However, according to the above observation, it only happens when 
 both $E_1$ and $E_2$ are contracted.
Moreover, it is easy to see that
 when the edge $E_1$ or $E_2$ is contracted, the whole loop must be contracted.

By the way, according to the description of the dual obstruction space associated
 to higher valent vertices (\cite[Theorem 60]{N1}), 
 it follows that the case when the edge $F$ is
 contracted (and either of the loops is not contracted) is
 also smoothable (provided a pre-log curve of the given type exists).

Thus, the space of smoothable tropical curves of this type (if not empty)
 consists of those which contract the loop or the edge $F$, and both have the
 same expected dimension to deform.
This is in contrast to the genus one case, where
 although the set of smoothable tropical curves in general contains
 those which contract the loop, such locus has positive codimension.
Therefore, while Assumption A does not give major restriction in the case of genus one,
 for the general study of higher genus curves, 
 we need to weaken Assumption A in an essential way, 
 see Remark \ref{rem:assumption}.  \qed
\end{example}

\begin{example}
In \cite[Example 48]{N1}, we considered the following tropical curve 
 (see Figure \ref{fig:0G}).
\begin{figure}[h]
	\includegraphics{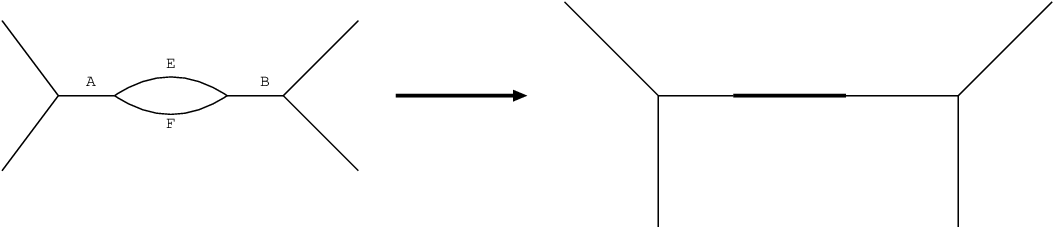}
	\caption{The abstract graph (the picture on the left) is mapped into $\Bbb R^2$.
		The edges $E$ and $F$ have the same image
		(the bold line in the picture on the right).}\label{fig:0G}
\end{figure}
\noindent
Plane tropical curves were studied by Mikhalkin \cite{M} in great detail.
There it was shown that any conventional (that is, immersive and
 without multiple edges) plane tropical curve is regular
 and smoothable.
On the other hand, if we only assume Assumption A, 
 then even 3-valent immersive (in the sense of Definition \ref{immersive}) 
 plane tropical curves
 can be superabundant and non-smoothable.

In Figure \ref{fig:0G}, all the edges except $E$ and $F$ have weight two, 
while the edges $E$ and $F$ have weight one.
The loop part is the union of $E$ and $F$, and the space $H$ is one dimensional.
By Theorem \ref{thm:immersive} or \ref{thm:general}, this tropical curve is smoothable 
if and only if the lengths of the images of the edges $A$ and $B$ are the same.

\begin{figure}[h]
	\includegraphics{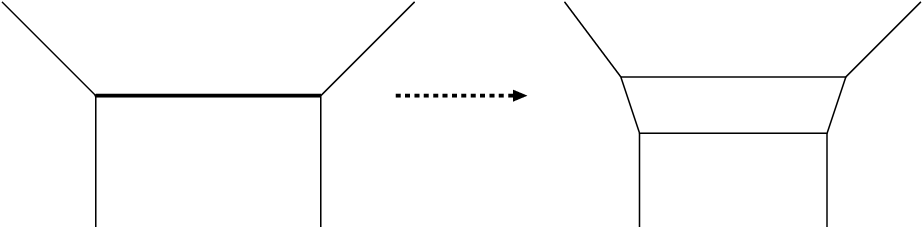}
	\caption{Tropical curves in Figure \ref{fig:0G}
		can be deformed into a tropical curve with two 4-valent vertices
		(the picture on the left).
	It can then further deformed into embedded tropical curves
  (the picture on the right).
 These tropical curves are smoothable.}\label{fig:0I}
\end{figure}

On the other hand, tropical curves in Figure \ref{fig:0G} can be deformed into 
 tropical curves in Figure \ref{fig:0I}.
These tropical curves are all smoothable.
Note that the tropical curves on the right of Figure \ref{fig:0G} fulfilling the 
 well-spacedness condition, and the tropical curves on the right of Figure \ref{fig:0I}
 both have four dimensional freedom to deform.
Since they are smoothable, corresponding algebraic curves also 
 has four dimensional freedom.
These give two components of the degeneration of the
 moduli space of curves of the given degree
 in the toric variety $\Bbb P(1, 1, 2)$.

As in \cite{M, NS}, we can impose incidence conditions to these tropical curves
 and to corresponding algebraic curves as well, 
 which in this case will be four points in generic position.
Those curves satisfying incidence conditions give Gromov-Witten like counting 
 invariants.
However, only the curves corresponding to the picture
 on the right of Figure \ref{fig:0I} contribute
 to the invariants.

\end{example}

\begin{example}
Here we give several examples in the case of cubic
 curves of genus one in $\Bbb P^3$.
The is one of the most typical superabundant cases, and
 various cases appeared in the main text of this paper
 can be realized in this example.

First, consider the following figure (Figure \ref{fig:cubic}).
There are three unbounded edges for each of the following directions
\[
(-1, 0, 0),\;\; (0, -1, 0),\;\; (1, 1, 1),\;\; (0, 0, -1).
\]
All of these edges have weight one.
The three black dots in the figure
 means the unbounded edges of direction $(0, 0, -1)$.
Cubic curves of genus one
 in $\Bbb P^3$ are known to be contained in some projective plane,
and it is true also in the tropical case.
The dotted lines represent the image of the one
dimensional skeleton of the
tropical hyperplane containing the tropical
cubic curve of genus one under the projection to the 
plane.

\begin{figure}[h]
\includegraphics{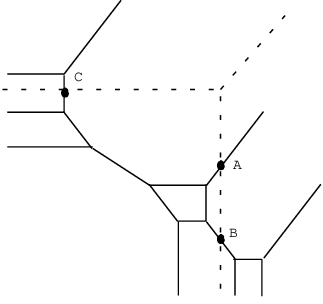}
\caption{Embedded tropical cubic curve in $\Bbb R^3$.}\label{fig:cubic}
\end{figure}

In this figure, the vertices $A, B, C$ are the 1-valent vertices of 
 $\overline\Gamma'$, 
 and the vertices $A$ and $B$ assure the well-spacedness condition.

When we slide the curve on the plane, the picture becomes as follows
 (Figure \ref{fig:cubicweight2}).
Note that we do not deform $\Gamma$, but only the map $h$.

\begin{figure}[h]
\includegraphics{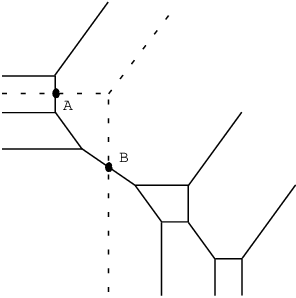}
\caption{Non-immersive
 tropical cubic curve.
 The dot $B$ corresponds to a multiple edge with weight $(1, 1)$.}
  \label{fig:cubicweight2}
\end{figure}
In this case, the unbounded edge of the direction $(0, 0, -1)$ 
 at the vertex $B$ has weight $(1, 1)$,
 and this is the special case of 
 Case (b-2) of Subsection \ref{subsec:general_genus_one}, 
 where there are 
 different edges of $\Gamma$ with 
 the same image.
Therefore, this satisfies the extended well-spacedness 
 condition of Definition \ref{well-spaced2}.

Next, consider the cases
 when some of the horizontal unbounded edges of $\Gamma$ are merged,
 as in the following picture (Figure \ref{fig:cubicmove}).
Note that the unbounded edges are merged \emph{in $\Gamma$}.
That is, the image of the merged edge has weight two, not (1, 1)
 (the latter case is not generic in the space of smoothable tropical curves).

\begin{figure}[h]
\includegraphics{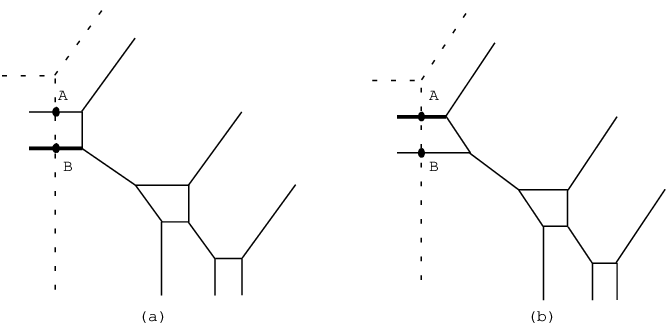}
\caption{}\label{fig:cubicmove}
\end{figure}

Forgetting the vertical direction and seen
 as curves on a plane, these 
 correspond to cubic curves with the condition that 
 they have one intersection point of multiplicity two with a toric divisor.
The edges drawn by bold lines in Figure \ref{fig:cubicmove}
 have weight two, and correspond to these intersection points.

On the other hand, the vertical unbounded edges emanating from the 
 vertices on the bold lines 
 (the dot $B$ of (a) of Figure \ref{fig:cubicmove}
 and the dot $A$ of (b) of Figure \ref{fig:cubicmove})
 also have total additive weight two (Definition \ref{def:trop}).

In the case of (a), the integral distance from $B$ to 
 the loop is shorter than that from $A$,
 so for the well-spacedness condition, it is necessary that the 
 vertical edge from $B$ has weight $(1, 1)$.
This is Case 2 of Subsection \ref{subsec:general_genus_one}.

In the case of (b), the integral distance from $A$ or $B$ to the loop is the same.
Since the length of the edge with weight two should be halved, the leading contribution
 to the obstruction comes from the vertex $A$.
There are two cases.
\begin{enumerate}[(i)]
\item The vertical edge from $A$ has weight (1, 1).
\item The vertical edge from $A$ has weight two. 
In this case, the tropical curve is immersed.
\end{enumerate}
The former case is smoothable (again
 Case 2 of Subsection \ref{subsec:general_genus_one}).
The latter case corresponds to a genus one cubic curve which has two 
 intersection points of multiplicity two with toric divisors, but not smoothable.
\begin{rem}
As this example shows, not every cubic tropical curve contained in 
 a tropical hyperplane is smoothable.
In fact,  this can happen
 even when we restrict our attention to immersive tropical curves.
See Figure \ref{fig:cubicnon}. 
\end{rem}

\begin{figure}[h]
\includegraphics{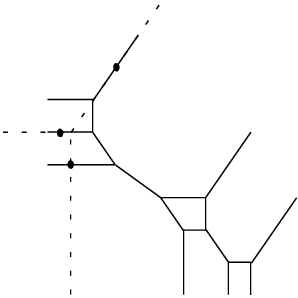}
\caption{An immersive tropical cubic curve on a tropical hyperplane
 which is not smoothable.}\label{fig:cubicnon}
\end{figure}
\end{example}

\begin{example}
Finally, we give an example that clearly 
 shows why we have to divide the integral length of 
 an edge by its weight.

\begin{figure}[h]
\includegraphics{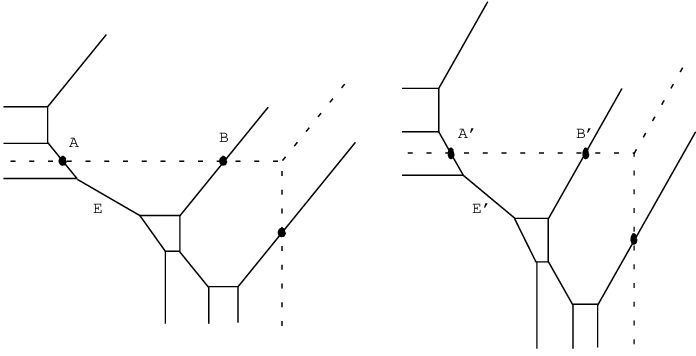}
\caption{}\label{fig:weight}
\end{figure}

In Figure \ref{fig:weight}, 
 the left tropical curve corresponds a standard 
 genus one cubic curve
 in $\Bbb P^3$.
The right tropical curve is obtained from the left 
 by the linear transformation by 
 the matrix
 $\begin{pmatrix}
 1&0\\
 0&2
 \end{pmatrix}
$
 on the plane.
As a result, this tropical curve is also smoothable
 (namely, the smoothing is given as the image of 
 genus one cubic curves in $\Bbb P^3$ by the map 
 between toric varieties associated to the above linear map).
Note that the edge $E'$ has the weight two.

In the left tropical curve,
 the integral lengths of the paths from the vertices $A$ and $B$ to the
 loop are the shortest ones.
While, the integral length of the edge $E'$ in the right tropical curve
 is the twice of that of the edge $E$ in the left tropical curve.
Therefore, the integral length of the path
 from the vertex $A'$ to the loop is greater than that of 
 the path from the vertex $B'$ to the loop.
However, when we halve the integral length of $E'$,
 these two are the same, 
 adjusting to the well-spacedness condition.
\end{example}

\end{document}